\documentclass[fleqn,10pt]{wlscirep}
\usepackage[utf8]{inputenc}
\usepackage{algorithm2e}
\usepackage{stix}
\usepackage{multirow}
\usepackage[symbol]{footmisc}
\usepackage{comment}

\usepackage{epsfig}
\usepackage{soul}

\usepackage{tikz}
\usepackage{tkz-fct}
\usetikzlibrary{arrows,shapes,trees}

\RestyleAlgo{ruled}
\usepackage{titletoc}
\SetLabelAlign{Center}{\hfil#1\hfil}
\SetLabelAlign{CenterWithParen}{\hfil(\makebox[1.0em]{#1})\hfil}

\definecolor{irishgreen}{RGB}{0,154,73}

\usepackage[T1]{fontenc}
\title{
Classifying {absence} seizure generation mechanisms:\\ A critical transitions framework\\
\vspace{2mm}
{\small 
New insights from combining brain voltage recordings, models of critical transitions and machine learning
}
}

\author[1,2,*]{Andrew Flynn}
\author[3]{Cian McCafferty}
\author[4,5,6]{Klaus Lehnertz}
\author[7]{Fran\c{c}ois David}
\author[8]{Vincenzo Crunelli}
\author[2,9]{William P. Marnane}
\author[1]{Sebastian Wieczorek}
\affil[1]{School of Mathematical Sciences, University College Cork, T12 XF62 Cork, Ireland.}
\affil[2]{INFANT Research Centre, University College Cork, T12 DC4A Cork, Ireland.}
\affil[3]{Department of Anatomy \& Neuroscience, University College Cork, Cork, Ireland.}
\affil[4]{Department of Epileptology, University of Bonn Medical Centre, 53127 Bonn, Germany.}
\affil[5]{Helmholtz-Institute for Radiation and Nuclear Physics, University of Bonn, 53115 Bonn, Germany.}
\affil[6]{Interdisciplinary Center for Complex Systems, University of Bonn, 53175 Bonn, Germany.}
\affil[7]{Center for Interdisciplinary Research in Biology, Coll{\`e}ge de France, 75005 Paris, France}
\affil[8]{Department of Pharmacology and Neuroscience, Faculty of Medicine, University of Lisbon, Lisbon, Portugal}
\affil[9]{Electrical and Electronic Engineering, University College Cork, T12 YF78 Cork, Ireland.}

\affil[*]{andrew\_flynn@umail.ucc.ie}

\begin{abstract}
Understanding how the brain switches from normal activity to an epileptic seizure is essential for improving seizure therapy, yet the underlying seizure generation mechanisms remain largely unknown. In particular, while seizure onset has been described as a critical transition (CT), there is no consensus on whether (i) bifurcation-induced, (ii) noise-induced, or (iii) bifurcation/noise-induced CTs are responsible. To clarify this, we develop a versatile CT-classification framework that can be applied to seizures in both animals and humans. First, we identify a canonical mathematical model which displays CTs that closely resemble voltage recordings of real seizures and can be of the three types mentioned above. We then identify distinctive properties of each CT-type in the model’s output and use them to train a machine learning CT-type classifier. Finally, we apply the model-trained classifier to voltage recordings from epileptic rodents which consist of thousands of real absence seizures. We find that the largest proportion of analysed seizures are classified as noise-induced CTs. In other words, our results on absence seizures in rodents are in contrast to the conventional view that seizures are predominantly bifurcation-induced, and indicate that different CT mechanisms may dominate different seizure types.
\\
\textbf{Relevance to Life Sciences: }
We view the onset of epileptic seizures as a CT from a non-seizure to a seizure brain state, where different CT-types reflect distinct neural processes involved in generating seizures. For the first time, we have successfully classified seizure onset as different CT-types by developing a new mathematical framework and applying it to thousands of real seizures in epileptic rodents. Our analysis shows that the largest proportion of analysed seizures are consistent with noise-induced CTs, suggesting that these particular seizures may predominantly be triggered by intermittent disturbances, either exogenous or endogenous, to particularly sensitive groups of neurons. This is in contrast to the prevailing view that seizures are predominantly triggered by bifurcations of the non-seizure state, which would require a gradual change in the conditions that cause the shift to the seizure state. These new insights could have implications for tailored therapies designed to prevent or terminate different types of seizures.
\\
\textbf{Mathematical Content: }
We develop a mathematical framework for classifying CT-types in noisy time series. This framework combines mathematical modelling, bifurcation analysis, time series feature extraction, and machine learning tools. While our motivation was to improve our understanding of seizure onset, the framework is universal and can be applied to other seizure types, other CTs in the brain, as well as to CTs in complex systems beyond the brain.

\end{abstract}

\begin{document}

\flushbottom
\maketitle


\section*{Introduction}


Epileptic seizures are states of high-amplitude, synchronous electrical activity in the brain with defined onset and offset\cite{beniczky2025_LeagueAgainstEpil_def}. Seizures have a wide range of associated symptoms but all constitute a large change in brain state and have a commensurate impact on the life of a person with epilepsy\cite{kaye2025_seizimpact}. Why the brain suddenly changes from normal activity to an epileptic seizure is largely unknown\cite{gonzalez2019_WhySeizHappen}. 
A better understanding of this process may lead to the more optimal use of existing therapeutic interventions and advance the development of new ones that either prevent onset of seizure or accelerate its end.
A sudden and large change in the state of a complex system exactly fits Ashwin~\textit{et al.}'s~\cite{ashwinwieczorek2017CTdef} definition of a critical transition (CT).
The mathematical theory of CTs is perhaps best known for describing and foreseeing sudden and large changes in climate and ecological systems~\cite{lenton2008tippingpolicymaking,lenton2020tippingpositivechange}.  
It has provided a better understanding of key tipping elements in the Earth system and their impacts on our lives, and has also informed global policymaking~\cite{armstrong2022exceeding,lenton2025global}.
Inspired by its impact on the environmental sciences, we apply the theory of CTs to analyse seizure activity in the brain. Specifically, we model the onset of seizure activity as a CT between two different states of the brain, the non-seizure state (NS state) and the seizure state (S state), and develop a framework to classify the CT type.

{
Previous studies focused mainly on identifying early warning signals prior to the seizure onset. They also interpreted these early warning signals as clues to the types of CTs involved in generating seizures.  
} 
Scheffer \textit{et al.}~\cite{scheffer2009early} reported that, in an example of an electroencephalography (EEG) time series with seizure activity  taken from McSharry \textit{et al.}~\cite{mcsharry2003prediction}, the variance of the time series increased prior to  seizure onset. 
From the theory of CTs, this increase in variance is a signature of `critical slowing down' (CSD), a phenomenon that occurs prior to a bifurcation-induced CT, thus indicating that the seizure occurred due to a bifurcation taking place in the brain. 
Meisel and Kuehn~\cite{meiselkuehn2012scaling} later analysed an EEG time series from eight adults with epilepsy and found similar increases in variance before only two out of eight seizures. They also showed that such increases occur prior to a subcritical Hopf bifurcation in  a mathematical model. 
Milanowski and Suffczynski~\cite{milanowski2016seizures} expanded on this approach by analysing hundreds of EEG time series of various seizure types in adult humans. 
They restricted their analysis to CTs due to crossing different Hopf bifurcations and concluded that seizures start without common signatures of such CTs. 
Jirsa \textit{et al.}~\cite{Jirsa_2014_NatureOfSeizDyn} used a different modelling approach, based on noise-free dynamic bifurcations, to argue that seizures occur due to slow passage through a saddle-node bifurcation. 
In more recent times, further opposing reports emerged based on large-scale investigations of EEG time series: Wilkat \textit{et al.}~\cite{wilkat2019noCSDseizure} found no evidence of CSD before seizures, whereas Maturana \textit{et al.}~\cite{maturana2020yesCSDseizure} did. 
From a CT point of view, if a bifurcation does not appear to be responsible for the seizure onset, then another likely candidate is a noise-induced CT, which has long been theorised to play a role in seizure onset~\cite{DaSilva03EpilepsyDynDis}.
However, very little is known about the typical signatures of noise-induced CTs in a time series, nor how these signatures manifest themselves in voltage recordings of seizure activity~\cite{Dakos_2024_ESD}.

{
In this paper, we analyse thousands of seizures from three different brains and approach the question of `how are seizures generated?' from a broader perspective 
and with a different focus. 
Unlike previous studies, our focus is on identifying the type of CT rather than on early warning signals. Therefore, we examine the properties of the voltage recordings both before and shortly after seizure onset.}
Most importantly, we go beyond the realm of bifurcation-induced CTs and hypothesise that the mechanism of seizure generation is a CT from the NS to S state that can be: noise-induced, bifurcation-induced, or a combination of both. 
In other words, we consider these three CT-types as potential seizure generation mechanisms. 
To test this hypothesis, we introduce a new CT-based framework, summarised in five steps in Table~\ref{tab:StepsToClassifyData}, to analyse voltage recordings of real seizure activity. 
This framework consists of an algorithm to detect CTs between NS and S states, and a machine learning classifier to distinguish between the different CT-types. 
This enables us to test whether, at a given moment in time near a CT to the S state, the brain: (i) operates within a multistable or excitable regime whereby noise alone in the form of exogenous or endogenous disturbances can cause it to cross some threshold and trigger a seizure, (ii) slowly drifts towards and then past a bifurcation point that generates the seizure, or (iii) experiences a combination of both (i) and (ii).
Seen in this light, our framework can be used not only to address a fundamental question regarding seizure generation mechanisms, but also to inform future seizure prevention and treatment strategies - once we know the CT types that are responsible, we can relate them to different types of neural activity, and adapt therapeutic strategies accordingly.

The structure of the Results section reflects our five-step framework from Table~\ref{tab:StepsToClassifyData}. 
In the first step, we identify a set of two differential equations (an augmented Bautin bifurcation normal form) as a canonical mathematical model capable of generating the three CT-types listed above. We then show that the parameters of our mathematical model can be tuned so that its output resembles real seizure activity in voltage recordings taken from inside the brain of epileptic rodents. 
More specifically, we use local field potential (LFP) measurements of neural activity in Genetic Absence Epilepsy Rats from Strasbourg (GAERS) that were obtained and annotated by McCafferty \textit{et al.}~\cite{mccafferty2018cortical}. We refer to these LFP measurements as the `voltage recordings' throughout the paper; see part~\hyperref[md:CianSeizureAnnotate]{M1} of the Methods section for further details. 
We choose to work with GAERS because their absence seizures are fully generalised and morphologically distinctive in EEG~\cite{crunelli2020roots_of_epil} compared to other seizures which may be focal or have variable EEG signatures~\cite{zhou2024generalized_vs_focal}. In other words, GAERS seizures are both stable in amplitude and frequent in occurrence\cite{depaulis2006models}, making them ideal for investigating seizure onset dynamics and for constructing our framework. We see this as a proof of concept and an important step towards analysing more complex human seizures.
In the second step, we construct an algorithm based on the concept of a non-ideal relay to detect CTs between the NS and S states in noisy time series. 
We apply this detection algorithm to both the model's output and the voltage recordings. Specifically, we test the algorithm on our mathematical model and, crucially, show that it detects CTs between NS and S states in the actual voltage recordings from different GAERS in excellent agreement with expert annotations.
In the third step, we identify several distinctive properties in time series of the mathematical model's output that are associated with each CT-type. This is an important component of our study. It allows us, in the fourth step, to train a machine-learning CT-type classifier to distinguish between the three CT-types from the NS to S state in a noisy time series, without any prior knowledge of CT-types.
In the fifth and final step, we use the CT-type classifier, which has been trained using the model's output, to identify the CT-types responsible for real seizure onset in voltage recordings from different GAERS.

Our findings demonstrate that: (i)  seizures 
{can be}
generated by different types of CT and (ii) noise-induced CTs appear to be the dominant mechanism responsible for generating the seizures we analysed. 
These results align with recent experiments on seizure termination in GAERS~\cite{mccafferty2025interrupt_seiz},
{are in contrast to}
the conventional view that seizures occur predominantly due to bifurcation-induced CTs, and highlight the need to account for 
{different CT-types}
when developing  seizure therapies.
Furthermore, our framework is versatile in the following ways: (i) it can use other mathematical models of seizure activity, which may be more realistic, to train different CT-type classifiers, (ii) it can classify CT-types responsible for other types of seizures, including those in humans, and (iii) it can also be applied to similar instabilities found in complex systems beyond the brain.

\begin{table}[t]
\centering
\def\arraystretch{1.2}
\begin{tabular}{|c||l|}
\hline
    \multirow{4}{12em}{\centering \textbf{Step 1:
    Build mathematical model of CTs with real seizure characteristics.
    }
    } & \multirow{4}{35em}{Construct a canonical mathematical model that displays CTs from its NS to S state that have characteristics of real seizure activity and can be of three types:\\ (i) bifurcation-induced, (ii) noise-induced, (iii) bifurcation/noise-induced.
    }\\ & \\ & \\ & \\
\hline
    \multirow{4}{12em}{\centering \textbf{Step 2: 
    Design and build a CT detection algorithm.
    }
    } & \multirow{4}{35em}{Create a single CT detection algorithm to  detect CTs between the NS and S states, both in the model's output and the voltage recordings in agreement with expert annotations.}\\ & \\ & \\ & \\
\hline
    \multirow{4}{12em}{\centering \textbf{Step 3: 
    Identify distinctive properties of the model's output near each CT-type.
    }
    } & \multirow{4}{35em}{Use the CT detection algorithm from Step 2 to detect CTs from the NS to S state in the model's output, where the CT-type is known. For CTs of each type that satisfy selection criteria, identify their characteristic time series properties. 
    }\\ & \\ & \\ & \\
\hline
    \multirow{4}{12em}{\centering \textbf{Step 4: Build CT-type classifier.
    }
    } & \multirow{4}{35em}
    {Train a machine learning CT-type classifier using the characteristic time series properties identified for each CT-type in the model's output in Step 3. Optimise the classifier by maximising its ability to correctly classify the CT-types in the unseen model's output.} \\ & \\ & \\ & \\
\hline
    \multirow{5}{12em}{\centering \textbf{Step 5: Classify the CT-type of real seizures. 
    }
    } & \multirow{5}{35em}
    {Use the CT detection algorithm from Step 2 to detect CTs from the NS to S state in the voltage recordings, where the CT-type is not known. For detected CTs that satisfy  classification conditions, calculate the characteristic time series properties identified in Step 3. Use these properties as an input for the model-trained CT-type classifier from Step 4 to classify the CT-type of real seizures.}\\ & \\ & \\ & \\ & \\
\hline
\end{tabular}
\caption{
\label{tab:StepsToClassifyData} 
A summary of our proposed five-step framework for classifying CTs that are responsible for the onset of seizures in the voltage recordings of brain activity as: (i) bifurcation-induced, (ii) noise-induced or (iii) bifurcation/noise-induced. 
}
\end{table}


\section*{Results}


\subsection*{Step 1: A mathematical model with  three CT types and real seizure characteristics}

\begin{figure}[t]
    \centering
    \includegraphics[width=\textwidth]{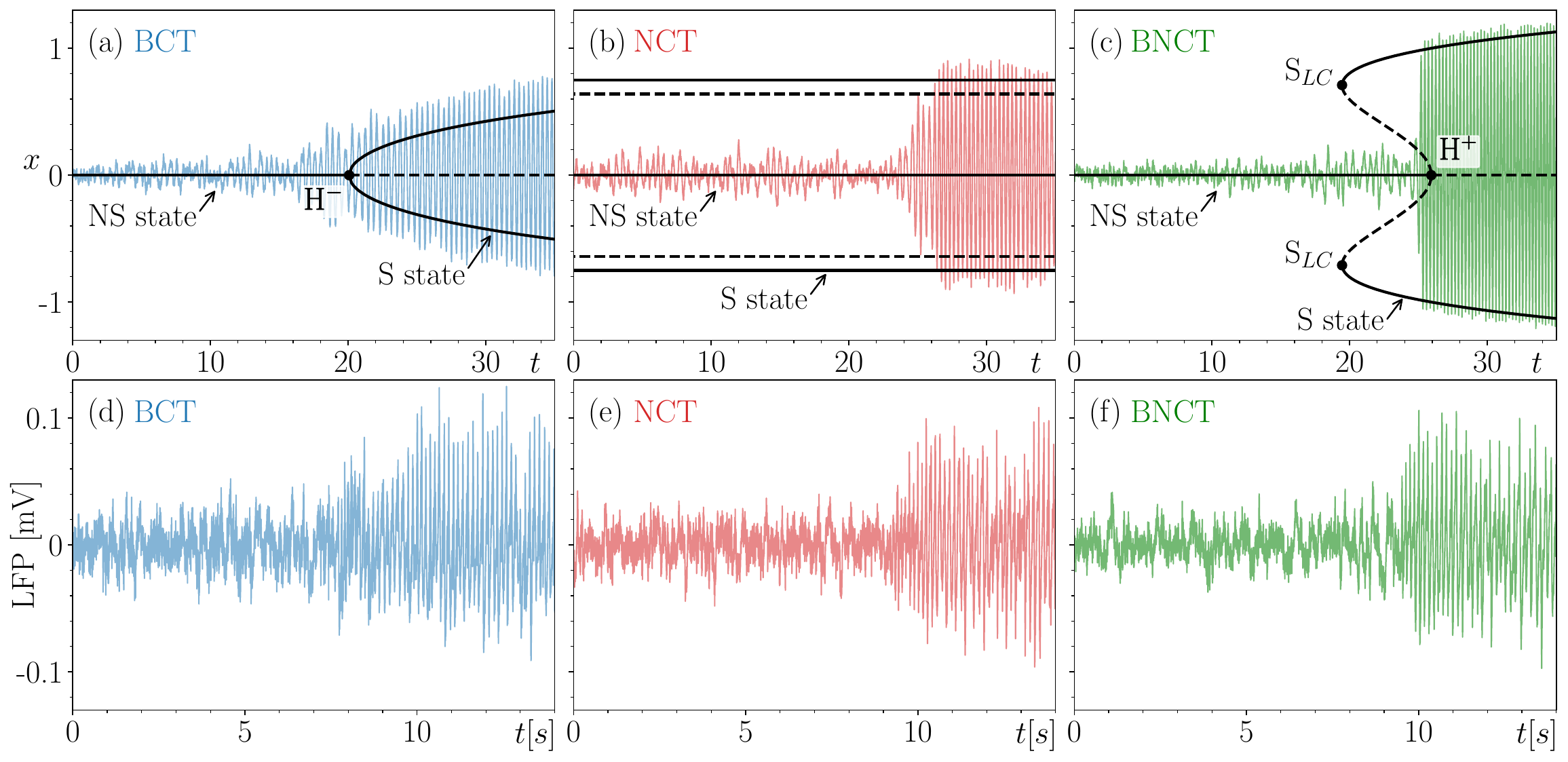}
    \caption{
    The top row shows examples of the three CT-types from the NS to S state in the model's output. (a) Bifurcation-induced CT due to drifting across a supercritical bifurcation point H$^-$.
    (b) Noise-induced CT in the bistability region between the stable  NS and S states. 
    (c) Bifurcation/noise-induced CT due to a combination of drifting across a subcritical bifurcation point, H$^+$, and bistability.
    The solid (dashed) lines represent the stable (unstable) NS and S states of the noise-free system~\eqref{eq:BautinShearNF_complex}, and S$_{LC}$ in (c) denotes the fold of limit cycles bifurcation point. The coloured time series show the outputs of the mathematical model~\eqref{eq:BautinShearNF_cartesian_noise}.
    The bottom row shows the corresponding CTs in the voltage recordings from GAERS, as classified by our framework; see Fig.~{S-19} in [\citeonline{Flynn25_Supplement}] for further details.
    }
    \label{fig:CTtypes}
\end{figure}

In most mathematical models of seizure dynamics, NS states are represented by noisy oscillations about a stable equilibrium point near $0~\text{mV}$, and S states are represented by noisy motion along a non-stationary stable state, such as a stable limit cycle, with amplitude much larger than the amplitude of the noisy oscillations near $0~\text{mV}$.
These mathematical models can be separated into the following groups: process-based models which are derived from a number of physical processes occurring in the brain to produce seizure-like activity~\cite{DaSilva03EpilepsyDynDis,DaSilva03dynamicaldisease,suffczynski04EpilepsyDynamics,suffczynski05EpilepticTransitions}, phenomenological models which are designed to mimic real seizure activity in voltage recordings without taking into account the underlying processes in the brain~\cite{meiselkuehn2012scaling,milanowski2016seizures,Jirsa_2014_NatureOfSeizDyn,jungesTerry20epilepsySurgery_SubCSupC_Network,harringtonTerry24_SubCSupC_Network,qinBassett24_SubSupC_Network}, and a mixture of both\cite{byrne2022modelEEG}.
CTs from the NS to S state in a mathematical model can be related to processes in the brain and largely  categorised into:
\begin{itemize}
    \item
    Bifurcation-induced CTs (BCTs): The brain begins in the NS state, which is the only stable state. However, the NS state becomes unstable when the brain slowly drifts across a supercritical Hopf bifurcation point due to changing environmental or structural parameters. This causes the brain to transition to the stable S state, which  gradually emerges from the bifurcation point and is the only stable state beyond it. See Fig.~\ref{fig:CTtypes}~(a).
    \item 
    Noise-induced CTs (NCTs): The brain is either bistable, meaning that it has two stable states, NS and S, that coexist, or excitable, in which case the S state is metastable (transient). A sudden transition from the NS to S state is triggered by noise  alone, in the form of endogenous or exogenous disturbances. See Fig.~\ref{fig:CTtypes}~(b) for the bistable case.
    \item 
    Bifurcation/noise-induced CTs (BNCTs): The brain becomes bistable as it slowly drifts towards a subcritical Hopf bifurcation, at which point the NS state loses stability and the S state becomes the only stable state. A combination of noise, bistability, and the imminent instability of the NS state causes a sudden transition to the S state. See Fig.~\ref{fig:CTtypes}~(c).
\end{itemize}

\noindent
In part \hyperref[md:model]{M2} of the Methods section, we introduce our mathematical model, which is phenomenological and takes the form of two coupled stochastic differential equations. 
As shown in Fig.~\ref{fig:CTtypes}, our model has equivalent NS and S states and can generate the three CT-types listed above. Most importantly, a comparison of the top and bottom rows of Fig.~\ref{fig:CTtypes} shows that the CT types in our model closely align with those observed in the voltage recordings from GAERS. More generally, the model mimics three important characteristics of real seizure activity  in the voltage recordings from GAERS:
(i) the amplitude and correlation of oscillations in the NS and S states, (ii) the range of residence times in the NS and S states, and (iii) the intrinsic timescales of the NS and S states. 
Part \hyperref[md:tuning_model]{M3} of the Methods section explains how we tune the model's parameters to produce realistic seizure characteristics.

\subsection*{Step 2: CT detection algorithm for the model's output and voltage recordings 
\label{sssec:UsingCTSeizAlgorithm}}

To ensure consistency in detecting CTs between the NS and S states,  we develop a single algorithm to detect CTs in both the model's output and the voltage recordings. In doing so, we account for the fact  that both are noisy time series.

\begin{figure}[t]
    \centering
    \includegraphics[width=0.95\textwidth]{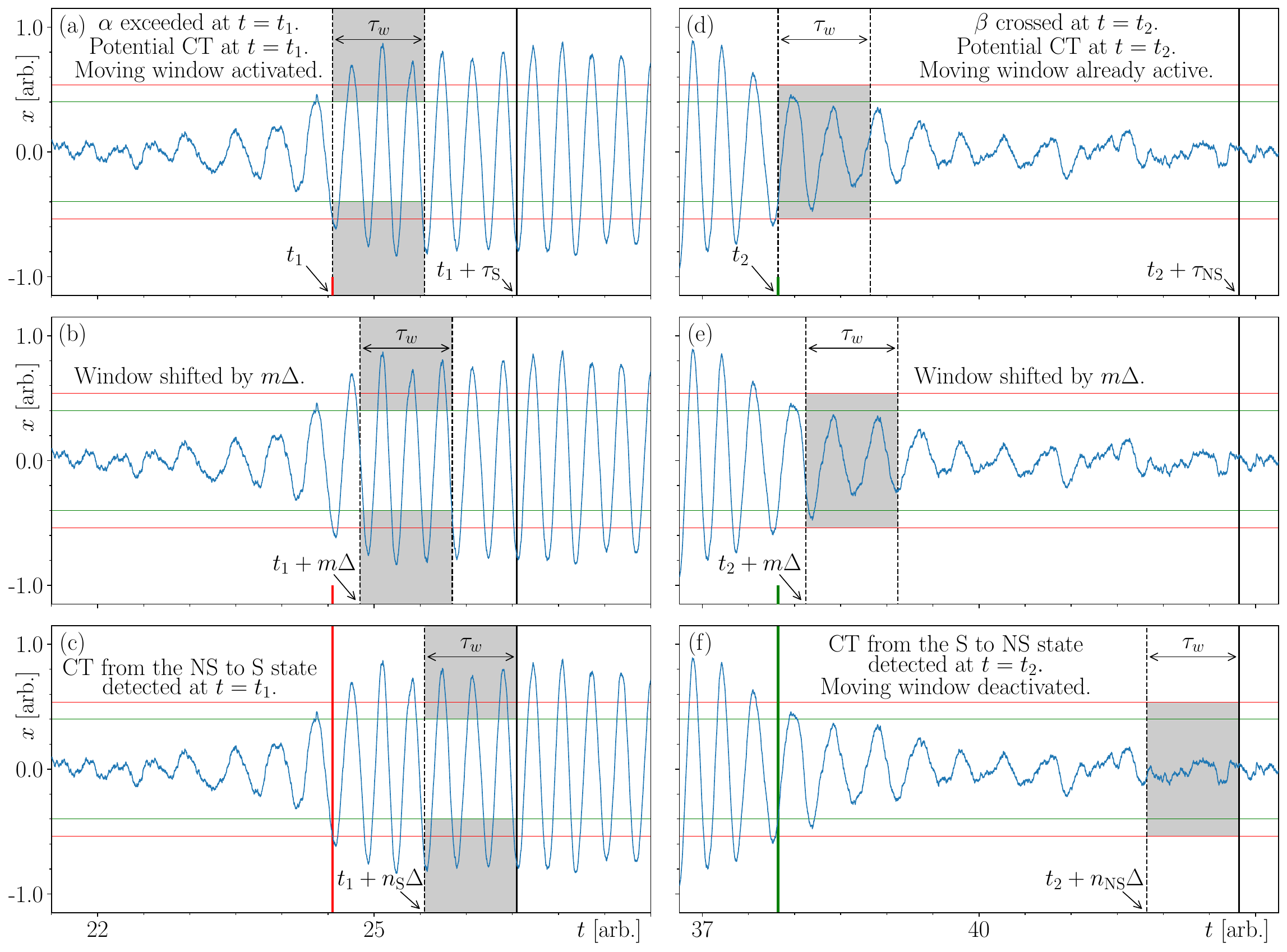}
    \caption{
    Illustrating how the CT detection algorithm works using  the model's output from Fig.~\ref{fig:_Model_vs_Data_3characteristics_}~(a) as an example. 
    (a)-(c) Detecting a CT from the NS to S state at time $t = t_1$. (d)-(f) Detecting a CT from the S to NS state at time $t=t_2$. 
    The red and green vertical lines indicate the times of detected  CTs from the NS to S state, and vice versa, respectively.
    The red and green horizontal lines indicate the  thresholds $\alpha$ and  $\beta$, respectively. 
    The detection algorithm  details and parameters are specified in part \hyperref[md:ConstructingAlgorithm]{M4} and Table~\ref{tab:AlgorithmParameters} of the Methods section, respectively.
    }
    \label{fig:Algorithm_in_use_}
\end{figure}

\begin{figure}[t!]
    \centering
    \includegraphics[width=0.8\textwidth]{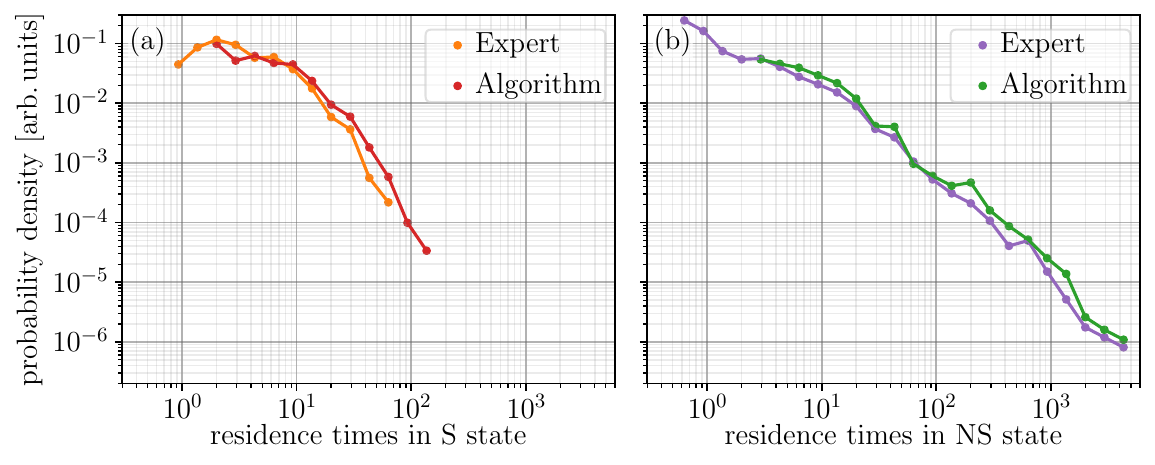}
    \caption{Probability density of residence times in the (a) S state and (b) NS state for rat S according to the (orange and purple) expert's annotations and (red and green) our CT detection algorithm from Step 2.}
    \label{fig:ResTimes_RatS_Cian_vs_Algorithm}
\end{figure}

To detect CTs in a noisy time series {of an observable}, one typically specifies a {threshold value for the observable} and defines a CT to a different state as a crossing of this threshold {by the time series}. However, such a single threshold approach has its shortcomings. 
{
A single threshold can be crossed multiple times within a short time interval due to noise, oscillatory states, or a combination of both, which could result in false detections. These shortcomings
}
arise in the case of CTs to the S state and we address them as follows:
\begin{itemize}
    \item 
    We use the concept of a `non-ideal relay' with two thresholds, which we call $\alpha$ and $\beta$, and define CTs between the NS and S states in terms of successive crossings of the two thresholds~\cite{Pokrovskii12systemswHysteresis}. 
    \item 
    We use a `moving window analysis' with a time window of duration $\tau_w$ to assess how long the system remains in the new state following successive crossings of the two thresholds. 
\end{itemize}
In short, the algorithm detects a CT from the NS to S state at time $t=t_{1}$ if the time series exceeds the upper threshold $\alpha$ and then continues to exceed the lower threshold $\beta$ frequently enough for a period of at least $\tau_{\text{S}}$. Similarly, the algorithm detects a CT from the S to NS state at time $t=t_{2}$ if the time series falls below the lower threshold $\beta$ and then does not exceed the upper threshold $\alpha$ for a period of at least $\tau_{\text{NS}}$. 
Figure~\ref{fig:Algorithm_in_use_} provides an example that illustrates how the detection algorithm works.
Part~\hyperref[md:ConstructingAlgorithm]{M4} of the Methods section provides precise definitions of when the algorithm detects a CT and specifies the algorithm's technical details. Furthermore,  we refer to Secs.~{S1}-{S5} in [\citeonline{Flynn25_Supplement}] for additional details on how we apply the algorithm consistently to both the model's output and the voltage recordings, and how we account for artefacts that sometimes appear in the voltage recordings. 
In particular, in Sec.~{S5}, we tune the algorithm parameters to obtain excellent agreement between the algorithm and the expert annotations of electrical seizure onset and offset times in the voltage recordings from different rats.
This is illustrated in Fig.~\ref{fig:ResTimes_RatS_Cian_vs_Algorithm} for rat S, and in Figs.~{S-8} and {S-9} in [\citeonline{Flynn25_Supplement}]  for rats T and K. 
{Further details on the performance of our algorithm, including a detailed Receiver Operating Characteristic (ROC) analysis, can be found in Flynn \textit{et al.}~\cite{flynn2026detecting}.}


\subsection*{Step 3: Identifying distinctive properties of the model's output that differentiate between the three CT-types \label{ssec:ClassificationResultsModel}}

When a system is perturbed away from its stable state — whether stationary or not — it will attempt to return to this state. However, the closer the system is to a bifurcation point at which  the state loses stability, the less stable it becomes and the longer it takes to return.
Therefore, as a noisy system gradually approaches a bifurcation point, one would expect to see an increase in the variance and autocorrelation of the observed time series.
This phenomenon is known as `critical slowing down' (CSD) and the accompanying increases in the variance and autocorrelation are commonly referred to as `early warning signals' (EWSs). 
Scheffer \textit{et al.}~\cite{scheffer2009early} recognised that EWSs can be detected in time series recordings of a single observable of a complex system, allowing researchers to detect CSD without knowledge of the entire state of the system.
This approach has been widely successful in studies of climate and environmental systems, enabling researchers to forecast when BCTs may occur~\cite{ditlevsen2023ews,vandijkstra2024ewsAMOC}. 
A recent review and advances on this topic can be found in the work of Dakos \textit{et al.}~\cite{Dakos_2024_ESD}, Ashwin \textit{et al.}~\cite{ashwin2025_ews_skill}, and Bury \textit{et al.}~\cite{bury2021deepEWS}.

In contrast, CTs in the brain occur on timescales that are several orders of magnitude faster than those in climate and environmental sciences. Additionally, we also consider NCTs that lack the EWSs of BCTs.
Therefore, our approach differs from previous work on EWSs in three key ways:
First, we analyse the time series around a CT, that is, shortly before and shortly after a CT, rather than just prior to a CT. This is because our primary goal is to identify the CT-type rather than to foresee one that is about to occur within seconds. 
Second, we consider the following properties of the time series and their slopes around a CT from the NS to S state:
\begin{enumerate}[labelindent=20pt,itemindent=1em,leftmargin=!]
    \item[(TSP1)\label{TSP1}] Gaussian variance (GV). 
    \item[(TSP2)\label{TSP2}] Base-10 logarithm of the Gaussian variance (log$_{10}$GV). 
    \item[(TSP3)\label{TSP3}] Lag-1 autocorrelation (AC). 
    \item[(TSP4)\label{TSP4}] Base-10 logarithm of the Gaussian variance of the lag-1 autocorrelation (log$_{10}$GV(AC)).
\end{enumerate}
We collectively refer to \hyperref[TSP1]{(TSP1)}-\hyperref[TSP4]{(TSP4)} as the `time series properties' or TSPs for short. 
We also consider the `slopes of the TSPs', and denote them by $m$(TSPs, $t_m$),  where $t_m$ is the chosen 
{duration of the}
time interval over which the slope is calculated. 
In other words, the $m$(TSPs, $t_m$) carry information about the trends in TSPs over the past $t_m$ seconds.
Third, and most importantly, we use the TSPs and  $m$(TSPs, $t_m$) to classify CTs from the NS to S state as (i) bifurcation-induced, (ii) noise-induced, or (iii) bifurcation/noise-induced.
{
To distinguish between behaviours that appear similar on a linear scale but are distinctly different near zero, we also work with the logarithm of certain TSPs.
}

\textbf{Selection criteria for CTs:}
We select 100 examples of each CT-type from the NS to S state in the model's output according to the following {\em selection criteria}, bearing in mind that the CT detection algorithm detects each CT at a time $t=t_1$, 
\begin{itemize}
    \item[SC1: \label{SC1}] 
    The model is in the NS state for all $t \in \left[ t_{1}+T^{-}, \, t_{1} \right)$, where $T^{-}=-30$.\\
    $|T^{-}|$ is the {duration} of the time series considered prior to each detected CT.
   \item[SC2: \label{SC2}] 
    The model is in the S state for all $t \in \left[ t_{1}, \, t_{1}+T^{+} \right]$, 
    where $T^{+}=10$.\\
    $T^{+}$ is the {duration} of the time series considered after each detected CT.
    \item[SC3: \label{SC3}]
    There is no crossing of the `on threshold' of the CT detection algorithm, $\alpha$, for any $t \in \left[ t_{1}+T^{-}, \, t_{1} \right)$, where $T^{-}=-30$.
\end{itemize}
Note that these criteria exclude some CTs detected by the detection algorithm. 
This is because our CT-type classifier requires input from a time series with sufficient time prior to  and after a CT from the NS to S state. There is also a biological reason for excluding CTs before which insufficient time 
has been spent in the NS state: the brain may not have 
fully returned to the NS state beyond what can be seen in the voltage recording. 
Hence the criteria~\hyperref[SC1]{SC1} and~\hyperref[SC2]{SC2}. The criterion~\hyperref[SC3]{SC3} excludes CTs that are preceded by  crossings of the `on threshold' that are too brief to be detected as a CT. We refer to these as `almost-occurring CTs' in part~\hyperref[md:ConstructingAlgorithm]{M4} of the Methods section; see Sec.~{S10} in [\citeonline{Flynn25_Supplement}] for further information.
{
We choose $T^{-}=-30$ and $T^{+}=10$ to examine the behaviour of TSPs and $m($TSPs, $t_{m})$ versus $T$ around different CTs and for different durations $t_{m}$.
Additional justification for the choice of $T^-$ is provided in Sec.~{S4} in [\citeonline{Flynn25_Supplement}]. 
}
%

\begin{figure}[t]
    \centering
    \includegraphics[width=0.99\textwidth]{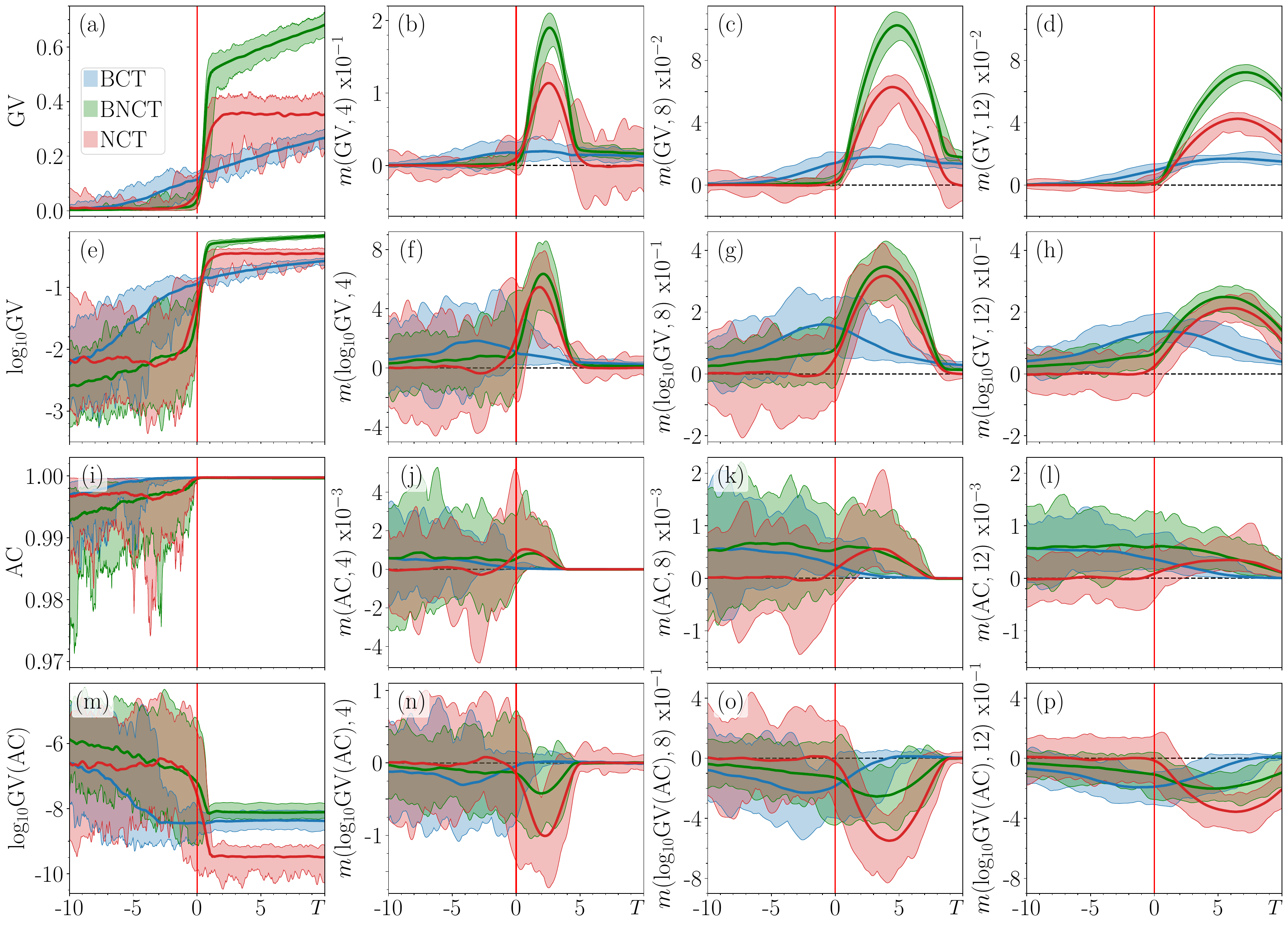}
    \caption{
    Mean values (solid curves) and spread (shaded regions surrounding each curve) of the different TSPs and $m$(TSPs, $t_m$), before and after (in {blue}) BCTs, (in {red}) NCTs, and (in {green}) BNCTs, plotted versus $T=t-t_{1}$, where $t_1$ is the time at which the detection algorithm detects each CT. The first column corresponds to the TSPs, the remaining three columns correspond to the $m(\text{TSPs},~t_{m})$ for $t_{m}=4$, $8$, and $12$ respectively. 
    }
    \label{fig:EWSslopes_4_8_12}
\end{figure}

{\textbf{Computing TSPs and $m$(TSPs, $t_m$) near selected CTs:}}  
When computing TSPs and $m$(TSPs, $t_m$) near each CT, it is convenient
to work with a new time, $T=t-t_{1}$, where $T \in \left[ T^{-}, T^{+} \right]$ and each CT is detected at $T=0$. 
{In other words, $T$ is a signed distance to a CT, taking negative values before a CT and positive values after.}
A single value of each TSP is computed over a time interval {whose duration is} denoted by $t_{w}$ and referred to as the `window {duration}'. We consider a fixed $t_{w} = 1$ for both the model's output and voltage recordings; 
see Fig.~{S-11} in [\citeonline{Flynn25_Supplement}] for justification.
A single value of each $m$(TSPs, $t_m$) is computed over a longer 
{duration}, $t_{m} > t_{w}$, referred to as the `slope {duration}', and we examine {the effect of} different
$t_m$.  
Thus, for a given choice of $T^{-}$ and $T^{+}$, we obtain values of \hyperref[TSP1]{(TSP1)}-\hyperref[TSP3]{(TSP3)} in the interval $\left[ T^{-} + t_{w}, T^{+} \right] = \left[ T^{-} + 1, T^{+} \right]$, \hyperref[TSP4]{(TSP4)} in the interval $\left[ T^{-} + 2t_{w}, T^{+} \right] = \left[ T^{-} + 2, T^{+} \right]$, and all $m$(TSPs, $t_m$) in the interval $\left[ T^{-} + 2t_{w} + t_m, T^{+} \right] = \left[ T^{-} + 2 + t_m, T^{+} \right]$.
Figure~\ref{fig:EWSslopes_4_8_12} summarises the results obtained for all the TSPs and $m$(TSPs, $t_m$) for the 100 selected examples of each CT-type.
There, we plot the mean values (solid curves) and spread (shaded bands) as a function of $T$, and use a different colour for each CT-type. 
Additionally, Figs.~\ref{fig:_compution_of_GV_mGV_picture_} and \ref{fig:_compution_of_logGVAC_mlogGVAC_picture_} in part~\hyperref[md:TSPs]{M5} of the Methods section illustrate the process of computing TSPs and $m$(TSPs, $t_m$). Section~{S6} in [\citeonline{Flynn25_Supplement}] provides additional details, and Fig.~{S-10} in [\citeonline{Flynn25_Supplement}] shows how the TSPs behave for a single example of each CT-type.

\textbf{Differentiating between the three CT-types based on TSPs and $m$(TSPs, $t_m$):}
We now focus on how the TSPs and $m$(TSPs, $t_m$) can be used to distinguish between the three CT-types.
For instance, we see in Fig.~\ref{fig:EWSslopes_4_8_12}~(c) that
three separate (non-overlapping) bands of $m(\text{GV},~t_{m})$ emerge near $T = 5$, where each band corresponds to a different  CT-type. 
This allows one to visually distinguish between each CT type. 
Comparing Figs.~\ref{fig:EWSslopes_4_8_12}~(b), (c), and (d), we see that as $t_{m}$ increases, the width of these bands decreases, thereby improving our ability to visually distinguish between different CT types. 
As another example, we can distinguish between the three CT types by combining the following two observations when $T>2$: 
BNCTs can be distinguished from BCTs and NCTs based on the separation between corresponding $\text{GV}$ bands in Fig.~\ref{fig:EWSslopes_4_8_12}~(a), and NCTs can be distinguished from BCTs and BNCTs based on the separation between corresponding $\text{log}_{10}\text{GV(AC)}$ bands in Fig.~\ref{fig:EWSslopes_4_8_12}~(m).
By combining several TSPs and $m$(TSPs, $t_m$) we may distinguish between the three CT-types at earlier $T$ values, even before a CT occurs.


\subsection*{Step 4: Training and optimising a Support Vector Machine on CTs of known type in the model's output.}

It is convenient to automate the process of differentiating between the three CT-types based on their distinctive TSPs and $m$(TSPs, $t_m$) using a machine learning classifier. 
We chose a Support Vector Machine (SVM) because it is an explainable {and deterministic} machine learning technique. {There is no ambiguity about how it works and, if the same CT is analysed repeatedly, it will always be assigned to the same CT-type}. 
An SVM classifier with $F$ features can separate different classes in an $F$-dimensional space of features.
In our case, the features are the TSPs and/or $m$(TSPs, $t_m$), and the different classes correspond to the three CT-types.
In other words, the feature space is divided into three regions, with   each region corresponding to a different CT-type. Each CT is then classified as a BCT, NCT or BNCT, depending on which region of the feature space its features fall into.

\begin{figure}[t]
    \centering
    \includegraphics[width=\textwidth]{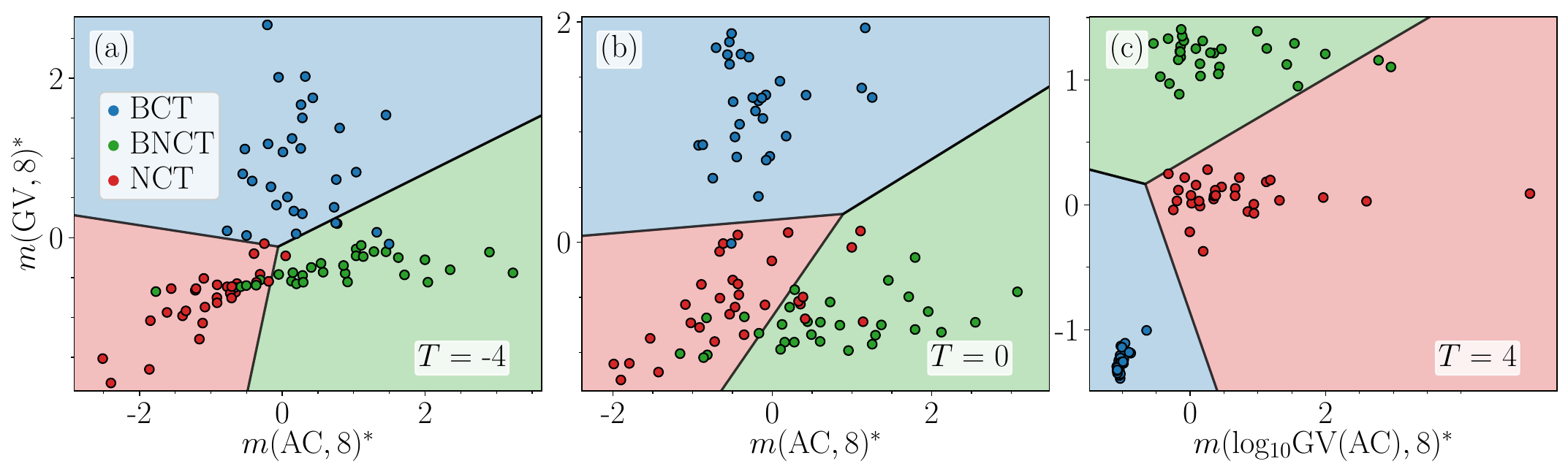}
    \caption{
    An illustrative example of how a two-feature SVM works, with  the features specified on each axis. The different coloured points correspond to ({blue}) the true  BCTs, ({green}) the true BNCTs, and ({red}) the true NCTs. This information is known to us but not to the SVM.
    The three different coloured regions are where the SVM classifies CTs as ({blue}) BCTs, ({green}) BNCTs, and ({red}) NCTs. The accuracy of the classification depends on the chosen features and whether classification is performed  (a) before a detected CT, (b) at a detected CT, or (c) after a detected CT. 
    }
    \label{fig:SVM_optimal_lines_illustration_tm8_}
\end{figure}

\textbf{Illustrating how SVM classifiers work: }
For illustrative purposes, we consider a simple SVM classifier with only two features ($F=2$), to show how CTs in the model's output are classified and how the classification changes for different values of {the signed distance to a CT,} $T$. 
Each panel in Fig.~\ref{fig:SVM_optimal_lines_illustration_tm8_} is divided into three regions separated by the optimal lines, and these regions are coloured according to each classification class: {blue} for BCT, {red} for NCT and {green} for BNCT.
The different coloured points correspond to the true BCTs in {blue}, 
true NCTs in {red}, and true BNCTs in {green}. 
In the case of the model's output, these are known to us but not to the SVM. 
The features in each panel are chosen from Fig.~{S-12}~(b) in [\citeonline{Flynn25_Supplement}] as the two most important features at the given value of $T$. 
{
By `most important features', we mean those with the highest  `mean permutation importance' (MPI), as described in more detail in part~\hyperref[md:SVMtrain]{M6} of the Methods section.
}
Figure~\ref{fig:SVM_optimal_lines_illustration_tm8_} shows that even a simple two-feature SVM classifier can successfully distinguish between the three seizure generation mechanisms in noisy time series, with fewer CTs being misclassified as $T$ increases. 

\textbf{Training SVM classifiers: }
We use the TSPs and $m$(TSPs, $t_m$) of the model's output, where we know how each CT is generated, as features to train different SVMs to classify the three CT-types from the NS to S state.
To determine how classifications vary depending on what combination of features are used, we train three types of SVM classifier:
\begin{itemize}
    \item 
    {\bf Type-1} using the four TSPs, so that $F=4$.
    \item
    {\bf Type-2} using the four $m$(TSPs, $t_m$), so that that $F=4$.
    \item
    {\bf Type-3} using the four TSPs and the four $m$(TSPs, $t_m$), so that $F=8$.
\end{itemize}
The training procedure is the same for each SVM classifier and is outlined in part~\hyperref[md:SVMtrain]{M6} of the Methods section.

{\textbf{Testing and optimising SVM classifiers: }
We test each SVM classifier by examining its ability to classify CTs in unseen model's output.
We refer to the fraction of correctly classified CTs as the \textit{SVM accuracy}, and examine how it changes with {the signed distance to a CT,} $T$, for each SVM type. 
{For SVMs type-2 and 3, we also examine the dependence on the duration, $t_m$, of the time interval over which $m$(TSPs, $t_m$) are calculated.}
When the SVM accuracy reaches $1$, we say that the SVM classifier achieves perfect accuracy at this value of $T$. 
In Fig.~\ref{fig:SVM_Accuracy_slopes_slopesANDvalues} we plot the SVM accuracy versus $T$ for SVMs type-1 and 2 in (a) and SVM type-3 in (b).
The value of $T$ at which a given SVM classifier achieves perfect accuracy for the first time is indicated by a coloured dot. 
\begin{figure}[t]
    \centering
    \includegraphics[width=\textwidth]{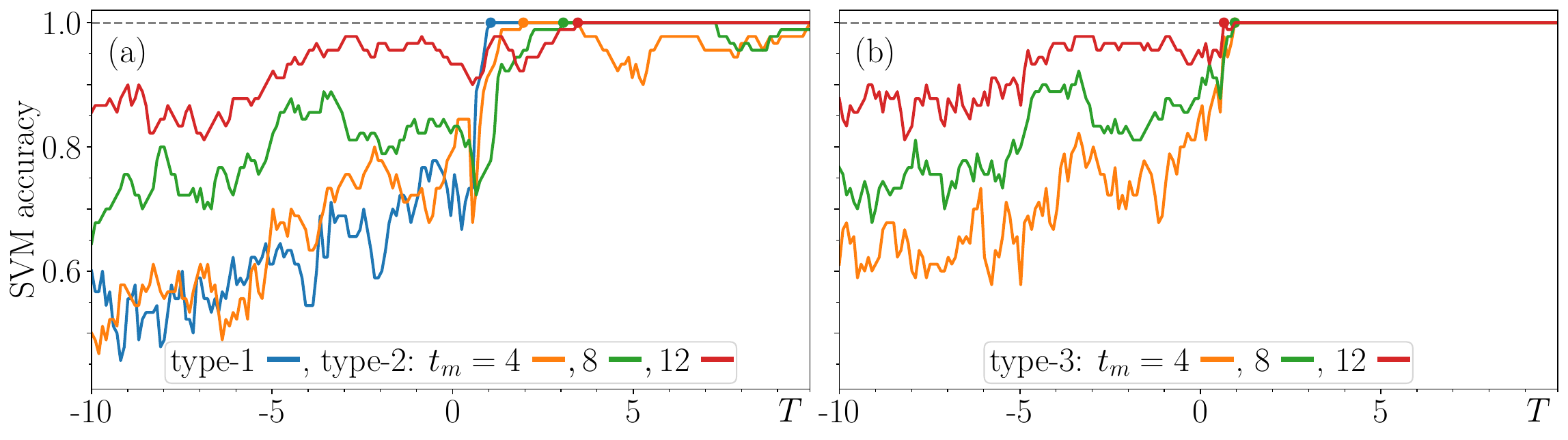}
    \caption{
    SVM accuracy in classifying CTs in unseen model's output vs. $T$ for (a) SVM type-1 and 2 and (b) SVM type-3.
    For SVM type-2 and 3, we include the dependence on different {durations} $t_{m}$ specified in the plot legends.
    The dots indicate when each SVM achieves perfect accuracy for the first time.
    }
    \label{fig:SVM_Accuracy_slopes_slopesANDvalues}
\end{figure}
SVM type-1 performs relatively poorly for most $T < 1$, but achieves perfect accuracy around $T=1$, before SVM type-2 for each $t_{m}$, and maintains its perfect accuracy for $T > 1$. 
Increasing $t_{m}$ increases the accuracy of SVMs type-2 and 3. 
While SVM type-3 maintains its perfect accuracy after first achieving it, SVM type-2 does not for 
{shorter durations} $t_m$.
It is clear that SVM type-3 benefits from the larger accuracy of SVM type-2 for $T < 1$ and the consistently perfect accuracy of SVM type-1 for $T > 1$. 

We therefore conclude that SVM classifier type-3 is the most suitable choice for classifying seizure generation mechanisms. 
However, an optimal choice of $t_{m}$ must be made in the context of the voltage recordings. 
Although SVM type-3 is more accurate for $T < 1$ when $t_{m}=12$ than for $t_{m}=8$, we can be less restrictive with our choice of $T^{-}$ and $T^{+}$ for 
{shorter durations} $t_{m}$, allowing more CTs to be classified. 
This results in a trade-off between maximising either the SVM accuracy or the number of CTs that can be classified.  
In the interest of classifying as many seizure generation mechanisms as possible with a reasonably high level of accuracy, we choose to use the SVM classifier type-3 with $t_{m}=8$ and refer to it as the `model-trained CT-type classifier' to emphasise that it is trained on the model's output.


\subsection*{Step 5: Using model-trained CT-type classifier to classify unknown seizure generation mechanisms in the voltage recordings of real seizure activity 
\label{ssec:ClassificationResultsGAERS}}

We now present our main results:
The classification of seizure generation mechanisms in voltage recordings of real seizures in GAERS, using a CT-type classifier trained on the model's output. 
To describe the voltage recordings near CTs from the NS to S state, we use the same notation as for the model's output. 
Specifically, we use {the signed distance form a CT, namely} 
$T \in \left[ T^{-}, T^{+} \right]$, where $T=t-t_{1}$, and $T=0$ is the time when the CT detection algorithm detects a CT from the NS to S state.
{
In part~\hyperref[md:CianSeizureAnnotate]{M1} of the Methods section, we  describe all the  voltage recordings that were available and explain why only three GEARS, namely S, T, and K, were included in our analysis.
}

\textbf{Stage 1 (Detection of CTs): } 
We begin by using our detection algorithm to detect CTs between the NS and S states in the voltage recordings of rats S, T, and K;
The algorithm parameters specified in Table~\ref{tab:AlgorithmParameters} of the Methods section are chosen to maximise agreement with expert annotations of seizure onset and offset times, according to the procedure outlined in Sec.~{S5} in [\citeonline{Flynn25_Supplement}]; 
{see also Flynn \textit{et al.}~\cite{flynn2026detecting} for a detailed ROC analysis of the algorithm's performance.}
We denote the total number of detected CTs from the NS to S state for a given rat by $N_{\text{det}}$, where `det' stands for detected. We refer to this set of CTs as the `detected CTs'.

\textbf{Stage 2 (Filtering of detected CTs): } 
We classify a detected CT only if it satisfies the following {\em classification conditions}: 
\begin{enumerate}[labelindent=10pt,itemindent=1em,leftmargin=!]
    \item[C1: \label{C1}] 
    The brain remains in the NS state for all $T \in \left[ T^{-}, \, 0 \right)$. 
    \item[C2: \label{C2}] 
    The brain remains in the S state for all $T \in \left[ 0, \, T^{+} \right]$. 
    \item[C3: \label{C3}] 
    There is some overlap between what the algorithm identifies as the S state and the expert annotation of the S state. 
    \item[C4: \label{C4}] There is no crossing of the `on threshold' of the CT detection algorithm,
    $\alpha$, for any $T \in \left[ T^{-}, \, 0 \right)$.
\end{enumerate}
Conditions~\hyperref[C1]{C1} and \hyperref[C2]{C2} are consistent with the selection criteria~\hyperref[SC1]{SC1} and \hyperref[SC2]{SC2}}. They ensure that the analysed TSPs and $m$(TSPs, $t_m$) correspond to an individual CT. 
Condition~\hyperref[C3]{C3} ensures that we only classify CTs where our algorithm agrees with the expert annotation that there is a CT. 
Condition~\hyperref[C4]{C4} is consistent with the selection criterion~\hyperref[SC3]{SC3}. 
Additionally, in the case  of voltage recordings, condition~\hyperref[C4]{C4} is used to exclude scenarios in which `artefacts' that obscure the features used by the classifier occur just prior to a detected CT.
Artefacts are absent from the model's output, but appear in the voltage recordings;
see Figs.~{S-15} and {S-20} in [\citeonline{Flynn25_Supplement}] for more details.
The subset of detected CTs that satisfy all four conditions  are referred to as the `filtered CTs'. The number of filtered CTs for a given rat  and choice of $T^{-}$ and $T^{+}$ is denoted by $N_{\text{filt}}(T^{-}, T^{+})$, where `filt' is short for filtered.

\textbf{Stage 3 (The main experiment: Classification of filtered CTs from voltage recordings of real seizures in GAERS): } 
For our choice of $T^-$, $T^+$, $t_w=1$ and $t_{m}=8$, we  obtain all TSPs and $m$(TSPs, $t_m$) of the filtered CTs for $T \in \left[ T^{-}+2t_w + t_m, T^{+} \right] = \left[ T^{-} + 10, T^{+} \right]$. This is the time interval where we can apply our model-trained CT-type classifier to the filtered CTs. Next, we use $N_{\text{type}}(T)$ to denote the number of CTs classified at time $T$ as type = BCT, BNCT or NCT , and compute the proportions of: 
\begin{itemize}
    \item Filtered CTs to the total number of detected CTs: $N_{\text{filt}}(T^{-}, T^{+})/N_{\text{det}} \in [0,1]$.
    \item Filtered CTs classified as a particular type at time $T$ to the total number of filtered CTs: $N_{\text{type}}(T)/N_{\text{filt}}(T^{-}, T^{+}) \in [0,1]$, bearing in mind that $T \in \left[ T^{-} + 10, T^{+} \right]$.
\end{itemize}
To demonstrate the robustness of our findings with respect to the selection of time intervals around the CTs, we consider a fixed value for $T^{+}$ and different values for $T^{-}$. 
More specifically, we set $T^{+}=\tau_{\text{S}}=2$ and consider $T^{-} = -16,-14,-12,-10,$ and $-8$. 
To maximise $N_{\text{filt}}(T^{-}, T^{+})$, the values of $T^{+}$ and $T^{-}$ 
are chosen to be closer to the CT than the values used for training the CT-type classifier in selection criteria~\hyperref[SC1]{SC1} and~\hyperref[SC2]{SC2}. 
The proportions $N_{\text{type}}(T)/N_{\text{filt}}(T^{-}, 2)$, which is the main quantity of interest, can be obtained at any $T\in[-6,\, 2]$ if $T^{-} = -16$, and at just one value of $T=2$ if $T^{-} = -8$. 
Therefore, to make a meaningful comparison of CT classifications for different values of $T^{-}$, we focus on the classifications at $T=2$, that is two seconds after the time at which a CT is detected. In other words, we examine $N_{\text{type}}(2)/N_{\text{filt}}(T^{-}, 2)$ for $T^{-} = -16,-14,-12,-10,$ and $-8$, where type = BCT, BNCT or NCT is given by the model-trained CT-type classifier.
We refer to Figs.~{S-16},~{S-17}, and~{S-18} in
[\citeonline{Flynn25_Supplement}] for details of classification dependence on $T$.


\subsubsection*{Key findings}

\begin{figure}[t]
    \centering
    \includegraphics[width=\textwidth]{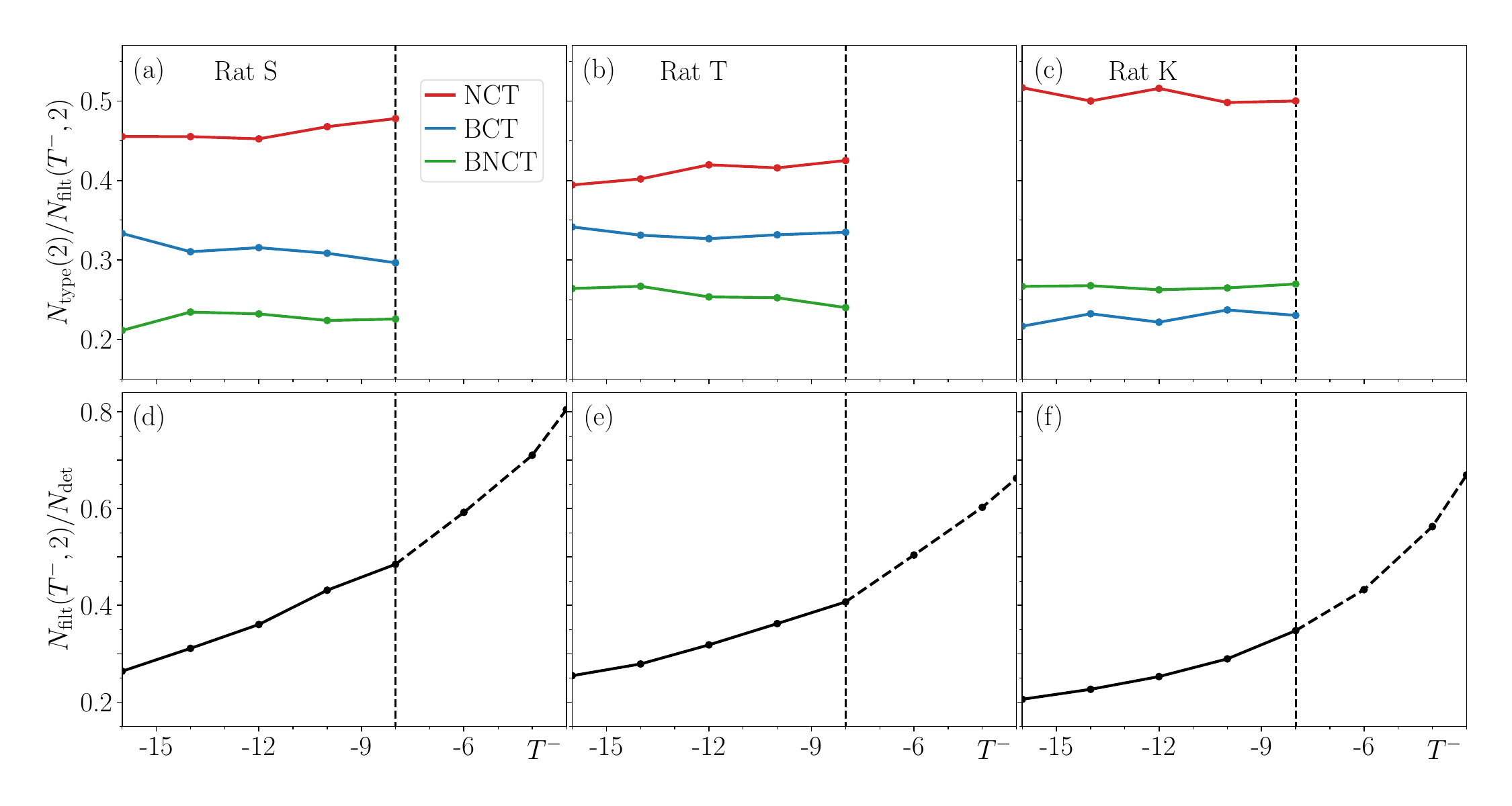}
    \caption{
    (Top row) The proportion of filtered CTs from the voltage recordings in GAERS  that are classified by the model-trained SVM, i.e. SVM type-3 with $t_w=1$ and $t_m=8$, at time $T=T^+=2$ as ({red}) NCT, ({blue}) BCT and ({green}) BNCT vs. $T^-$ for (a) rat S, (b) rat T and (c) rat K.
    Note that NCT emerges as the dominant CT-type responsible for the onset of seizures in all three rats.
    (Bottom row) The proportion of filtered CTs to the total number of  CTs detected in the voltage recordings vs. $T^-$ for (d) rat S, (e) rat T, and (f) rat K.
    Other parameters are specified in Table~\ref{tab:AlgorithmParameters}, where $\alpha$ is chosen according to the method outlined in Sec.~{S5} in [\citeonline{Flynn25_Supplement}].
    }
    \label{fig:RatSTK_Tm_vs_classfiltCTs_filtdetCTs_Tm_16_8_}
\end{figure}

The top row of Fig.~\ref{fig:RatSTK_Tm_vs_classfiltCTs_filtdetCTs_Tm_16_8_} shows how the proportions of all filtered CTs that are classified as ({red}) NCTs, ({blue}) BCTs, and ({green}) BNCTs vary with $T^-$. We show this for rat S in panel (a), rat T in panel (b), and rat K in panel (c).
The key finding is that NCTs are the dominant seizure generation mechanism in all three rats for all values of $T^-$. More specifically, 
\begin{itemize}
    \item NCTs account for approximately {$42-50\%$} of the filtered CTs,
    \item BCTs account for approximately {$23-34\%$} of the filtered CTs, and
    \item BNCTs account for approximately {$22-27\%$} of the filtered CTs. 
\end{itemize}

\begin{table}[t]
\centering
\def\arraystretch{1.2}
\begin{tabular}{|c||c|c|c|c|c|c|}
\hline
 & Rat S & Rat T & Rat K \\
\hline
\multirow{2}{17em}{\centering No. of CTs detected by expert} & \multirow{2}{7em}{\centering 621} & \multirow{2}{7em}{\centering 1593} & \multirow{2}{7em}{\centering 1136}  \\ 
& & &  \\
\hline 
\multirow{2}{17em}{\centering No. of CTs detected by the algorithm ($N_{\text{det}}$)} & \multirow{2}{7em}{\centering 466} & \multirow{2}{7em}{\centering 1115} & \multirow{2}{7em}{\centering 874}  \\ 
& & &  \\
\hline
\multirow{2}{10em}{\centering No. of filtered CTs ($N_{\text{filt}}(-8,2)$)} & \multirow{2}{7em}{\centering 226} & \multirow{2}{7em}{\centering 454} & \multirow{2}{7em}{\centering 304}  \\ 
 & & &  \\
\hline
\multirow{3}{15em}{\centering Percentage of filtered CTs classified as type = BCTs, BNCTs, or NCTs ($(N_{\text{type}}(2)/N_{\text{filt}}(-8,2)) \times 100$)} & \multirow{3}{7em}{\centering $29.6$, $22.6$, $\mathbf{47.8}$} & \multirow{3}{7em}{\centering $33.5$, $24.0$, $\mathbf{42.5}$} & \multirow{3}{7em}{\centering $23.0$, $27.0$, $\mathbf{50.0}$}  \\ 
& & &  \\ 
& & &  \\ 
\hline
\end{tabular}
\caption{
\label{tab:FinalResults} 
Our key findings in numbers. In the bottom row, the dominant CT-type for each rat is highlighted by the bold font.
}
\end{table}

In the bottom row of Fig.~\ref{fig:RatSTK_Tm_vs_classfiltCTs_filtdetCTs_Tm_16_8_},
we show how the proportion of filtered CTs to the total number of detected CTs increases with $T^-$ for a fixed $T^{+}=2$. 
For the 
{maximum}
value of $T^{-}$, which is $-8$ and is a technical restriction of our classifier, approximately $50\%$ of the detected CTs are retained by filtering (i.e., are classifiable) for rat S in panel (d), compared to approximately $40\%$ for rat T in panel (e), and approximately $35\%$ for rat K in panel (f), partially due to many artefacts appearing. 
This raises the question of how  the exclusion of detected CTs that did not satisfy the classification conditions~\hyperref[C1]{C1}-\hyperref[C4]{C4} impacted the key findings. To address this question, we make three observations. 
First, if $T^-$ were increased to $-3$, over $80\%$ of the detected CTs for rat S and  $65\%$ for rats T and K, would meet the classification criteria \hyperref[C1]{C1}-\hyperref[C4]{C4}.
Second, the top row in Fig.~\ref{fig:RatSTK_Tm_vs_classfiltCTs_filtdetCTs_Tm_16_8_} shows that the proportion of classified CTs of a given type in all three rats remains largely unchanged as $T^-$ is varied from $-16$ to $-8$. This suggests that the trend continues towards $T^-=-3$, but this cannot be verified due to the technical requirements of our classifier.
Third, McCafferty \textit{et al.}\cite{mccafferty2025interrupt_seiz} found that approximately 50$\%$ of seizures in GAERS can be interrupted by applying a suitably chosen auditory stimulus. This suggests that, like in our dominant case of noise-induced CTs, and the least dominant case of bifurcation/noise-induced CT, these epileptic brains frequently operate in a multistable or excitable regime whereby noise alone, in the form of endogenous or exogenous disturbances, can force CTs to the S state. For these reasons, we conjecture that the key findings persist for most of the detected CTs, including those excluded by the classification conditions~\hyperref[C1]{C1}-\hyperref[C4]{C4}. 
A numerical summary of the key findings is provided in Table~\ref{tab:FinalResults}.


\section*{Discussion}

This paper presents a general framework for classifying seizure generation mechanisms as different types of critical transition (CT)
between the non-seizure state (NS state) and the seizure state (S state) in the epileptic brain.
The framework consists of five steps. 
The first step is to construct a canonical mathematical model which displays CTs that (i) closely resemble voltage recordings of real seizures and (ii) can be bifurcation-induced (BCTs), noise-induced (NCTs), or a combination of both (BNCTs). 
The second step is to construct a detection algorithm that detects CTs between the NS and S states in both the model's output and voltage recordings of real seizures. 
The third step is to identify distinctive properties of the model's output near each CT-type. 
The fourth step is to construct a machine learning CT-type classifier. 
The classifier is trained using selected time series properties of the model's output where the CT-types are known. 
In the fifth step, we {
apply} the model-trained CT-type classifier 
{
to} voltage recordings, where CT-types are unknown, in order to classify the onset of real seizures as BCTs, NCTs or BNCTs.

We applied the above framework to voltage recordings taken from inside the brain of \textit{Genetic Absence Epilepsy Rats from Strasbourg} (GAERS) - a well-established model of absence epilepsy - and, for the first time, demonstrated that consensus can be reached on which CT-types are predominantly responsible for seizure onset 
{in a given dataset.}
We found NCTs to be the dominant seizure generation mechanism - they accounted for almost $50\%$ of the seizures we classified, with BCTs and BNCTs accounting for the remainder. 
Our findings
{indicate}
that (i) seizures 
{can}
be generated by different CT-types~\cite{mccormick2001cellular,vezzani2011role,kuhlmannlehnertz2018_seizurepred_newera,gonzalez2019ionic,jiruska2023update,lehnertz2023seizcontrol}, and (ii) multistability or excitability and noise 
{can}
play a significant role in 
{epileptic brains}~\cite{DaSilva03EpilepsyDynDis}.
They 
{are also in contrast to}
the view that seizures are predominantly bifurcation-induced~\cite{maturana2020yesCSDseizure}.

An important practical implication of our findings is that different CT-types may reflect distinct neural processes involved in generating seizures. From this perspective, a seizure that is classified as an NCT is likely to start locally with respect to a recording electrode. Furthermore, there are constitutively pro-ictal neurons in the epileptic brain that can initiate a seizure based solely on noisy input in the form of exogenous or endogenous disturbances to these neurons.  The hyper-excitable neurons found in layer V of the perioral somatosensory cortex\cite{polack2007seiz_neuron}, from which GAERS seizures appear to originate\cite{meeren2002GAERS}, may be one such example. 
{Interestingly, seizures that are likely to be generated by NCTs can be abolished by a single, appropriately timed electrical pulse~\cite{jalife1979phaseresetting,glass1988clocks,winfree1980geometry}, as demonstrated in both a computational model of absence seizures~\cite{suffczynski04EpilepsyDynamics} and an in vivo experiment~\cite{osorio2009seizure}. These results suggest potential therapies for terminating seizures generated by NCTs.}
On the other hand, a seizure that is classified as a BCT is likely to emerge gradually, or propagate to the recording electrode having started elsewhere.  In this case, the slow drift in brain state is considered as the ictogenic culprit. The observed changes in brain state tens of seconds prior to seizure onset in rats and humans~\cite{mccafferty2023decreased,bai2010GAERS} may reflect such a gradual drift towards seizure onset. Therefore, 
{this slow drift and the physiological mechanisms responsible for it (e.g. changes in potassium, calcium, chloride, or GABA concentrations~\cite{scalmani2026GABA}) may suggest particular avenues for therapeutic exploration towards prevention of seizures generated by BCTs.}
BNCTs may represent combinations of the neural processes discussed above. 
Therefore, classifying seizure generation mechanisms in a given brain is not only an important fundamental question. It could transform how we model, detect, prevent, and treat seizures, with the ultimate aim of informing new clinical decision support systems and tailored seizure prevention\cite{lehnertz2023seizcontrol} and control\cite{gluckman2001seizcontrol} strategies. 

Our results provide proof of concept. They also represent an important first step towards classifying more complex human seizures. 
While studies in GAERS have identified specific neurons capable of seizure initiation, human epilepsy may have multiple potential cortical sites of origin, likely related to cortical-thalamic interactions\cite{crunelli2020roots_of_epil}. Furthermore, these sites may differ between children and adults. It remains to be seen whether different CT-types correspond to how people with epilepsy experience their seizures, or whether NCTs play a dominant role. 
In other words, additional work is required to apply our framework to the less well-understood forms of human epilepsy. 
This may involve updating our mathematical model of seizure activity with a more realistic one, augmenting our CT detection algorithm with additional time series properties such as the dynamical eigenvalue early warning signal proposed by Grziwotz \textit{et al.}\cite{grziwotz2023_dyn_e_value}, incorporating different experimental data such as heart rate variability\cite{Kimia24_HeartSeizDetect}, developing more advanced filtering techniques to remove artefacts\cite{zhang2024review}, and considering alternative definitions of NS and S states~\cite{Aravind25_preprint}.
Furthermore, our framework can be extended to classify seizure termination mechanisms~\cite{kramer2012_SeizTerm,schindler2007assessing} to gain new insight into how the brain changes in response to a seizure and how an ongoing seizure can be terminated.

More generally, our framework can encompass other types of CT, such as rate-induced CTs~\cite{ritchie2023rate}, and identify new forms of early warning signals for different CT-types. Finally, it can also be applied to complex systems beyond the brain. 
{For instance, thermoacoustic systems~\cite{sujith2017thermoacoustic} and aircraft wings~\cite{kurths2025aircraft} can undergo similar CTs from small-amplitude fluctuations to larger-amplitude oscillations through different mechanisms and our framework could be adapted to classify such CTs. }


\section*{Appendix} 

\subsection*{A1: Mathematical Background and Further Reading \label{apx}}

This appendix provides a concise, non-technical overview of the key mathematical ideas and computational methods used in this paper, together with references for further reading.

The central mathematical concept used throughout this paper is a `critical transition' (CT), which describes a sudden and large change in a system's behaviour in response to small changes in input to the system. 
We model seizure onset as a CT between two brain states, a non-seizure (NS) and seizure (S) state, and construct a framework to classify the CT-types present in voltage recordings of real seizure activity.

We consider three types of CTs: bifurcation-induced CTs (BCTs), noise-induced CTs (NCTs), and bifurcation/noise-induced CTs (BNCTs).
BCTs arise from slow changes in system parameters that cause the system's current state to lose stability and a new state to gain stability, NCTs arise from random perturbations to the system's current state which push the system across a threshold to a new state, and BNCTs arise from a combination of both effects. 
In our case, we consider the current state as the NS state and the new state as the S state, and these three CT-types provide distinct ways in which seizure onset can occur.

Our CT-type classifier relies on the concept of `critical slowing down' (CSD), a phenomenon exhibited prior to CTs involving bifurcations.
Specifically, CSD describes the phenomenon whereby a system takes an increasingly longer time to recover from small perturbations to its state when approaching a bifurcation point at which its current state loses stability.
CSD appears as increases in variance and autocorrelation in time series of the system's state. In our framework, these time series properties help us to differentiate between CT-types: BCTs and BNCTs exhibit CSD prior to a CT while NCTs typically do not. 

Our framework involves training a CT-type classifier using time series from a mathematical model which can produce BCTs, NCTs, and BNCTs under different parameter settings and whose output resembles voltage recordings of real seizure activity. Our mathematical model is based on a normal form. 
A normal form is a simplified mathematical model that captures the essential behaviour of a system near a transition. 
In our case, it provides a minimal model of how the brain transitions between NS and S states. 
We extract features from time series of the model's output associated with CSD and use these features to train a machine learning classifier which we refer to as the CT-type classifier. This classifier is then applied to voltage recordings of real seizure activity to infer the transition mechanism underlying different seizures.

While our framework relies on mathematical concepts and computational methods, the insights it provides are grounded in observable seizure behaviour. 
In this context, BCTs may correspond to gradual changes in the conditions that cause the shift to the S state, NCTs may reflect seizures triggered by intermittent disturbances to particularly sensitive groups of neurons, and BNCTs may involve both effects. This links the classification of CTs to underlying neural processes and provides a basis for relating seizure dynamics to seizure generation mechanisms.
\\
\\
\textbf{Further Reading: }
We consider the following as entry points to the biological and mathematical foundations of this work,
\begin{itemize}
    \item Scheffer \textit{et al.}~\cite{scheffer2009early}: Introduction to critical transitions and critical slowing down.
    \item Ashwin \textit{et al.}~\cite{ashwinwieczorek2017CTdef}: Advanced mathematical descriptions of critical transitions.
    \item da Silva \textit{et al.}~\cite{DaSilva03EpilepsyDynDis}: Overview of dynamical systems approaches to modeling epileptic seizures and brain activity.
\end{itemize}
    
%



\section*{Methods} 


\subsection*{M1: Obtaining and annotating seizure data from GAERS \label{md:CianSeizureAnnotate}} 

The data discussed here was obtained by Cian McCafferty, Fran\c{c}ois David, and Vincenzo Crunelli and was presented in \citeonline{mccafferty2018cortical}.
\\
\textbf{Obtaining the data:} The electrophysiological data was acquired from male \textit{Genetic Absence Epilepsy Rats from Strasbourg} (GAERS) aged between 4-7 months when in a state of `relaxed wakefulness' where seizures were experienced more often. Silicon-site electrodes were used to sample voltages (20000/second) from the \textit{ventrobasal thalamus} while the rats were able to move freely and alternate between waking, sleeping, and seizing states. This data was processed with a Plexon HST/32V-G20 VLSI-based preamplifier and associated digitization system and subsequently down-sampled to 1000 samples/second. 
\\
\textbf{Labelling the recording session:}
Each recording session was labelled using the following convention described through example, `S1K': voltage recordings from rat S on day 1 of recording during the K$^{th}$ recording session on that day.
\\
\textbf{Annotating the data:}
Spike-wave discharges (SWDs) are the electrical hallmark of absence seizures, they define electrical seizure onset and offset times. These were identified using Cambridge Electronic Design's `Spike2' software and the following procedure. In all cases, EEG at $1000\text{Hz}$ was used for this step. EEG was acquired at $1000\text{Hz}$ for fMRI (see McCafferty \textit{et al.}\cite{mccafferty2023decreased}) and behaviour was down-sampled by averaging for neuronal activity. 
Briefly, smoothening (voltage at time $t$ is set to the mean of voltages from time $t-10\text{ms}$ to time $t+10\text{ms}$) and DC removal (voltage at time $t$ is set to the original voltage minus the mean of voltages from time $t-0.1\text{s}$ to time $t+0.11\text{s}$) functions were used to reversibly visually clean the frontoparietal differential EEG. 
Then, a negative amplitude threshold (mean voltage minus $5$-$7$~standard deviations of baseline non-SWD EEG) was used to detect putative spike-wave crossing points, defined as whenever the signal crossed this amplitude threshold. 
The crossing points were then grouped into events based on the intervals between them (maximum time between initial two crossings $0.2\text{s}$, maximum time between any two crossings within an event $0.35\text{s}$, minimum of 5 crossings per event) and the defined properties of SWDs (minimum duration $0.5\text{s}$, minimum inter-SWD interval $0.5\text{s}$ merging any SWDs with shorter intervals), and subsequently, using a frequency threshold, these events were classified as SWDs (if $>75\%$ of intercrossing intervals were within a 5–12\text{Hz} range) or other (e.g., noise, sleep).
Labelled SWDs were then visually inspected for accuracy of SWD detection, as well as seizure onset and offset times. 
Periods of sleep were identified based on sharp increases in the 1–4\text{Hz} frequency band and were excluded from analysis.
Periods of non-REM sleep were rare in the recordings and not informative due to their relatively short duration.
\\
\textbf{Artefacts and how they are accounted for:}
Artefacts are high-amplitude and frequency events in the voltage recordings which can appear for different reasons, the most common being (i) movement of electrodes/wiring and (ii) movement of muscles close to the electrodes (generally for chewing). 
Some artefacts are generated by signal overload and noise from electrical devices around the recording setup, however, these are less common as the recording process was generally appropriately amplified and shielded. 
Artefacts are accounted for through the annotation method described above since high-amplitude events that do not align with the frequency profile of SWDs. 
Additionally, the visual inspection described above also includes the manual rejection of artefacts which were not excised from annotation. Entire recording sessions were excluded if more than $5\%$ of the session consisted of artefact activity.
\\
{
\textbf{Inclusion criteria for our study:}
We received voltage recordings from eight GAERS. However, recording durations and the number of recording sessions varied substantially between GAERS. To ensure sufficient statistical significance of our main results, we analysed only the three GAERS with the largest numbers of seizures according to expert annotations. These were rats S, K, and T, which had 621, 1136, and 1593 seizures respectively. 
For example, while rat S had 621 expert-identified seizure intervals, our algorithm detected 466 CTs from the NS to S state, with 226 of these meeting the criteria for classification.
In contrast, the excluded GAERS had significantly less seizures, rat Z had 66, rat H had 253, rat R had 340, rat N had 347, and rat I had 375. 
The limited number of seizure intervals in the excluded animals would have resulted in substantially fewer CTs available for classification and meaningful comparison between GAERS.  
}
\\
\textbf{Remark (comparison with human data):}
GAERS can express multiple absence seizures per minute\cite{powell2014seizure}, far exceeding the frequency of any absence seizure syndrome in humans. 
For instance, Gregor{\v{c}}i{\v{c}} \textit{et al.}\cite{gregorvcivc2022difficult} found that, for children with treatment-resistant childhood absence epilepsy, the median number of seizures per day was three.
For a review of the GAERS model, with attention to its similarities and differences to human absence seizures, see Depaulis \textit{et al.}\cite{depaulis2016genetic}.


\subsection*{M2: The model and its dynamics \label{md:model}} 

Our starting point is the Bautin (generalised Hopf) bifurcation normal form for the time evolution of a complex-valued variable $z$ as described in Kuznetsov~\cite{kuznetsov2013elements}. We augment this normal form with a shear term and describe the rate of change of $z$ as follows, 
\begin{equation}\label{eq:BautinShearNF_complex}
    \frac{dz}{dt} = \gamma \, \left( \left( \mu + i \, \omega \right) \, z + s\,\left(1 + i\,\sigma\right)\,| z |^{2} \, z - b \, | z |^{4} \, z \right),
\end{equation}
where $\mu \in \mathbb{R}$ and $s \in \mathbb{R}$ are bifurcation parameters, $\omega \in \mathbb{R}$ is the frequency of small-amplitude oscillation, $\sigma \in \mathbb{R}$ in the augmented term is the shear parameter, and $b=1$ remains fixed throughout our investigations and is no longer explicitly referred to. The timescale parameter $\gamma > 0$ allows us to adjust the timescale of the solutions. 
The parameters $s$ and $\sigma$ are chosen to be of the same sign. Although {the shear term is} not typically present in normal forms, it is essential for analysing a system's response to external disturbances~\cite{LinYoung2008shear,Wieczorek2009shear,BlackbeardWieczorek2011shear}. In our case, the shear term enables the model's output to more closely resemble the real seizure activity seen in the voltage recordings.

We now discuss the existence of candidates for NS and S states and then conduct linear stability and bifurcation analysis to uncover different parameter paths in the  $(\mu,s)$  parameter plane that give rise to the different types of CT between these states. To simplify the discussion, we transform Eq.\,\eqref{eq:BautinShearNF_complex} to a polar coordinate system by setting $z=r\,e^{i\,\theta}$ where $r \ge 0$ and 
$\theta \in \left[ 0, 2\,\pi \right)$,
\begin{align}\label{eq:BautinShearNF_polar}
\begin{split}
    \frac{dr}{dt} &= \gamma \, r \, \left( \mu + s \, r^{2} - \, r^{4} \right), \\
    \frac{d\theta}{dt} &= \gamma \, \left( \omega + \sigma \, s \, r^{2} \right).
\end{split}
\end{align}
Note that the $dr/dt$ equation is decoupled from the $d\theta/dt$ equation in the sense that it does not depend on $\theta$ and can thus be solved on its own. 
An equilibrium point $r_e$ for the $dr/dt$ equation corresponds to an equilibrium point in the full system if $r_e=0$, and to a limit cycle in the full system if $r_e > 0$.
We find System~\eqref{eq:BautinShearNF_polar} has an equilibrium point, $e_{1}$, that is located at the origin and exists for all parameter values.
Noisy motion around $e_1$ corresponds to an NS state.
Furthermore, System~\eqref{eq:BautinShearNF_polar}  has two limit cycles, $L_{+}$ and $L_{-}$, whose radii are given by
\begin{align}
\label{eq:LC_radius}
r_{+}=\sqrt{\frac{1}{2} \left( s + \sqrt{s^2 + 4\,\mu} \right)} 
\quad\text{and} \quad
r_{-}=\sqrt{\frac{1}{2} \left( s - \sqrt{s^2 + 4\,\mu} \right)},
\end{align}
respectively. $L_{+}$ exists for $\mu > -s^2/4$ if $s > 0$ and for $\mu > 0$ if $s < 0$, whereas $L_{-}$ exists for $-s^2/4 < \mu < 0$ and $ s > 0$.
Noisy motion around $L_{+}$ corresponds to an S state.

\begin{figure}[t]
\centering
\includegraphics[width=0.75\textwidth]{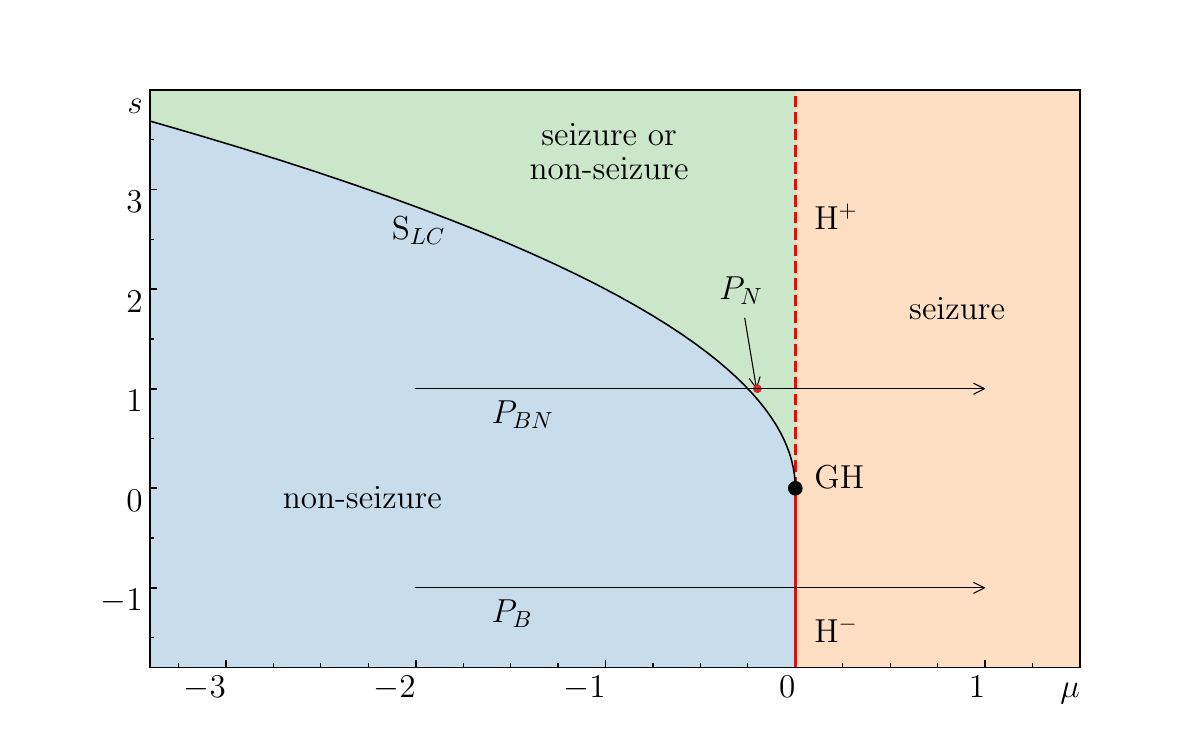}
\caption{
Two-parameter stability diagram for  System\,\eqref{eq:BautinShearNF_polar} in the $\left(\mu, s \right)$ parameter plane showing the regions of (blue) stable NS state, ({orange}) stable S state, and (green) bistability between the stable NS and S states.
S$_{LC}$ denotes the curve of saddle-node of limit cycle bifurcations, H$^{-}$ and H$^{+}$ denote the curves of supercritical and subcritical Hopf bifurcations, respectively, and GH denotes the Bautin (generalised Hopf) bifurcation point. 
The parameter paths, $P_{B}$ and $P_{BN}$, and the point $P_{N}$, denote the values of $\mu$ and $s$ used to generate BCTs,  BNCTs, and NCTs, respectively, in the model's output.
Note that the bifurcation curves do not depend on any other system parameters.
}
\label{fig:mu_s_bifplot}
\end{figure}

Linear stability analysis of equilibria for the 
$dr/dt$
equation reveals the following.
The equilibrium $e_1$ is stable for $\mu < 0$ and unstable for $\mu > 0$. 
The limit cycle $L_+$ is stable while  $L_-$ is unstable.
Bifurcations of $e_{1}$, $L_{-}$, and $L_{+}$ for changes in $\mu$ and $s$ are shown in a two-parameter bifurcation diagram in Fig.~\ref{fig:mu_s_bifplot}. For $s<0$ and as $\mu$ changes from negative to positive values, $e_{1}$ becomes unstable and $L_{+}$ is born in a supercritical Hopf bifurcation H$^{-}$ that occurs when crossing the solid part of the vertical red line at $\mu = 0$; see also the black solution branches in Fig.~\ref{fig:CTtypes}~(a).
For $s>0$ and as $\mu$ changes from positive to negative values, $e_{1}$ becomes stable and the unstable $L_{-}$ is born in a subcritical Hopf bifurcation H$^{+}$ that occurs when crossing the dashed part of the vertical red line at $\mu = 0$.
By decreasing $\mu$ further, $L_{+}$ and $L_{-}$ collide and disappear in a saddle-node of limit cycles bifurcation $S_{LC}$ that occurs when crossing the black solid curve given by $\mu = -s^2/4$ and $s>0$; see also the black solution branches in Fig.~\ref{fig:CTtypes}~(c).
Note that $e_{1}$, $L_{+}$, and $L_{-}$ coincide at a special Bautin (generalised Hopf) bifurcation point GH at $(\mu, s) = (0,0)$. 

As emphasised by the different coloured regions in Fig.~\ref{fig:mu_s_bifplot}, the bifurcations of $e_{1}$, $L_{-}$, and $L_{+}$ define three qualitatively different types of activity.
In the blue region, where $e_{1}$ is the only stable state, 
{only non-seizure-like activity can occur.}
In the orange region, where $L_{+}$ is the only stable state, 
{only seizure-like activity can occur.}
In the green region of bistability, where both  $e_{1}$ and $L_{+}$ are stable, either non-seizure-like or seizure-like activity can occur.

In Fig.\,\ref{fig:mu_s_bifplot} we also indicate two parameter paths and one special point. 
The paths $P_{B}$ and $P_{BN}$ are used to generate many BCTs and BNCTs from the NS to S state, respectively. This is done by repeatedly increasing $\mu$ from $-2$ to $1$ for $s=1$ and $s=-1$, respectively, according to
\begin{align}\label{eq:MuForcingFunc}
    \mu(t) = -2 + t/20, 
\end{align}
for $t \ge 0$.
The point $P_{N}$ at $s=1$ and $\mu = -0.22$ is used to generate a long time series containing many NCTs between the NS and S states.

In order to provide a direct comparison with the real seizure activity seen in voltage recordings such as Fig.~\ref{fig:RatData_Properties_C_}, we transform Eq.\,\eqref{eq:BautinShearNF_complex} into its Cartesian form by setting $z = x + i\,y$, where $x$ and $y$ are real-valued variables, and further augment it by including additive noise terms to generate the three different CTs mentioned above, 
\begin{align}\label{eq:BautinShearNF_cartesian_noise}
\begin{split}
   \frac{dx}{dt} 
    &= \gamma \, \left( \mu \, x + s \, \left( x^{2} + y^{2} \right) \, \left( x - \sigma \, y \right) - \omega \, y - \left( x^{2} + y^{2} \right)^{2} \, x \right) + \nu \, \eta_{x}(t), \\
    \frac{dy}{dt} 
    &= \gamma \, \left( \mu \, y + s \, \left( x^{2} + y^{2} \right) \, \left( y + \sigma \, x \right) + \omega \, x - \left( x^{2} + y^{2} \right)^{2} \, y \right) + \nu \, \eta_{y}(t).
\end{split}
\end{align}
Here, $\eta_{x}(t)$ and $\eta_{y}(t)$ are independent Gaussian random variables (white Gaussian noise with zero mean and unit variance) and $\nu$ is the noise level (standard deviation of the noise). 
To generate solutions for $x$ and $y$ at a given time $t$, we discretise System~\eqref{eq:BautinShearNF_cartesian_noise} according to the Euler-Maruyama method for a fixed time step $\delta_{t}$. 
In our numerical experiments, we set $\delta_t=0.001$, $x(0)=x_{0}=0.1$, and $y(0)=y_{0}=0.1$. This value of $\delta_t$ is chosen to prevent numerical errors from arising when $\sigma$ is relatively large.
The choice of initial condition $(x_0,y_0)$ ensures that all simulations of system~\eqref{eq:BautinShearNF_cartesian_noise} begin from the NS state.

Throughout this paper, we refer to system~\eqref{eq:BautinShearNF_cartesian_noise} as the `model', and to a time series of $x(t)$ as the `model's output'. 



\subsection*{M3: Tuning the model parameters to mimic voltage recordings \label{md:tuning_model}}

We begin by briefly describing three important characteristics of real seizure activity in the voltage recordings that we tune our mathematical model to mimic. These characteristics allow us to (i) use the same algorithm to detect CTs between the NS and S states in the model's output and in the voltage recordings, and (ii) train a machine learning classifier using the model's output and then apply it to the voltage recordings. 

\begin{figure}[t!]
    \centering
    \includegraphics[width=\linewidth]{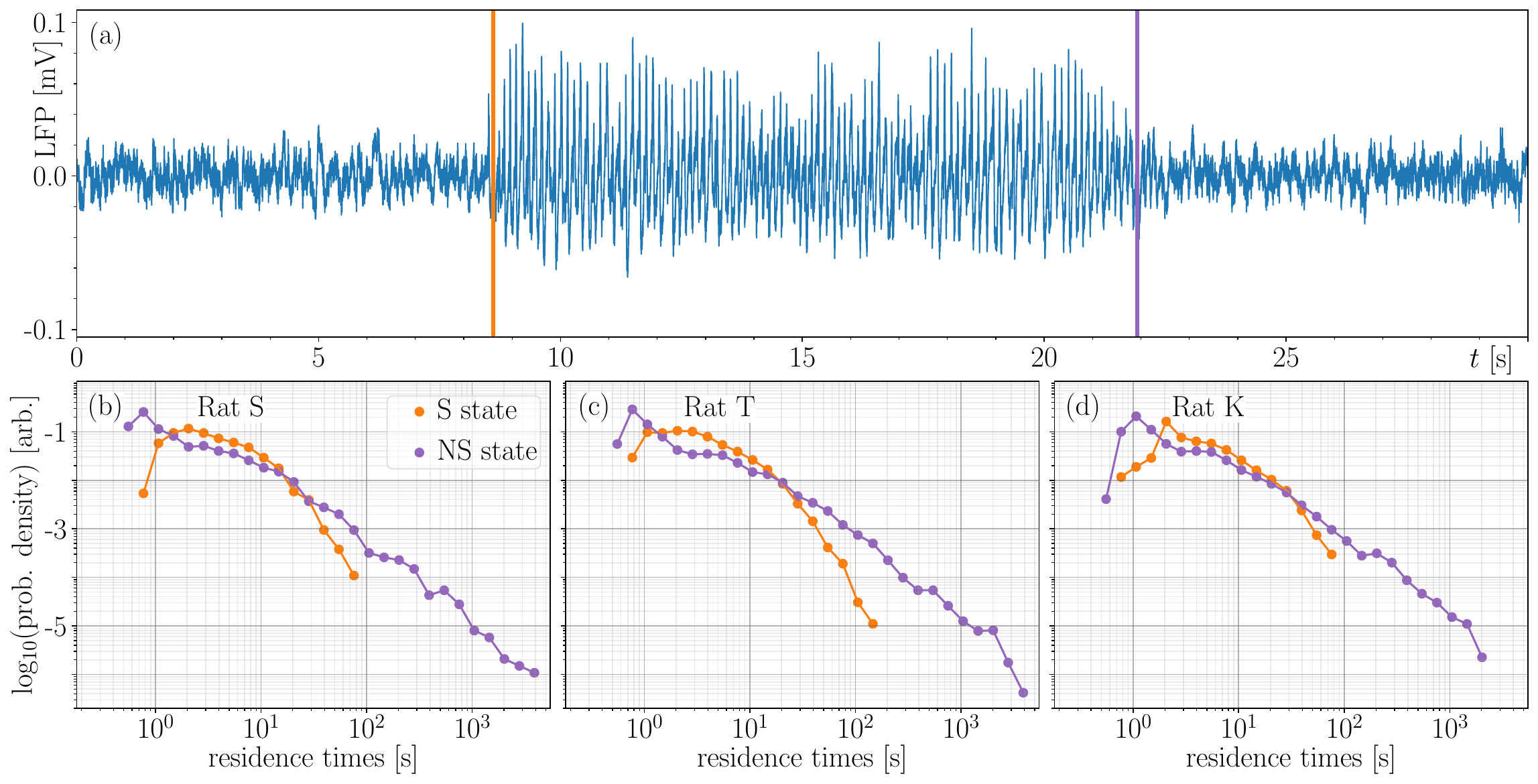}
    \caption{Characteristics of real seizure activity in GAERS. (a) An example of voltage recordings of real seizure activity, where the vertical lines indicate times of  electrical seizure (orange) onset and (purple) offset according expert annotations. (b)-(d) The probability density of residence times in the S and NS states of rat S, T, and K respectively.}
    \label{fig:RatData_Properties_C_}
\end{figure}

Figure~\ref{fig:RatData_Properties_C_}~(a) shows a typical example of real seizure activity in GAERS, taken from recording session `S5A' (rat S, day 5 of recording, A$^{th}$ recording session). This voltage recording was obtained for $t \geq 0$ with  time resolution of $0.001$~seconds via the procedure outlined in part \hyperref[md:CianSeizureAnnotate]{M1} of the Methods section. 
The vertical orange and purple lines in Fig.~\ref{fig:RatData_Properties_C_}~(a) indicate the electrical seizure onset and offset times, respectively, according to the expert. For the $i^{th}$ seizure, we denote the onset time by $t_{1}^{(i)} > 0$ and the offset time by  $t_{2}^{(i)} > t_{1}^{(i)}$. 
The same notation is used to denote times where our algorithm detects CTs between the NS and S states in our  mathematical model.
Note that the (orange)  electrical seizure onset is characterised by the sudden change in the voltage recordings from a small-amplitude and weakly-correlated oscillation to a larger-amplitude and highly-correlated oscillation. The opposite occurs for the (purple) electrical seizure offset. These sudden changes and the timescale of the oscillations are  two important characteristics that we tune our  mathematical model to mimic. 

The seizure in Fig.~\ref{fig:RatData_Properties_C_}~(a) is just one of several seizures that have been annotated by the expert in this recording session. Looking at all seizures in the session, there is significant variation in the {duration} of the different intervals of real seizure activity ($t_{2}^{(i)}-t_{1}^{(i)}$), which we refer to as {\em residence times in the S state,} and real non-seizure activity ($t_{1}^{(i+1)}-t_{2}^{(i)}$), which we refer to as {\em residence times in the NS state}. 
Figures~\ref{fig:RatData_Properties_C_}~(b)-(d) show the probability density of residence times in the S and NS states from `rat S', `rat T', and `rat K', using all the expert annotations of  electrical seizure onset and offset times. The information presented here is based on 621, 1593, and 1143 seizures from each respective rat. Note that the residence times in the S state are noticeably shorter than those in the NS state, and that residence times for both states have a similar range for different rats. The range of residence times is the third important characteristic we tune our mathematical model to mimic for NCTs. 

In summary, the three characteristics of real seizure activity that we want our mathematical model to mimic are: 
\begin{enumerate}[label=\roman*,align=CenterWithParen,labelindent=10pt,itemindent=1em,leftmargin=!]
    \item The amplitude and correlation of oscillation in the NS and S states.
    \item The range of residence times in the NS and S states. 
    \item Intrinsic time scales of the NS and S states. 
\end{enumerate}
When tuning the model's parameters we found that trade-offs emerge and we have to make choices on what characteristics to prioritise.
Optimising the bifurcation parameter, $\mu$, for characteristic (i) makes it difficult to find a noise level, $\nu$, to mimic characteristic (ii). 
Similarly, optimising the timescale parameter, $\gamma$, for characteristic (iii) makes it difficult to mimic characteristic (ii), and vice-versa.
We  find that characteristics (i) and (ii) are the most critical for training the machine learning classifier to accurately classify seizure generation mechanisms in time series. 
We therefore choose to prioritise characteristics (i) and (ii) and accept a trade-off in characteristic (iii).
The model parameters in Table~\ref{tab:ModelParams} have been chosen based on the above discussion.

\begin{table}[t!]
\centering
\def\arraystretch{1.2}
\begin{tabular}{|c||c|c|c|}
\hline
 & BCT & BNCT & NCT \\
\hline
\multirow{2}{10em}{\centering $\mu$: bifurcation parameter} & \multirow{2}{7em}{\centering $\mu(t)=-2+t/20$ for $t\in[0,60]$} & \multirow{2}{7em}{\centering $\mu(t)=-2+t/20$ for $t\in[0,60]$} & \multirow{2}{4em}{\centering $-0.22$} \\ & & & \\
\hline
\multirow{2}{10em}{\centering $s$: bifurcation parameter} & \multirow{2}{4em}{\centering $-1$} & \multirow{2}{4em}{\centering $1$} & \multirow{2}{4em}{\centering $1$} \\ & & & \\
\hline
\multirow{2}{10em}{\centering $\sigma$: shear} & \multirow{2}{4em}{\centering $-1$} & \multirow{2}{4em}{\centering $1$} & \multirow{2}{4em}{\centering $1$} \\ & & & \\
\hline
\multirow{2}{10em}{\centering $\nu$: noise level} & \multirow{2}{4em}{\centering $0.18$} & \multirow{2}{4em}{\centering $0.18$} & \multirow{2}{4em}{\centering $0.18$} \\ & & & \\
\hline
\multirow{2}{10em}{\centering $\omega$: oscillation freq.} & \multirow{2}{4em}{\centering $1.3$} & \multirow{2}{4em}{\centering $1.3$} & \multirow{2}{4em}{\centering $1.3$} \\ & & & \\
\hline 
\multirow{2}{10em}{\centering $\gamma$: time scale} & \multirow{2}{4em}{\centering $10$} & \multirow{2}{4em}{\centering $10$} & \multirow{2}{4em}{\centering $10$} \\ & & & \\
\hline
\end{tabular}
\caption{\label{tab:ModelParams} 
The parameter values chosen for our model in System~\eqref{eq:BautinShearNF_cartesian_noise} to generate different CT-types and mimic characteristics of seizure activity in voltage recordings.}
\end{table}

\begin{figure}[t]
    \centering
    \includegraphics[width=\linewidth]{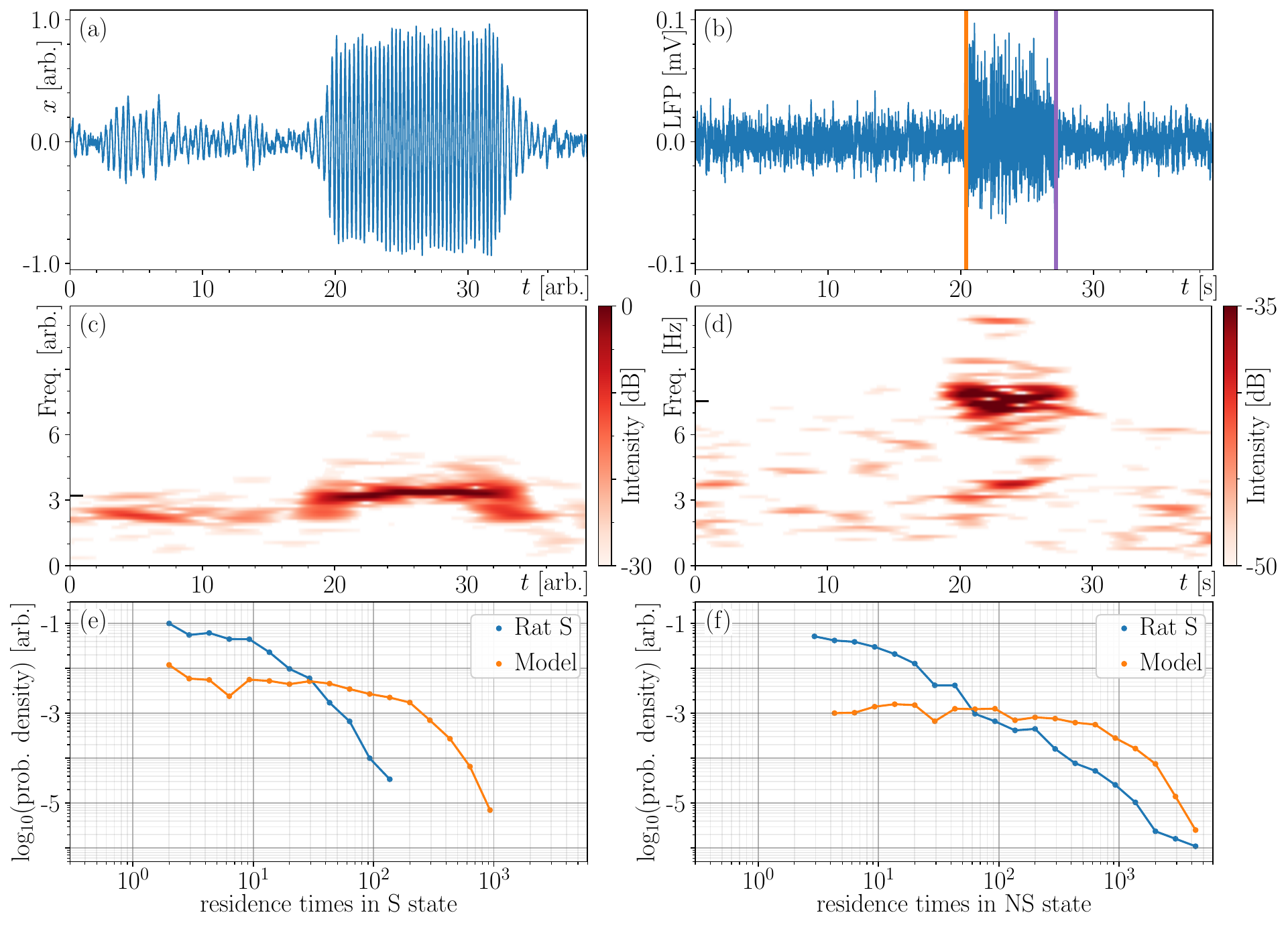}
    \caption{
    A comparison of the characteristics of (a) the model's output and (b) the actual voltage recordings for rat S.
    (c)-(d) show spectrograms of the time series from (a) and (b), respectively.
    (e)-(f) show the probability density of residence times in the S and NS state, respectively, for (blue) rat S and (orange) the model's output.
    }
    \label{fig:_Model_vs_Data_3characteristics_}
\end{figure}

\noindent
In Fig.~\ref{fig:_Model_vs_Data_3characteristics_} we show how the model's output compares to the voltage recordings when using the model parameters specified in Table~\ref{tab:ModelParams}. We discuss the contents of Fig.~\ref{fig:_Model_vs_Data_3characteristics_} below in terms of the three important characteristics.
\\
\textbf{Characteristic (i): }
The top row in Fig.~\ref{fig:_Model_vs_Data_3characteristics_} shows  that  the model's output (a) closely resembles real voltage recordings in GAERS (b). Characteristic (i) is reasonably well satisfied since the relative change in the amplitude of oscillation for a CT from the NS to S state in the model's output is similar to the seizure onset in the voltage recordings. The same can be said for the CT from the S to NS state and the seizure offset. 
Figures~\ref{fig:_Model_vs_Data_3characteristics_}~(a) and (b) also highlight one benefit of including shear, quantified by $\sigma$ in the model:  the model's output exhibits similar variations in its local maxima/minima to those observed in voltage recordings and these variations are beyond what can be achieved with noise alone. 
\\
\textbf{Characteristic (ii): }
In the bottom row, Fig.~\ref{fig:_Model_vs_Data_3characteristics_}~(e) shows the probability density of residence times in the S state for rat S (blue) and the model (orange). Figure~\ref{fig:_Model_vs_Data_3characteristics_}~(f) shows the same in the NS state. These residence times were obtained by applying our CT detection algorithm, introduced in the next part of the Methods section, to (i) all voltage recordings for rat S and (ii) a sufficiently long time series of the model's output where a similar number of CTs between the NS and S states occur, all of which were NCTs. 
Figures~\ref{fig:_Model_vs_Data_3characteristics_}~(e) and (f) show that the model exhibits a range of residence times in the NS and S state, which is the same or wider than those in the voltage recordings. 
However, the model tends to spend longer in the S state than the voltage recordings, resulting in some differences in the magnitude of the probability density. 
\\
\textbf{Characteristic (iii): }
In the middle row, Figs.~\ref{fig:_Model_vs_Data_3characteristics_}~(c) and (d) show spectrograms of the time series from Figs.~\ref{fig:_Model_vs_Data_3characteristics_}~(a) and (b), respectively. 
While in the S state, the model's output has a component at the same frequency as the voltage recording. However, the dominant component in the voltage recording is about twice that of the model's output.

\subsection*{M4: Constructing a CT detection algorithm \label{md:ConstructingAlgorithm}}



The time resolution of the voltage recordings determines the time step size $\delta$, i.e., the time interval between two consecutive points in the time series. The model's output is generated with the same time step size. In the case of voltage recordings in GAERS, we have $\delta=0.001$.

Additionally, we introduce six algorithm parameters: the ``on'' threshold $\alpha > 0$, the ``off'' threshold $0 <\beta < \alpha$, the size of the moving window $\tau_w > 0$, the time step size of the moving window $\delta \leq \Delta \leq \tau_w$, the minimum time duration of larger-amplitude oscillations, $\tau_{\text{S}} > \tau_{w}$, expressed as $\tau_{\text{S}} = n_{\text{S}}\Delta + \tau_w$, where $n_\text{S}$ is an integer, and the minimum time duration of lower-amplitude oscillations, $\tau_{\text{NS}} \geq \tau_{\text{S}}$, expressed similarly as $\tau_{\text{NS}} = n_{\text{NS}}\Delta + \tau_w$, where $n_{\text{NS}}$ is an integer. 
We list these parameters and specify their values in 
Table~\ref{tab:AlgorithmParameters}.
For convenience, we describe the algorithm below in terms of the model's output, $x(t)$. 
\\
{\bf The starting point:} 
We start in the NS state, where $|x(t)| < \beta$ for the time duration of at least $\tau_{\text{NS}}$.\\
{\bf The moving window:} 
When the system is in the NS state and $|x(t)|$ exceeds $\alpha$ at time $t = t_j$, the moving window  is activated and $|x(t)|$ is examined within consecutive windows of duration $\tau_w$  that are shifted in time by $\Delta$, starting with $[t_j,\, t_j + \tau_w]$, then $[t_j + \Delta,\, t_j  + \tau_w + \Delta]$, $[t_j + 2 \, \Delta,\, t_j  + \tau_w + 2 \, \Delta]$, and so on. 
The moving window is deactivated in two cases: (i) the system is in the NS state, the window is activated, but the algorithm does not detect a CT to the S state, and (ii) the system is in the S state and the algorithm detects a CT to the NS state. 
\\
{\bf Critical transitions:}
The algorithm detects a {\em CT from the NS to S state} at time $t=t_1$ if:
\begin{itemize}[labelindent=10pt,itemindent=1em,leftmargin=!]
    \item[(a1)\label{a1}]
    The system is in the NS state just before $t_1$.
    \item[(a2)\label{a2}]
    $|x(t)|$ exceeds $\alpha$ at time $t=t_1$,  i.e., $|x(t_1)|=\alpha$ and $|x(t_1 + \delta)| > \alpha$.
     \item[(a3)\label{a3}]
     Each of the $n_{\text{S}}$ consecutive positions of the moving window contains an $|x(t)| > \beta$.
\end{itemize}
The algorithm detects a {\em CT from the S to NS state} at time  $t=t_2$ if:
\begin{itemize}[labelindent=10pt,itemindent=1em,leftmargin=!]
    \item[(b1)\label{b1}]
    The system is in the S state just before $t_2$.
    \item[(b2)\label{b2}]
    $|x(t)|$ falls below $\beta$ at time $t=t_2$, i.e., $|x(t_2)|=\beta$ and $|x(t_2 + \delta)| < \beta$.
     \item[(b3)\label{b3}]
     Each of the $n_{\text{NS}}$ consecutive positions of the moving window contains no $|x(t)| \geq \alpha$.
\end{itemize}
In other words, the algorithm detects a CT from the NS to S state if $|x(t)|$ exceeds the upper threshold $\alpha$ and then continues to exceed the lower threshold $\beta$ frequently enough for a period of at least $\tau_{\text{S}}$. Similarly, the algorithm detects a CT from the S to NS state if $|x(t)|$ falls below the lower threshold $\beta$ and then does not exceed the upper threshold $\alpha$ for a period of at least $\tau_{\text{NS}}$. 
In Fig.~\ref{fig:Algorithm_in_use_} we provide an illustrative example of how the above algorithm detects CTs in the model's output from Fig.~\ref{fig:_Model_vs_Data_3characteristics_}~(a) using the algorithm parameters specified in Table~\ref{tab:AlgorithmParameters}. 
\\
\textbf{Almost-occurring critical transitions: }
The algorithm detects an {\em almost-occurring CT from the NS to S state} at time  $t=\tilde{t}_{1}$ if \hyperref[a1]{(a1)} and \hyperref[a2]{(a2)} are satisfied but \hyperref[a3]{(a3)} is not. 
Similarly, the algorithm detect an {\em almost-occurring CT from the S to NS state} at time  $t=\tilde{t}_{2}$ if \hyperref[b1]{(b1)} and \hyperref[b2]{(b2)} are satisfied but \hyperref[b3]{(b3)} is not.
\\
\textbf{Artefacts: }
We refer to Sec.~{S2} in [\citeonline{Flynn25_Supplement}] for details on how we adapt the above algorithm to account for artefacts that sometimes appear in the voltage recordings. 
\\
\textbf{Tuning algorithm parameters: }
We refer to Sec.~{S5} in [\citeonline{Flynn25_Supplement}] for details on how we tune the algorithm's parameters to maximise the agreement between the algorithm and the expert's annotations of electrical seizure onset and offset times in voltage recordings of real seizure activity. 
{We also refer to Flynn \textit{et al.}~\cite{flynn2026detecting} for a detailed ROC analysis of the algorithm's performance.}

\begin{table}[t]
\centering
\def\arraystretch{1.2}
\begin{tabular}{|c||c|c|}
\hline
Time series and detection algorithm parameters & Model's output & Voltage recordings \\
\hline
\multirow{2}{10em}{\centering $\delta$: time step size} & \multirow{2}{4em}{\centering 0.001} & \multirow{2}{6em}{\centering 0.001} \\ & & \\
\hline
\multirow{2}{10em}{\centering $\alpha$: on threshold} & \multirow{2}{4em}{\centering 0.55} & \multirow{2}{6em}{\centering [0.03, 0.1]} \\ & & \\
\hline 
\multirow{2}{10em}{\centering $\beta$: off threshold} & \multirow{2}{4em}{\centering 0.45} & \multirow{2}{6em}{\centering $\alpha$ - 0.01} \\ & & \\
\hline
\multirow{2}{14em}{\centering $\tau_{w}$: size of the moving window} & \multirow{2}{4em}{\centering 1} & \multirow{2}{6em}{\centering 1} \\ & & \\
\hline
\multirow{2}{18em}{\centering $\Delta$: time step size of the moving window} & \multirow{2}{4em}{\centering 0.001} & \multirow{2}{6em}{\centering 0.001} \\ & & \\
\hline
\multirow{2}{18em}{\centering $\tau_{\text{S}}$: minimum time duration of the S state} & \multirow{2}{4em}{\centering 2} & \multirow{2}{6em}{\centering 2}\\ & & \\
\hline
\multirow{2}{18em}{\centering $\tau_{\text{NS}}$: minimum time duration of the NS state} & \multirow{2}{4em}{\centering 5} & \multirow{2}{6em}{\centering 3} \\ & & \\
\hline
\end{tabular}
\caption{\label{tab:AlgorithmParameters} 
The values of the time step size and CT-detection algorithm parameters used for the model's output and voltage recordings. $\alpha$ is chosen differently for each voltage recording session according to the method outlined in Sec.~{S5} in [\citeonline{Flynn25_Supplement}].}
\end{table}


\subsection*{M5: Computing the time series properties and their slopes \label{md:TSPs}} 

\begin{figure}[t]
\centering
\includegraphics[width=0.88\textwidth]{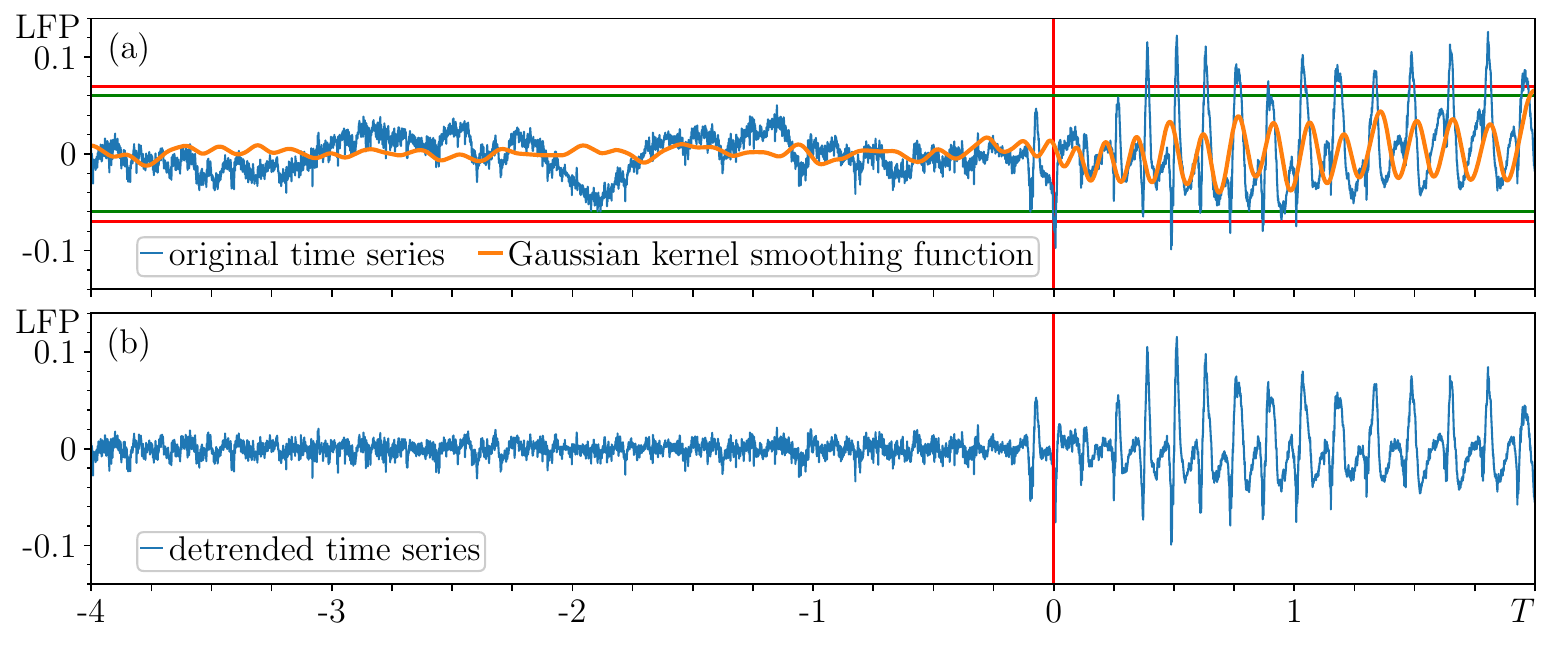}
\caption{Detrending voltage recordings using a Gaussian kernel smoothing function. (a) shows (in blue) an example of voltage recordings of real seizure activity before detrending and (in orange) the corresponding function fitted to these voltage recordings, and (b) shows the detrended voltage recordings. 
The vertical red solid line indicates when $T=0$. 
Horizontal lines indicate the values of (in red) $\alpha$ and (in green) $\beta$ used by the algorithm. The other algorithm parameters are specified in Table~\ref{tab:AlgorithmParameters}.}
\label{fig:_Compare_timeseries_detrendedimeseries_}
\end{figure}

\begin{figure}[t!]
\centering
\includegraphics[width=0.88\textwidth]{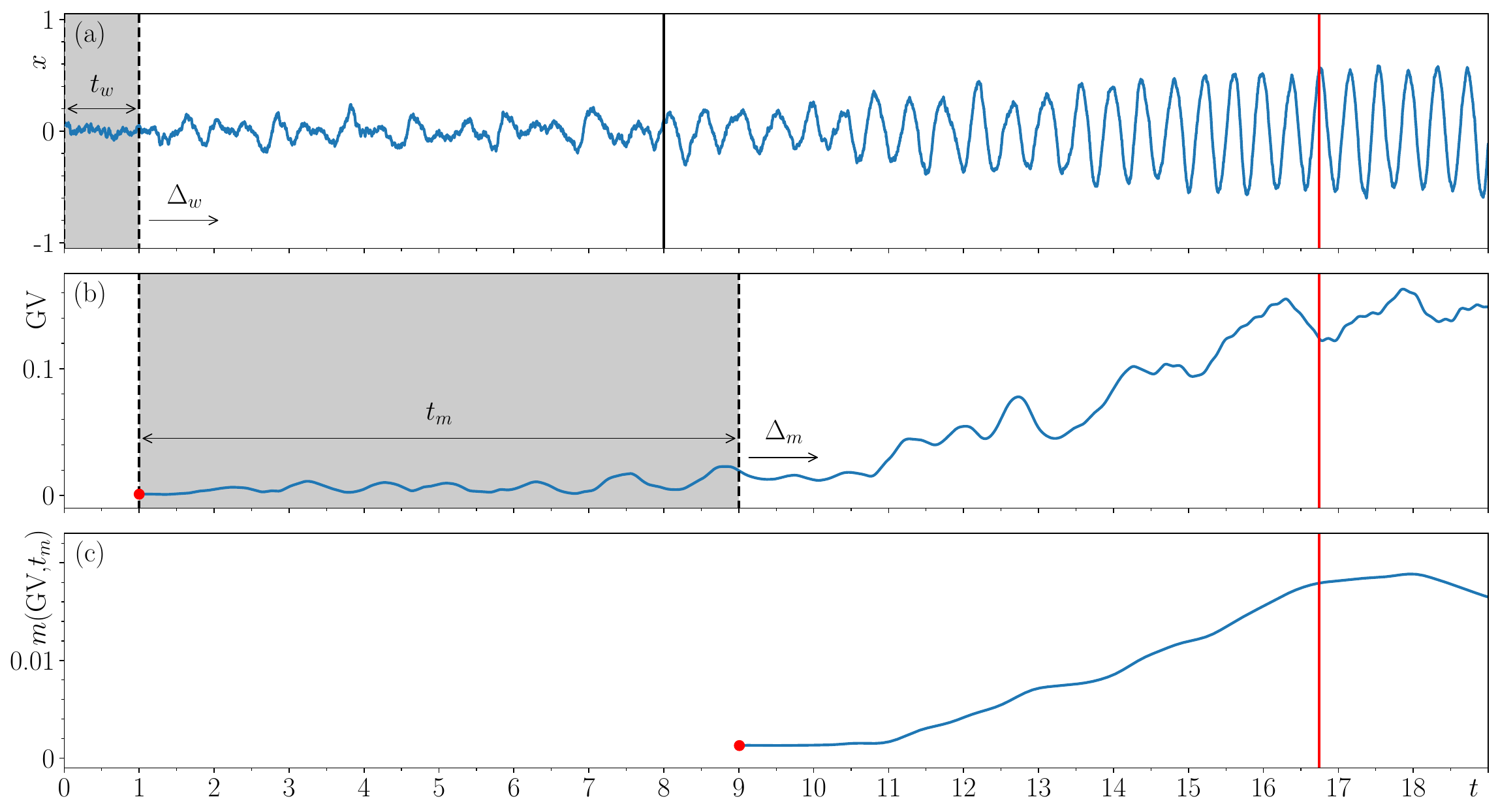}
\caption{Illustrating process for computing values of GV and $m$(GV, $t_m$) using moving windows of {duration} $t_{w}$ and $t_{m}$, which are advanced along the respective time series in steps of $\Delta_{w}$ and $\Delta_{m}$, {respectively.} 
The shaded region shows values of: (a) the model's output used to compute (red point in (b)) the first GV value and (b) GV used to compute (red point in (c)) the first $m$(GV, $t_m$) value in (c). 
The time series in (a) is an example of the model's output around a BCT from the NS to S state at $t=8$ ({time of bifurcation is indicated by vertical black line, time of detected CT is indicated by vertical red line}), generated using parameters specified in Table~\ref{tab:ModelParams} and plotted here for a new time $t=t'+38$, where $t'$ is the original time, to better illustrate the above process. $t_{w}=1$, $\Delta_{w}=0.001$, $t_{m}=8$, and $\Delta_{m}=0.1$. {The detection algorithm parameters are specified in Table~\ref{tab:AlgorithmParameters}.}}
\label{fig:_compution_of_GV_mGV_picture_}
\end{figure}

\begin{figure}[h!]
\centering
\includegraphics[width=0.88\textwidth]{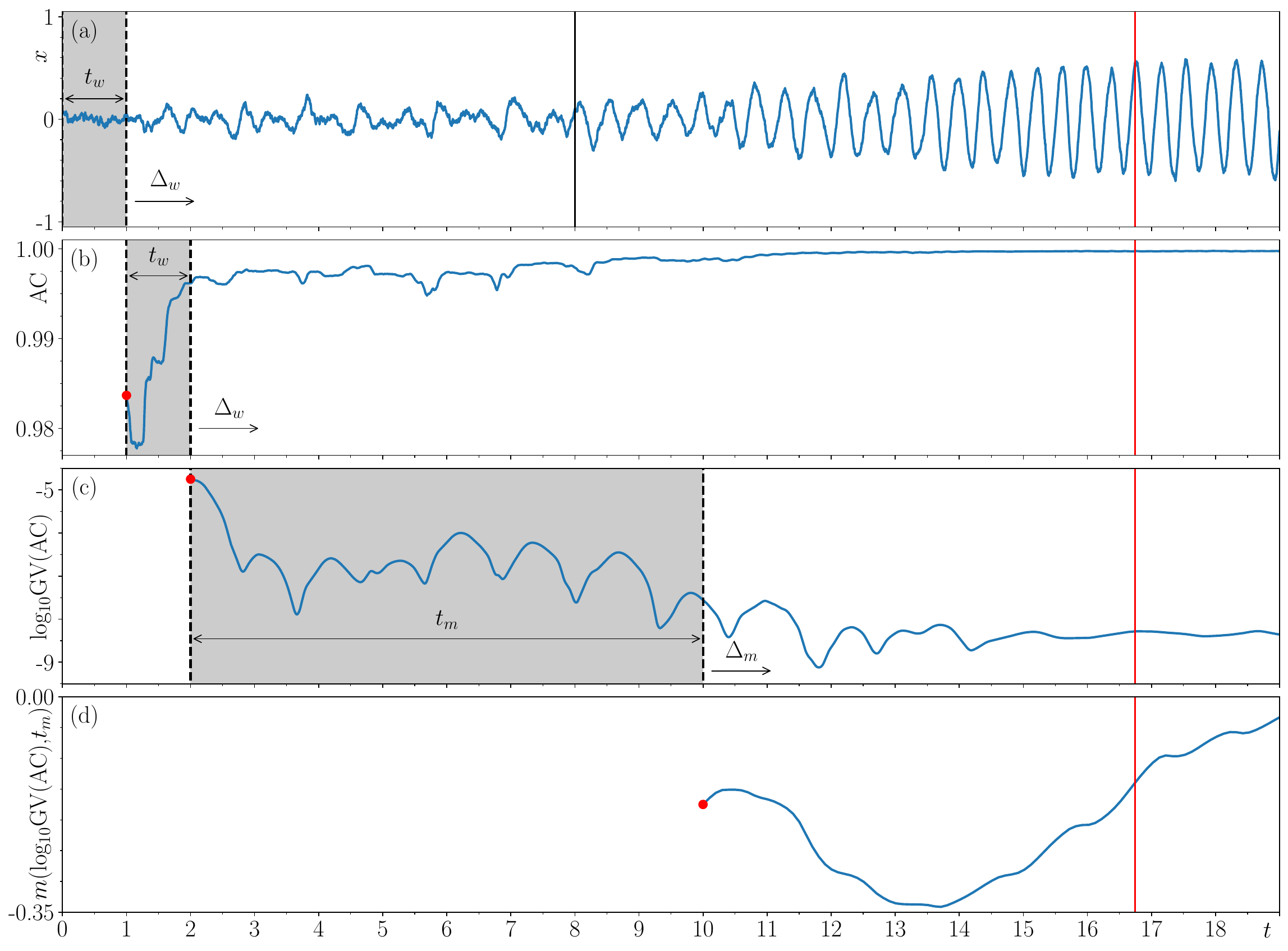}
\caption{Illustrating process for computing values of log$_{10}$GV(AC) and $m$(log$_{10}$GV(AC), $t_m$) using moving windows {duration} $t_{w}$ and $t_{m}$, which are advanced along the respective time series in steps of $\Delta_{w}$ and $\Delta_{m}$. 
The shaded regions indicate values of: (a) the model's output used to compute (red point in (b)) the first AV value, (b) AC used to compute (red point in (c)) the first log$_{10}$GV(AC) value, and (c) log$_{10}$GV(AC) used to compute (red point in (d)) the first $m$(log$_{10}$GV(AC), $t_m$) value. 
The time series in (a) {and algorithm parameters are} the same as {used in} Fig.~\ref{fig:_compution_of_GV_mGV_picture_}~(a) (see caption for details). $t_{w}=1$, $\Delta_{w}=0.001$, $t_{m}=8$, and $\Delta_{m}=0.1$.}
\label{fig:_compution_of_logGVAC_mlogGVAC_picture_}
\end{figure}

{In this section, we summarise the procedures used to compute different properties in time series of both the model's output and the voltage recordings. These properties are then used 
to
classify different CT-types. }
\\
\textbf{Detrending the time series: } Time series of voltage recordings are detrended using the Gaussian filtering approach mentioned in Dakos {et al.}\cite{dakos2008detrend}. This involves fitting a Gaussian kernel smoothing function to a given time series and then subtracting the fit from the time series to obtain the detrended time series. As further detailed in Lenton \textit{et al.}\cite{lenton2012detrend}, this process involves choosing a suitable bandwidth (a fitting parameter) for the kernel. The bandwidth determines the degree of smoothening and is typically chosen to neither over-fit the data nor filter out low frequencies from the time series, in our case we choose a bandwidth of $30$. In our experiments we use the `gaussian\_filter1d' function from the `scipy.ndimage' Python library to detrend the time series. See Fig.~\ref{fig:_Compare_timeseries_detrendedimeseries_} for an example of a detrended time series of voltage recordings. 
In our experiments we find that, whether the voltage recordings are detrended or not, there are minor changes in the main results presented in Figs.~\ref{fig:RatSTK_Tm_vs_classfiltCTs_filtdetCTs_Tm_16_8_}~(a)-(c), these changes are on the order of 2-3$\%$ in terms of the percentage of CTs classified as BCTs, BNCTs, and NCTs specified in Table~\ref{tab:FinalResults}. 
\\
\textbf{Variance: } We compute the Gaussian rolling variance at time $t \geq t_{w}$ using data points of the time series within the interval, $[t-t_{w},\, t]$, which we call `the window' and $t_{w}>0$ is the `window {duration}'. 
The variance is calculated as the weighted variance of a Gaussian distribution of the data points centered at $t-t_{w}/2$ with mean equal to zero. 
This involves assigning weights to data points within the window and these weights are determined from the probability density function of a Gaussian distribution centered at the middle of the window with standard deviation taken as $t_{w}/6$. 
Based on this procedure we define GV$(t) : \mathbb{R} \to \mathbb{R}$ as the Gaussian variance of the time series at $t \geq t_{w}$. 
\\
\textbf{Logarithm of the variance: } Considering the results presented by Milanowski and Suffczynski \cite{milanowski2016seizures}, we also find it useful to record the base-10 logarithm of the GV in each of the moving windows described above. This allows us to more closely examine the smaller changes in the variance. We define log$_{10}$GV$(t): \mathbb{R} \to \mathbb{R}$ as the base-10 logarithm of the Gaussian variance of the time series at $t \geq t_{w}$. 
\\
\textbf{Autocorrelation: } We compute the autocorrelation at $t \geq t_{w}$ using a rolling window technique and data points of the time series within the same interval as above, $[t-t_{w}, t]$. More specifically, we compute the lag-1 autocorrelation at a given time $t$ by computing the Pearson correlation between data points in the following two intervals, $[t-t_{w}+t_{a}, t]$ and $[t-t_{w}, t-t_{a}]$ where $t_{a}>0$ is the time between two successive points in the time series. 
Based on this procedure we define AC$(t) : \mathbb{R} \to \mathbb{R}$ as the lag-1 autocorrelation of the time series at $t \geq t_{w}$. 
\\
\textbf{Variance of the autocorrelation: } Prior to a BCT or BNCT from the NS to S state we observe a decrease in the fluctuations of AC$(t)$ and, after the CT, AC$(t)$ remains relatively constant and close to $1$. 
To more closely examine this behaviour, we compute the base-10 logarithm of the Gaussian variance of the time series obtained for AC$(t)$ using the same procedure outlined above. 
From this we define log$_{10}$GV(AC)$(t): \mathbb{R} \to \mathbb{R}$ as the base-10 logarithm of the Gaussian variance of the autocorrelation of a given time series at $t \geq 2~t_{w}$ (since AC$(t)$ is defined for $t \geq t_{w}$, the resulting log$_{10}$GV(AC)$(t)$ is defined for $t \geq 2~t_{w}$).
\\
\\
Each window of {duration} $t_{w}$ is shifted forward in time by a factor $\Delta_{w} > 0$ to obtain time series for each of the TSPs. We set $\Delta_{w} = \delta = 0.001$, the time between two successive points in the original time series. 
From several trial runs we find that setting $t_{w}=1$ provides relatively robust estimates of the variance and autocorrelation of the time series. 
For convenience we refer to GV$(t)$ as GV, and similarly for log$_{10}$GV, AC, and log$_{10}$GV(AC). We collectively refer to these four quantities as the time series properties (TSPs). 
\\
\\
\textbf{Slopes: } We also analyse how the slopes of the TSPs change over time by fitting the function, $f(x) = m\,x + c$, to the values of, for example, GV within the interval $[t-t_{m}, t]$, where $t_{m}$ is the `slope {duration}'. We use the `curve\_fit()' function in the SciPy Python library to perform a least squares fit of this function to the data in {the time interval} $[t-t_{m}, t]$ via the Levenberg–Marquardt algorithm. Based on this procedure we define $m(\text{GV}, t_{m}) : \mathbb{R} \times \mathbb{R} \to \mathbb{R}$ as the slope of the Gaussian variance of a given time series at $t \geq t_{w} + t_{m}$. The same procedure is applied to all TSPs. Note, $m(\text{log$_{10}$GV(AC)},~t_{m})$ is defined for $t \geq 2~t_{w} + t_{m}$. We collectively refer to these four quantities as `the slopes of the TSPs for a given $t_{m}$', denoted by $m$(TSPs, $t_m$). 
Each window of {duration} $t_{m}$ is shifted forward in time by a factor $\Delta_{m} > 0$ to obtain time series for each of the $m$(TSPs, $t_m$). We set $\Delta_{m} = 100\Delta_{w}$ to reduce computational time.
\\
\\
In Figs.~\ref{fig:_compution_of_GV_mGV_picture_} and \ref{fig:_compution_of_logGVAC_mlogGVAC_picture_} we illustrate the different time intervals over which TSPs and $m$(TSPs, $t_m$) are computed using an example of the model's output nearby a BCT from the NS to S state. This helps clarify (i) how values of TSPs and $m$(TSPs, $t_m$) are obtained from moving windows and (ii) the correspondence in time between TSPs, $m$(TSPs, $t_m$), and the original time series from which they are obtained. More specifically, Fig.~\ref{fig:_compution_of_GV_mGV_picture_} illustrates the above for GV and $m(\text{GV},~t_{m})$, the same procedure is used to compute AC, log$_{10}$GV, and their corresponding slopes. 
Fig.~\ref{fig:_compution_of_logGVAC_mlogGVAC_picture_} illustrates the above for log$_{10}$GV(AC) and $m(\text{log$_{10}$GV(AC)},~t_{m})$.
In both Figs.~\ref{fig:_compution_of_GV_mGV_picture_} and \ref{fig:_compution_of_logGVAC_mlogGVAC_picture_}, $t_{w}=1$ and $t_{m}=8$.

\begin{table}[t]
\centering
\def\arraystretch{1.2}
\begin{tabular}{|c||c|c|c|c|}
\hline
 & \multirow{2}{3.9em}{\centering Gaussian variance} & \multirow{2}{8.7em}{\centering base-10 logarithm of the Gaussian variance} & \multirow{2}{5.8em}{\centering lag-1 autocorrelation} & \multirow{2}{14.2em}{\centering base-10 logarithm of the Gaussian variance of the lag-1 autocorrelation} \\ & & & & \\
\hline
\multirow{2}{9.8em}{\centering Time series properties (TSPs)} & \multirow{2}{3.9em}{\centering GV} & \multirow{2}{8.7em}{\centering log$_{10}$GV} & \multirow{2}{5.8em}{\centering AC} & \multirow{2}{14.2em}{\centering log$_{10}$GV(AC)} \\ & & & & \\
\hline
\multirow{2}{9.8em}{\centering Slopes of time series properties ($m(\text{TSPs},~t_{m})$)} & \multirow{2}{3.9em}{\centering $m(\text{GV},~t_{m})$} & \multirow{2}{8.7em}{\centering $m(\text{log$_{10}$GV},~t_{m})$} & \multirow{2}{5.8em}{\centering $m(\text{AC},~t_{m})$} & \multirow{2}{14.2em}{\centering $m(\text{log$_{10}$GV(AC)},~t_{m})$} \\ & & & & \\
\hline
\end{tabular}
\caption{\label{tab:TimeSeriesProperties} Quantities used as features to classify seizure generation mechanisms in a noisy time series. These quantities are dependent on time, i.e., $\text{GV} = \text{GV}(t)$, and are computed using data points of a given time series within the interval $[t-t_{w}, t]$ where $t_{w}>0$ is the `window {duration}'. $t_{m}$ denotes the slope {duration}, the duration of time over which the slope is calculated.}
\end{table}
%



\subsection*{M6: SVM training and testing \label{md:SVMtrain}}

We use the `sklearn.svm' Python library to construct a linear support vector machine (SVM) that classifies CTs from the NS to S state according to the CT-type where type = BCT, BNCT, or NCT. 
The SVM classifier is obtained by finding an optimal 
{hyperplane}
that maximises the distance between each class, i.e., CT-type, in an $F$-dimensional space where $F \in \mathbb{N}^{+}$ is the number of `features'. 
\\
\textbf{Preparing the data: } We use values of the TSPs and/or $m$(TSPs, $t_m$) specified in Table~\ref{tab:SVM_TSPs} at time $t$ as features. $S \in \mathbb{N}$ samples of each feature are used to train the SVM. 
Since (i) some features are orders of magnitude larger than others and (ii) we would like to construct a classifier which treats features with equal importance, each of the features are scaled to be of zero mean and unit variance using the `preprocessing.scale' function from the `sklearn' Python library. 
The data used to construct the SVM classifier is structured as follows, the `feature matrix' is defined as $X \in \mathbb{R}^{S \times F}$, and the `class vector' is defined as $y \in \mathbb{R}^{S}$. The information contained in $X$ and $y$ define our data set. 
\\
\textbf{Training and testing: } We split the above data set into a training and testing data set whereby $70\%$ of the data is used for training and $30\%$ for testing (a convention commonly followed when training and testing SVMs.)
The SVM classifier is constructed based on the training data set and the accuracy of the SVM classifier is determined based on the fraction of correctly classified samples in the testing data set. 
\\
\textbf{Feature importance: } We use the `permutation\_importance' function from the `sklearn.inspection' Python library to determine what features have the greatest impact on the accuracy. 
This function constructs different random permutations of a given feature (i.e., reorders elements in a column of $X$ reserved for testing while the corresponding $y$ remains unchanged) and returns the resulting decrease in accuracy; the greater the decrease the more dependent the SVM classifier's accuracy is on the feature. 
This process is repeated 100 times for each feature to obtain an estimate of the `mean permutation importance' (MPI) where $-1 \leq \text{MPI} \leq 1$. We also obtain the corresponding standard deviations from this process. 
Therefore, if, for a given feature, 
\begin{itemize}
    \item  $\text{MPI} > 0$, the SVM classifier's accuracy depends on this feature and the greater the $\text{MPI}$ the more important the feature.
    \item $\text{MPI} \approx 0$, the SVM classifier's accuracy does not depend on this feature. 
    \item $\text{MPI} < 0$, the SVM classifier's accuracy increases when using the permuted data set as opposed to the original data set.
\end{itemize}  

\begin{table}[t]
\centering
\def\arraystretch{1.2}
\begin{tabular}{|c||l|}
\hline
 & Features used \\
\hline
SVM type-1 & TSPs: GV, log$_{10}$GV, AC, and log$_{10}$GV(AC). \\
\hline
SVM type-2 & $m$(TSPs, $t_m$) for a given $t_{m}$: $m(\text{GV},~t_{m})$, $m(\text{log$_{10}$GV},~t_{m})$, $m(\text{AC},~t_{m})$, and $m(\text{log$_{10}$GV(AC)},~t_{m})$. \\
\hline
SVM type-3 & All TSPs and $m$(TSPs, $t_m$) for a given $t_{m}$.\\
\hline
\end{tabular}
\caption{\label{tab:SVM_TSPs} Features used to train and test different types of Support Vector Machines (SVMs). The time series properties (TSPs) and their slopes, $m$(TSPs, $t_m$), are specified in Table~\ref{tab:TimeSeriesProperties}. $t_{m}$ denotes the duration of time over which the slope is calculated.}
\end{table}




\section*{Authorship and Contributorship Statement}

All authors have made substantial intellectual contributions to the study conception, execution, and design of the work. All authors have read and approved the final manuscript.  In addition, the following contributions occurred: 
A.F. and S.W. conceived the experiments, A.F. conducted the experiments, A.F. and S.W. analysed the results, C.Mc~C., F.D., and V.C. provided the voltage recordings, C.Mc~C., K.L., and W.P.M. provided valuable discussions. A.F., C.Mc~C., K.L., W.P.M., and S.W. reviewed the manuscript. 

\section*{Conflicts of Interest}

The authors declare there are no conflicts of interest.

\section*{Data and Code Availability Statement}

The data and code used to conduct our experiments are available upon reasonable request.

\section*{Ethics Statement}

The experiments taken to obtain the voltage recordings used in this paper are outlined in McCafferty \textit{et al.}~\cite{mccafferty2018cortical}. As stated in \citeonline{mccafferty2018cortical}, these experiments were approved under the UK Animals (Scientific Procedures) Act 1986 by Cardiff University Animal Welfare and Ethics Committee and were conducted in accordance with recommendations for experimental work in epilepsy at the time when \citeonline{mccafferty2018cortical} was published.

\section*{Funding}

This publication has emanated from research conducted with the financial support of Taighde \'{E}ireann – Research Ireland under grant number [19/FFP/6782].  



\renewcommand{\thefigure}{S-\arabic{figure}}

\setcounter{figure}{0}

\newpage
\section*{Supplementary information\label{si:sec}}

\noindent The supplementary materials provide additional supporting content that complements and clarifies points made in the manuscript that are not required for understanding the main results. 


\renewcommand{\thefigure}{S-\arabic{figure}}
\setcounter{figure}{0}


\section*{S1: Applying the algorithm to detect seizure-like activity in the model's output \label{si:Alg_model}}

\begin{figure}[t]
\centering
\includegraphics[scale=0.325]{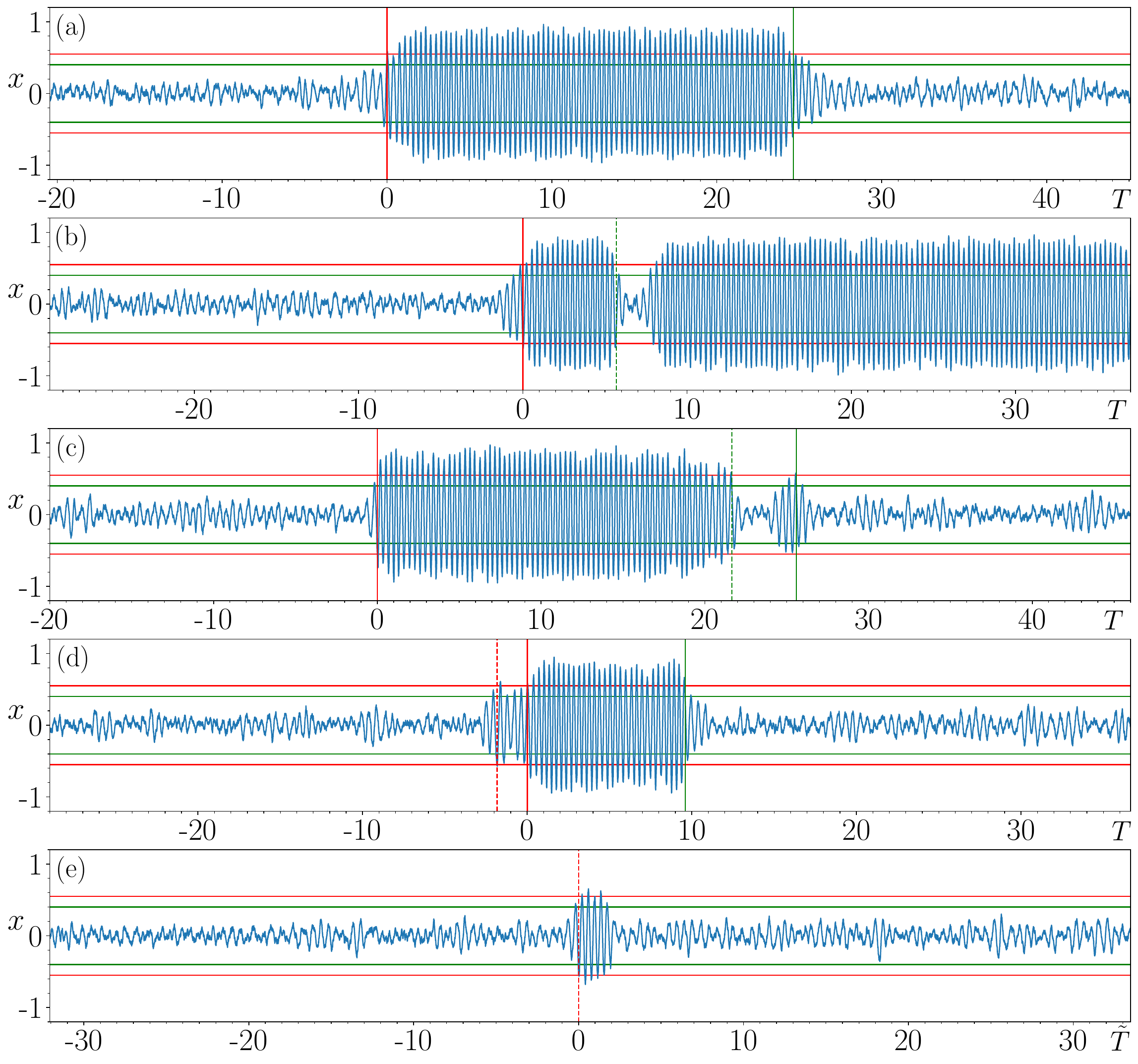}
\caption{Examples of NCTs in the model's output (the variable $x$ from System~(5) in \hyperref[seref:MainText]{[1]}) detected by the algorithm. The model's output is plotted versus $T = t - t_{1}$ for chosen values of $t_{1}$ in (a)-(d), and versus $\tilde{T} = t - \tilde{t}_{1}$ for a chosen $\tilde{t}_{1}$ in (e). 
Vertical lines indicate CTs from the (solid red) NS to S state and (solid green) S to NS state, and almost-occurring CTs from the (dashed red) NS to S state and (dashed green) S to NS state.
The model parameters are specified in Table~3 in \hyperref[seref:MainText]{[1]}. 
The algorithm parameters are specified in Table~4 in \hyperref[seref:MainText]{[1]}, the horizontal lines indicate the corresponding thresholds of (red) $\alpha$ and (green) $\beta$.}
\label{fig:SeizureDetectionExample1}
\end{figure}

In this section we illustrate some common outcomes that arise when using our algorithm to detect critical transitions (CTs) between non-seizure (NS) and seizure (S) states in noisy time series of our mathematical model's output (the variable $x$ from System~(5) in \hyperref[seref:MainText]{[1]}) for each CT-type it can produce, namely bifurcation-induced CTs (BCTs), noise-induced CT (NCTs), and bifurcation/noise-induced CTs (BNCTs). 
\\
\noindent
\textbf{NCTs:} We first integrate the model (System~(5) in \hyperref[seref:MainText]{[1]}) from $t=0$ to $t=200000$ using the procedure specified in part {M2} of the Methods section in \hyperref[seref:MainText]{[1]} with model parameters specified in Table~3 in \hyperref[seref:MainText]{[1]}.
We then apply our CT detection algorithm (see part {M4} of the Methods section in \hyperref[seref:MainText]{[1]}) to the model's output using the algorithm parameter values specified in Table~4 in \hyperref[seref:MainText]{[1]}. 
Bearing in mind that the algorithm detects a given CT from the NS (S) to S (NS) state at time $t=t_{1}$ ($t=t_{2}$), and an almost-occurring CT from the NS (S) to S (NS) state at time $t=\tilde{t}_{1}$ ($t=\tilde{t}_{2}$), the algorithm detected (i) 193 CTs between the NS and S states (i.e. 193 pairs of $t_{1}$ and $t_{2}$ values) and (ii) 86 almost-occurring CTs (i.e., a set of 86 values consisting of different $\tilde{t}_{1}$ and $\tilde{t}_{2}$ values), examples of these are shown in Fig.~\ref{fig:SeizureDetectionExample1}. 
For our convenience, in Figs.~\ref{fig:SeizureDetectionExample1}~(a)-(d) we plot $x$ versus $T=t-t_{1}$ for the chosen values of $t_{1}$, and in Fig.~\ref{fig:SeizureDetectionExample1}~(e) we plot $x$ versus $\tilde{T}=t-\tilde{t}_{1}$ for the chosen value of $\tilde{t}_{1}$. 
Solid vertical lines indicate (in red) when $T=0$ and (in green) when $T=t_{2}-t_{1}$. Similarly, dashed vertical lines indicate (in red) $\tilde{T}=0$ and (in green) $\tilde{T}=\tilde{t}_{2}-\tilde{t}_{1}$. 
The horizontal lines indicate the values of (in red) $\alpha$ and (in green) $\beta$ used by the algorithm.

\begin{figure}[t]
\centering
\includegraphics[scale=0.35]{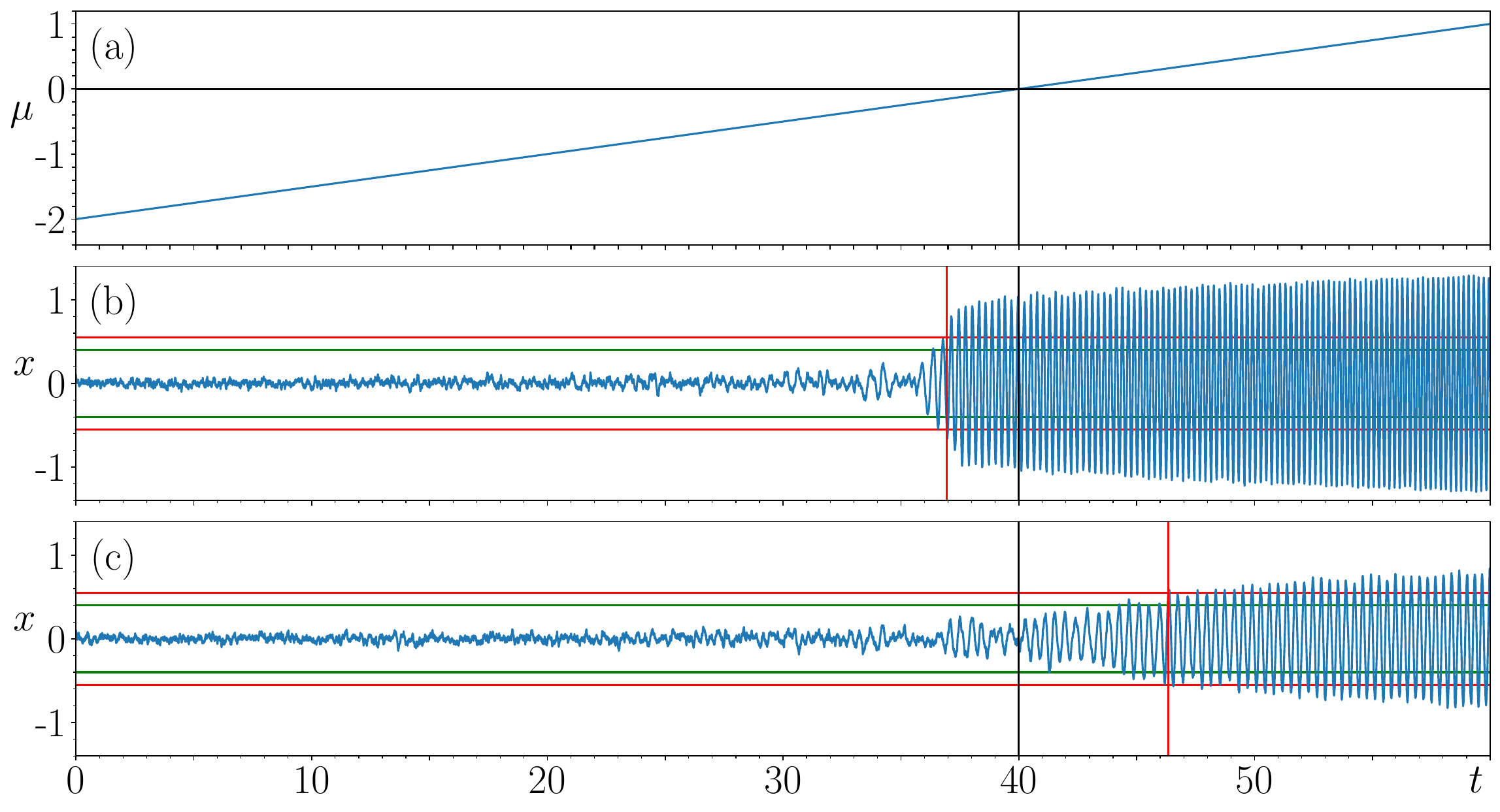}
\caption{Examples of BNCTs and BCTs from the NS to S state in the model's output (the variable $x$ from System~(5) in \hyperref[seref:MainText]{[1]}) that were detected by the algorithm. The evolution of $\mu$ versus $t$ is shown in (a), an example of a BNCT is shown in (b) and a BCT is shown in (c). Black vertical lines at $t=40$ indicate when the respective bifurcations occur. Red vertical lines indicate CTs from the NS to S state ($t_{1}$). The model parameters are specified in Table~3 in \hyperref[seref:MainText]{[1]}. The algorithm parameters are specified in Table~4 in \hyperref[seref:MainText]{[1]}, the horizontal lines indicate the thresholds of (in red) $\alpha$ and (in green) $\beta$.}
\label{fig:mu_sub_sup_vs_t_CTplots}
\end{figure}

Figure~\ref{fig:SeizureDetectionExample1}~(a) shows a clear-cut example of CTs between the NS and S states detected by the algorithm. 
Figure~\ref{fig:SeizureDetectionExample1}~(b) illustrates a common type of almost-occurring CT where, with reference to the conditions specified in part {M4} of the Methods section in \hyperref[seref:MainText]{[1]}, {(b1)} and {(b2)} are satisfied but {(b3)} is not.  
Another almost-occurring CT is shown in Fig.~\ref{fig:SeizureDetectionExample1}~(c), however, even though $x$ only exceeds $\alpha$ once after briefly falling below $\beta$ at the dashed vertical green line, the algorithm still considers the model to be in the S state until {(b1)}-{(b3)} are satisfied. 
Figure~\ref{fig:SeizureDetectionExample1}~(d) shows that the algorithm does not treat the opposite scenario to (b) in the same way. Regardless of how many almost-occurring CTs are detected, the algorithm still considers the model to be in the NS state until {(a1)}-{(a3)} are satisfied. 
The example in Fig.~\ref{fig:SeizureDetectionExample1}~(e) illustrates the most common almost-occurring CT from the NS to S state; {(a1)} and {(a2)} are satisfied but {(a3)} is not. 
\\
\noindent
\textbf{BCTs and BNCTs:}
We integrate the model from $t=0$ to $t=60$ with the parameters specified in Table~3 in \hyperref[seref:MainText]{[1]}, we set $s=1$ ($s=-1$) so that a BNCT (BCT) from the NS to S state can occur. Fig.~\ref{fig:mu_sub_sup_vs_t_CTplots}~(a) shows that $\mu$ increases linearly from $-2$ to $1$ according to Eq.~(4) in \hyperref[seref:MainText]{[1]}. 
We apply our algorithm to the model's output using the parameters values specified in Table~4 in \hyperref[seref:MainText]{[1]}. 
In Figs.~\ref{fig:mu_sub_sup_vs_t_CTplots}~(b) and (c) we illustrate an example of a BNCT and BCT detected by the algorithm. 
The vertical lines indicate (in black) when the bifurcation takes place and (in red) when the algorithm detects a CT from the NS to S state, i.e., the value of $t_{1}$. 
Figure~\ref{fig:mu_sub_sup_vs_t_CTplots}~(b) shows the algorithm may detect a BNCT before the bifurcation takes place. The opposite is shown in Fig.~\ref{fig:mu_sub_sup_vs_t_CTplots}~(c) for the BCT. We study this difference between values of $t_{1}$ in greater detail in Fig.~\ref{fig:probdens_CTtime_SuperSub_sigma_0_1_2_} in \hyperref[si:shear_t1times_bif]{S4}.

\section*{S2: Applying the algorithm to detect real seizure activity in the voltage recordings \label{si:Alg_data}}

\begin{figure}[t]
    \centering
    \includegraphics[width=0.7\linewidth]{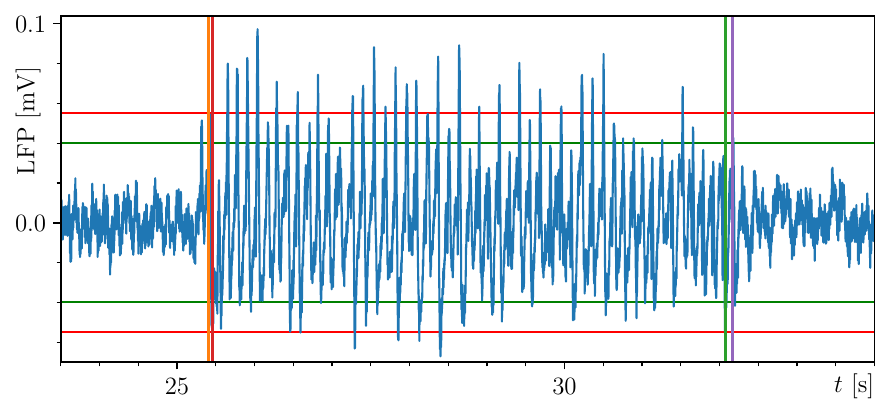}
    \caption{Expert annotations vs. Algorithm. Vertical lines indicate values of $t_{1}$ and $t_{2}$ according to the (orange and purple) expert's annotations and (red and green) algorithm for the same portion of the voltage recordings shown in Fig.~10~(b) in \hyperref[seref:MainText]{[1]}. Algorithm parameters used: $\alpha=0.055$, $\beta=0.04$, $\tau_{\text{NS}}=5$, $\tau_{\text{S}}=2$, $\tau_{w}=1$, and $\Delta = 0.001$, horizontal lines indicate the thresholds of (red) $\alpha$ and (green) $\beta$.}
    \label{fig:Compare_ModelOutput_to_RatData_detection}
\end{figure}

\begin{figure}[t]
    \centering
    \includegraphics[width=0.7\linewidth]{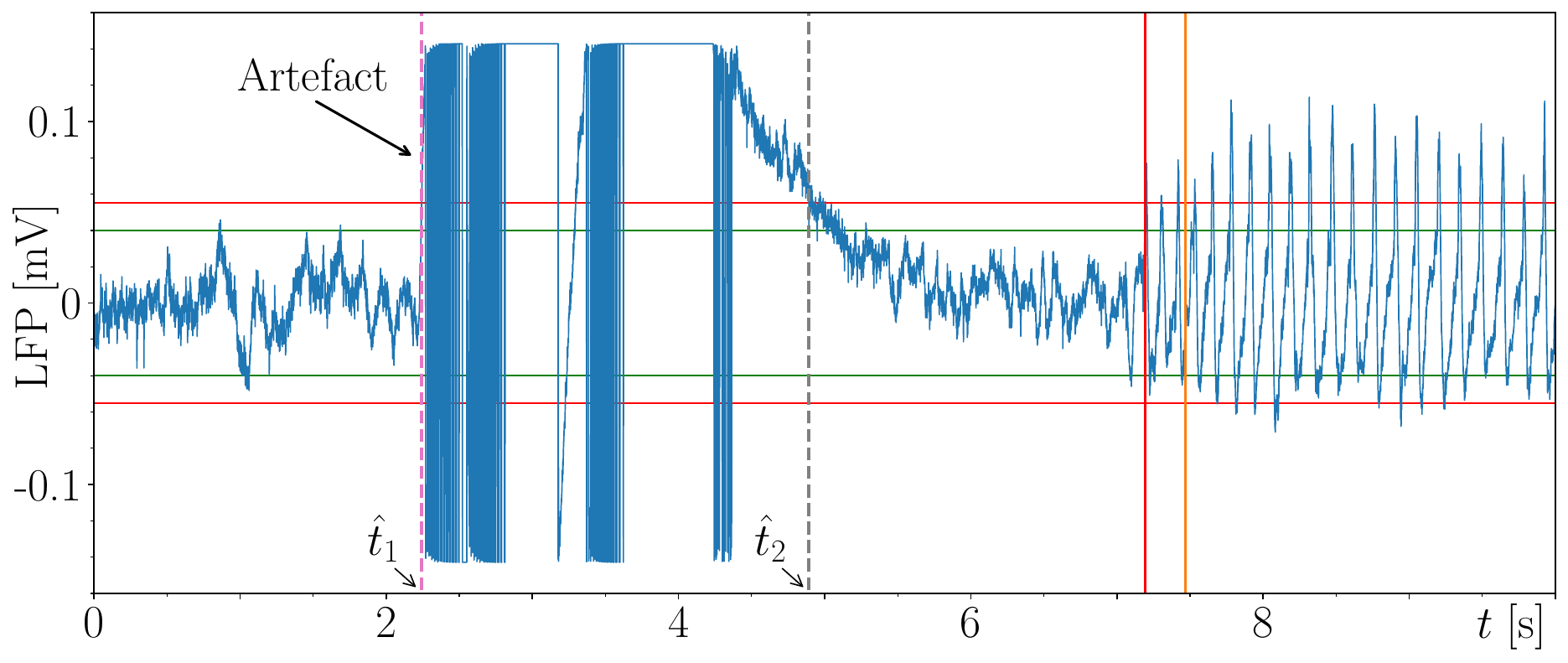}
    \caption{Accounting for artefacts in the NS state. Vertical dashed lines indicate the (pink) beginning and (grey) end of artefact activity according to our algorithm, i.e., the values of $\hat{t}_{1}$ and $\hat{t}_{2}$. Vertical solid lines indicate values of $t_{1}$ according to the (orange) expert's annotations and (red) algorithm when applied to a portion of the voltage recordings from session `T8C'. Algorithm parameters used: $\alpha=0.055$, $\beta=0.04$, $\tau_{\text{NS}}=5$, $\tau_{\text{S}}=2$, $\tau_{w}=1$, and $\Delta = 0.001$, horizontal lines indicate the thresholds of (in red) $\alpha$ and (in green) $\beta$.}
    \label{fig:Compare_detection_RatData_artefact}
\end{figure}

In this section we illustrate that our algorithm detects CTs between NS and S states in the voltage recordings at times which closely correspond with expert annotations of seizure onset and offset times. 
We also describe how our algorithm is adjusted to account for artefacts that sometimes appear in voltage recordings and are excluded from expert annotations. 
\\
\textbf{Algorithm vs. expert: } We first apply our algorithm to the portion of voltage recordings shown in Fig.~10~(b) in part {M3} of the Methods section in \hyperref[seref:MainText]{[1]}. To provide a comparison to the CTs detected in Figs.~\ref{fig:SeizureDetectionExample1} and \ref{fig:mu_sub_sup_vs_t_CTplots}, we use the same $\tau_{w}$, $\Delta$, $\tau_{\text{S}}$, and $\tau_{\text{NS}}$, and set $\alpha=0.055$ and $\beta=0.04$ since both the NS and S states in the model's output are approximately an order of magnitude larger than the corresponding states in the voltage recordings (see values on vertical axes). 
In Fig.~\ref{fig:Compare_ModelOutput_to_RatData_detection} we compare the values of $t_{1}$ and $t_{2}$ obtained by the algorithm to the corresponding electrical seizure onset and offset times provided by the expert. 
The vertical lines indicate the values of $t_{1}$ and $t_{2}$ according to (in orange and purple) the expert and (in red and green) the algorithm. 
The horizontal lines indicate the values of (in red) $\alpha$ and (in green) $\beta$ used by the algorithm. 
Figure~\ref{fig:Compare_ModelOutput_to_RatData_detection} shows that, for this choice of algorithm parameters, the expert's $t_{1}$ immediately precedes the algorithm's, and the algorithm's $t_{2}$ immediately precedes the expert's. 
In this case, we say the seizure interval defined by the algorithm's CT times is contained within the seizure interval defined by the expert's seizure onset and offset times. 
In Sec.~\hyperref[si:restimes_alg_expert]{S5} we describe how the algorithm's parameters are tuned to maximise the agreement between our algorithm and the expert in terms of these seizure intervals.
\\
\textbf{Accounting for artefacts: } 
We now describe how we alter the CT detection algorithm (presented in part {M4} of the Methods section in \hyperref[seref:MainText]{[1]}) to account for artefacts that sometimes appear in the voltage recordings; see part {M1} of the Methods section in \hyperref[seref:MainText]{[1]} for reasons why artefacts appear. 
In Fig.~\ref{fig:Compare_detection_RatData_artefact} we show a typical example of artefact activity for $t \in [2,5]$. 
We observe that, during this time interval, the LFP jumps between $-0.15$ and $0.15$ mV in short time intervals, sometimes within consecutive measurements ($0.001\text{s}$). 
Based on this observation, we adapt our algorithm to account for artefacts by altering condition {(a2)} of our algorithm, specified in part {M4} of the Methods section in \hyperref[seref:MainText]{[1]}, as follows: 
\begin{itemize}[labelindent=20pt,itemindent=1em,leftmargin=!]
    \item[(a2-1)\label{a2-1}] $|x(t)|$ exceeds $\alpha$ at time $t=t_{1}$ \textbf{and} $|x(t+\delta)-x(t)| < \xi$ for all $t \in \left[t_{1}, t_{1} + t_{w} \right]$, or,
    \item[(a2-2)\label{a2-2}] $|x(t)|$ exceeds $\alpha$ at time $t=t_{1}$ \textbf{and} $|x(t+\delta)-x(t)| \geq \xi$ for any $t \in \left[t_{1}, t_{1} + t_{w} \right]$.
\end{itemize}
If \hyperref[a2-1]{(a2-1)} is true then the algorithm continues as before to evaluate whether a CT from the NS to S state is detected at $t_{1}$. 
On the other hand, if \hyperref[a2-2]{(a2-2)} is true then we say artefact activity begins at time $t=\hat{t}_{1}=t_{1}$ and we continue to monitor $|x(t)|$ in moving windows. 
We say artefact activity ends at time $t=\hat{t}_{2}$ if $|x(\hat{t}_{2})| < \alpha$ \textbf{and} $|x(t+\delta)-x(t)| < \xi$ for all $t \in \left[\hat{t}_{2}, \hat{t}_{2} + t_{w} \right]$. 
We find most artefacts are accounted for by choosing $\xi = 0.2$. 

\begin{figure}[h!]
    \centering
    \includegraphics[width=0.7\linewidth]{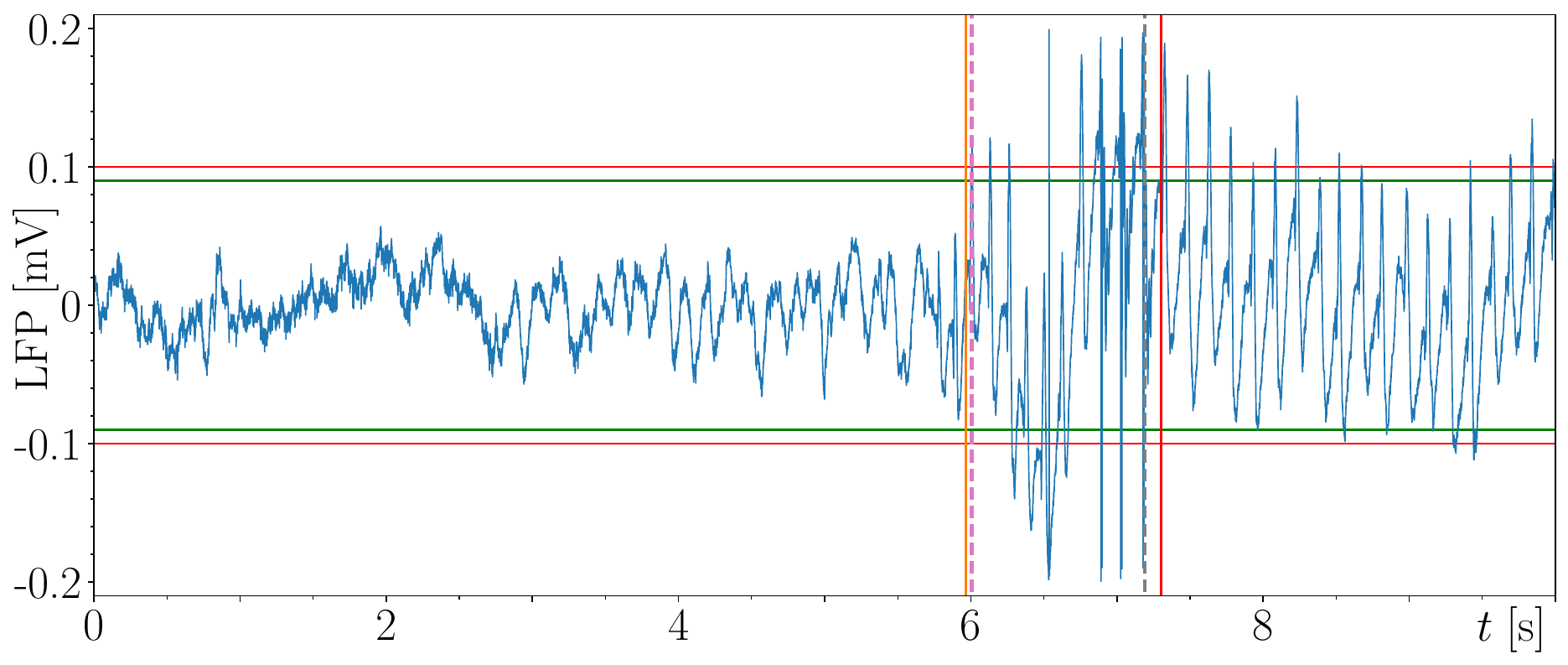}
    \caption{Accounting for exaggerated spikings which happen in conjunction with a CT. Vertical dashed lines indicate the (pink) beginning and (grey) end of exaggerated spiking activity according to our algorithm, i.e., the values of $\hat{t}_{1}$ and $\hat{t}_{2}$. Vertical solid lines indicate values of $t_{1}$ according to the (orange) expert's annotations and (red) algorithm when applied to a portion of the voltage recordings from session `K7E'. Algorithm parameters used: $\alpha=0.1$, $\beta=0.09$, $\tau_{\text{NS}}=3$, $\tau_{\text{S}}=2$, $\tau_{w}=1$, and $\Delta = 0.001$, horizontal lines indicate the thresholds of (in red) $\alpha$ and (in green) $\beta$.}
    \label{fig:Compare_detection_RatData_artefact_duringCT_AlgUpdate}
\end{figure}

While artefacts can appear at various locations in the voltage recordings, the example in Fig.~\ref{fig:Compare_detection_RatData_artefact} is chosen to show that for the current choice of $\tau_{\text{S}}$ and $\tau_{\text{NS}}$, the algorithm would have considered the time when artefact activity begins as the time when a CT from the NS to S state occurs. 
Furthermore, this alteration to the algorithm does not diminish the strong agreement between the expert and the algorithm; see the corresponding $t_{1}$ values in Fig.~\ref{fig:Compare_detection_RatData_artefact}. 

We find the following two scenarios unfold after an interval of artefact activity has ended according to the algorithm, $|x(t)|$ either remains below $\alpha$ or quickly re-exceeds $\alpha$. 
The second scenario indicates that the artefact activity happens in conjunction with (i) additional artefact activity, (ii) an almost-occurring CT from the NS to S state, or (iii) a CT from the NS to S state. 
In the case of (iii), what our algorithm considers as artefact activity preceding the CT typically coincides with what the expert may consider as `exaggerated spikings' near seizure onset/CT. 
Figure~\ref{fig:Compare_detection_RatData_artefact_duringCT_AlgUpdate} shows an example of artefact activity in the form of exaggerated spikings near CT. 
More specifically, Fig.~\ref{fig:Compare_detection_RatData_artefact_duringCT_AlgUpdate} shows an example of when the expert's $t_{1}$ (vertical orange line) is similar to $\hat{t}_{1}$ (vertical dashed pink line), the time where the algorithm says that artefact activity begins. It is only after this activity ends (vertical dashed grey line) that the algorithm says there is a CT from the NS to S state (vertical red line). 
While exaggerated spikings may not be related to artefacts in terms of their description in part {M1} of the Methods section in \hyperref[seref:MainText]{[1]}, and may instead be generated from particularly synchronous cortical activity for a few cycles of the spike-wave-discharge, it is important that we detect this activity and do not classify seizures in which it occurs around the time of seizure onset since such exaggerated spikings are not generated by our mathematical model; see Sec. \hyperref[si:TSPs_artefacts]{S10} for further details. 
A future iteration of our model may be able to accommodate such patterns.


\section*{S3: Influence of shear on residence times in NS and S states for noise-induced CTs \label{si:shear_restimes_noise}}

\begin{figure}[t]
\centering
\includegraphics[scale=0.34]{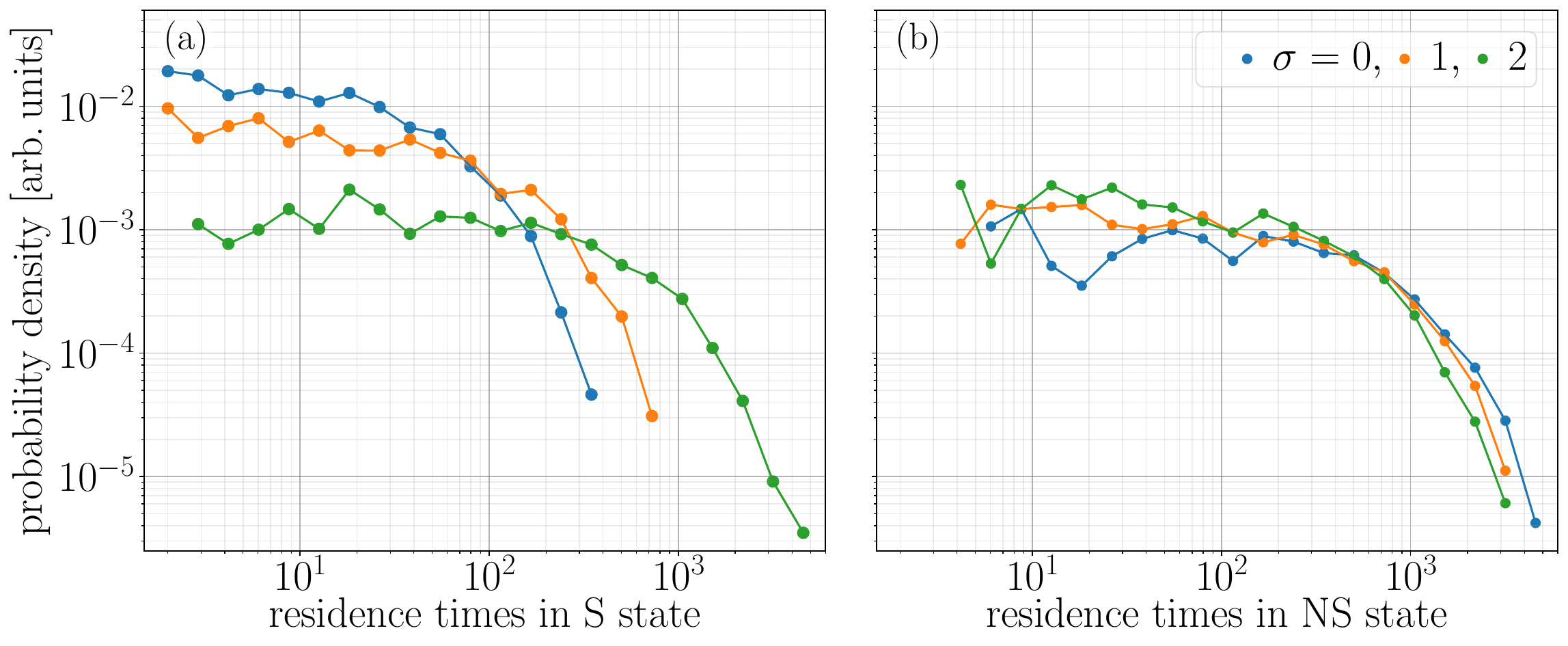}
\caption{Probability density of residence times in the model's (a) S state and (b) NS state for the values of $\sigma$ specified in the plot legend. The information presented here is based on 700 CTs between the NS and S states in the model as detected by the algorithm using the parameter values specified in Tables~3 and 4 in \hyperref[seref:MainText]{[1]}.}
\label{fig:ResTimesExample1}
\end{figure}

In this section we examine how $\sigma$, the shear parameter in our model, influences the model's residence times in the NS and S states. 
More specifically, for $\sigma = 0, 1,$ and $2$, we integrate the model forward in time up to $t=1 \times 10^{6}$ with the remaining model parameters as specified in NCT section of Table~3 in \hyperref[seref:MainText]{[1]}. 
We apply our algorithm to the model's output using the algorithm parameter values specified in Table~4 in \hyperref[seref:MainText]{[1]}. 
Our algorithm detected the following number of CTs between the NS and S states for different values of $\sigma$: 948 CTs for $\sigma = 0$, 984 CTs for $\sigma = 1$, and 715 CTs for $\sigma = 2$. 
We then compute residence times in the NS and S states based on the $t_{1}$ and $t_{2}$ values obtained for each choice of $\sigma$. 
Since the number of CTs vary non-linearly depending on $\sigma$, and we would like to compare residence times across different values of $\sigma$, we compute the probability density of residence times in the NS and S states for each $\sigma$ for 700 randomly chosen residence times, reflecting the smallest number of CTs detected when $\sigma = 2$. 
The resulting probability densities of residence times in the S state are shown in Fig.~\ref{fig:ResTimesExample1}~(a) and NS state in Fig.~\ref{fig:ResTimesExample1}~(b). 

Figures~\ref{fig:ResTimesExample1}~(a) and (b) show that the model is most likely to spend a short duration of time in a given state and least likely to spend a large duration of time.
Figure~\ref{fig:ResTimesExample1}~(a) shows that $\sigma$ has a significant influence on residence times in the S state, and Fig.~\ref{fig:ResTimesExample1}~(b) shows it has no significant influence on residence times in the NS state. 
More specifically, Fig.~\ref{fig:ResTimesExample1}~(a) shows that as $\sigma$ increases the probability that the model spends (i) a relatively short duration of time in the S state slightly decreases, and (ii) a relatively large duration of time in the S state increases significantly; for the largest recorded residence time in the S state when $\sigma = 0$, there is an approximately equal probability that the model is in the S state for nearly ten times longer when $\sigma = 2$. 


\section*{S4: Influence of shear on the detected \texorpdfstring{$t_{1}$}{TEXT} values for bifurcation and bifurcation/noise induced CTs \label{si:shear_t1times_bif}}

\begin{figure}[t]
\centering
\includegraphics[scale=0.35]{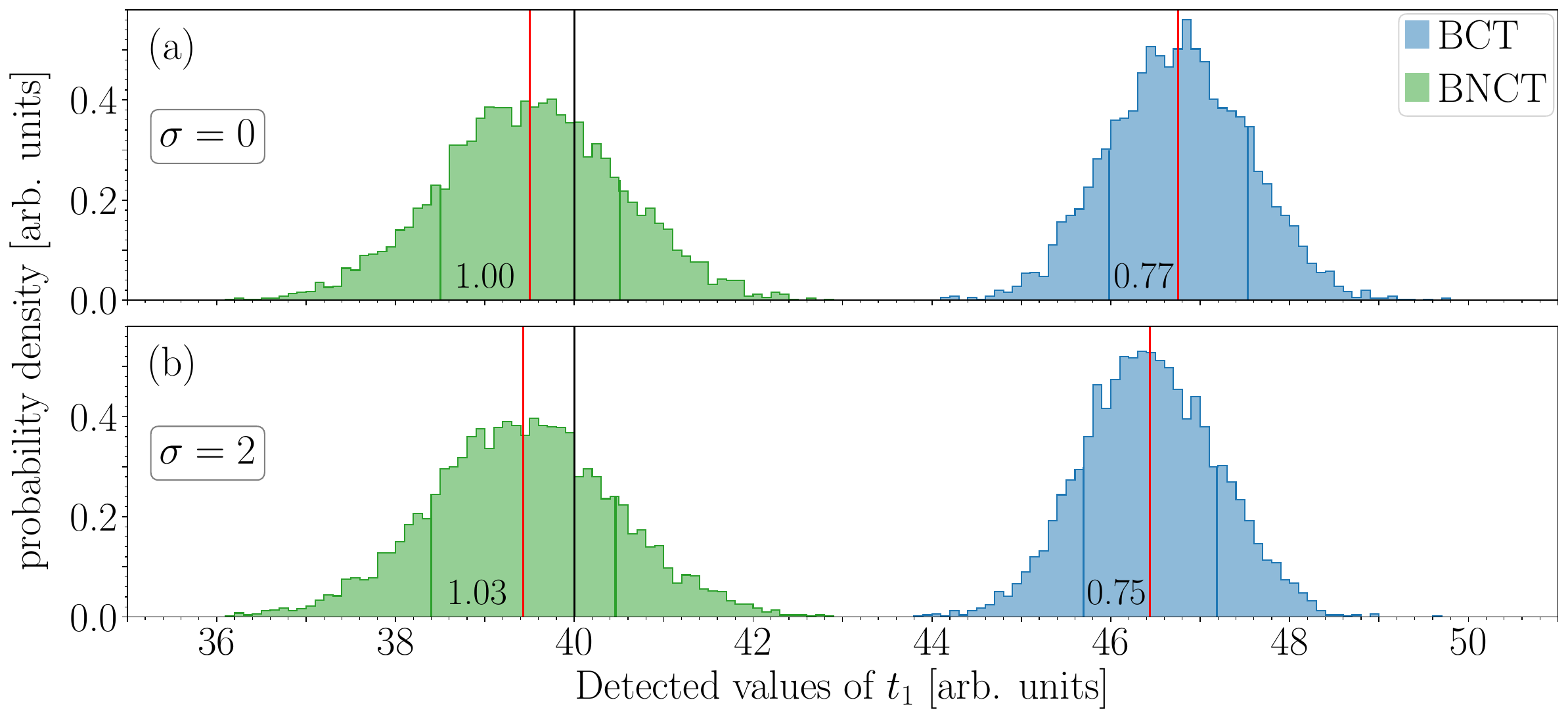}
\caption{Probability density of $t_{1}$ values for (blue) BCTs and (green) BNCTs with shear parameter set to (a) $\sigma=0$ and (b) $\sigma=2$. The vertical black line indicates when the respective bifurcations take place, vertical red lines indicates the mean $t_{1}$ for each CT-type, and the vertical blue and green lines indicate the corresponding standard deviation (quoted the left of each red line). The information presented here is based on 5000 examples of each CT-type.}
\label{fig:probdens_CTtime_SuperSub_sigma_0_1_2_}
\end{figure}

In this section we examine how $\sigma$, the shear parameter in our model, influences when BCTs and BNCTs from the NS to S state occur. 
More specifically, for $\sigma = 0$ and $2$, we generate 5000 different examples of BCTs and BNCTs for the model parameters specified in Table~3 in \hyperref[seref:MainText]{[1]}.
We apply our algorithm to the model's output using the algorithm parameter values specified in Table~4 in \hyperref[seref:MainText]{[1]}.
We then construct probability density diagrams (based on histograms) of the different $t_{1}$ values detected by our algorithm, these are shown in Fig.~\ref{fig:probdens_CTtime_SuperSub_sigma_0_1_2_}~(a) for $\sigma = 0$ and Fig.~\ref{fig:probdens_CTtime_SuperSub_sigma_0_1_2_}~(b) for $\sigma = 2$ where the data in blue corresponds to BCTs and in green corresponds to BNCTs. 
In both Figs.~\ref{fig:probdens_CTtime_SuperSub_sigma_0_1_2_}~(a) and (b) the vertical black line indicates when the respective bifurcations happen, the vertical red lines indicate the mean $t_{1}$ for the respective BCTs and BNCTs, the vertical blue and green lines indicate the corresponding standard deviation (which is quoted to the left of each red line). 
{Note that, in our model, we define the CT onset time $t_1$ as the time at which the on threshold value $\alpha$ is crossed when our algorithm detects a CT from the NS to S state of the model. This is why some CT onset times differ from the time at which a Hopf bifurcation occurs, which happens at $t=40$ for both BCTs and BNCTs (see Fig.~\ref{fig:mu_sub_sup_vs_t_CTplots})}

\textbf{BNCTs:} Figure~\ref{fig:probdens_CTtime_SuperSub_sigma_0_1_2_} shows that as $\sigma$ increases there are relatively small changes in the statistics of the detected $t_{1}$ values; the mean shifts slightly to the left and the standard deviation slightly increases. 
For both $\sigma$ values the mean $t_{1}$ value is relatively close when the bifurcation happens. 
From additional experiments we found that decreasing the noise level, $\nu$, results in the mean $t_{1}$ value increasing. If $\nu$ is small enough then the CTs becomes purely bifurcation-induced. 
If $\nu$ is too large then few CTs will take place after the bifurcation. 
The chosen value of $\nu$ in Table~3 in \hyperref[seref:MainText]{[1]} enables the model to produce BNCTs with an almost equal likelihood that noise or the bifurcation triggers the CT.

\textbf{BCTs:} Figure~\ref{fig:probdens_CTtime_SuperSub_sigma_0_1_2_} shows that as $\sigma$ increases there are relatively small, but more noticeable, changes in the statistics of the detected $t_{1}$ values than BNCTs. Specifically, the mean $t_{1}$ value and standard deviation both decrease.

Two aspects of Fig.~\ref{fig:probdens_CTtime_SuperSub_sigma_0_1_2_} are critical for the main aim of the paper of classifying CTs. 
The first is that there is no overlap between the different distributions; BNCTs are designed to cluster around the time of bifurcation while BCTs are designed to occur at a time long after the bifurcation.  
The second, there is no $t_{1} < 35$; this informs our choice of $T^{-} = -30$ when producing Fig.~4 in \hyperref[seref:MainText]{[1]}, allowing us to experiment with the relatively wide range of $t_{m}$ in the interest of finding the optimal way to differentiate between the three seizure generation mechanisms via our SVM classifier. 


\section*{S5: Tuning algorithm parameters \label{si:restimes_alg_expert}}

In this section we describe the procedure used to tune the algorithm's parameters to maximise the agreement between the CT detection algorithm and the expert's annotations of electrical seizure onset and offset times in voltage recordings.
By agreement we mean that, in the sense of a `receiver-operator-characteristic' (ROC) analysis, the proportion of correctly labelled seizure intervals is maximal and the proportion of incorrectly labelled non-seizure intervals is minimal, using the expert's annotations as the reference. 
Seizure and non-seizure intervals are defined by seizure onset and offset times provided by the expert and the values of $t_{1}$ and $t_{2}$ obtained by the algorithm when applied to voltage recordings from a given recording session. 
We consider a seizure interval defined by the algorithm's $t_{1}$ and $t_{2}$ values to be correctly labelled if it (i) is contained within, (ii) contains, or (iii) partially overlaps with a seizure interval defined by the expert's annotations. 
See Fig.~\ref{fig:Compare_ModelOutput_to_RatData_detection} for an example of (ii). 
Similarly, we consider an incorrectly labelled non-seizure interval to be a seizure interval defined by the algorithm's $t_{1}$ and $t_{2}$ values that is contained within a non-seizure interval defined by the expert's annotations. 

\begin{figure}[t]
    \centering
    \includegraphics[width=0.8\textwidth]{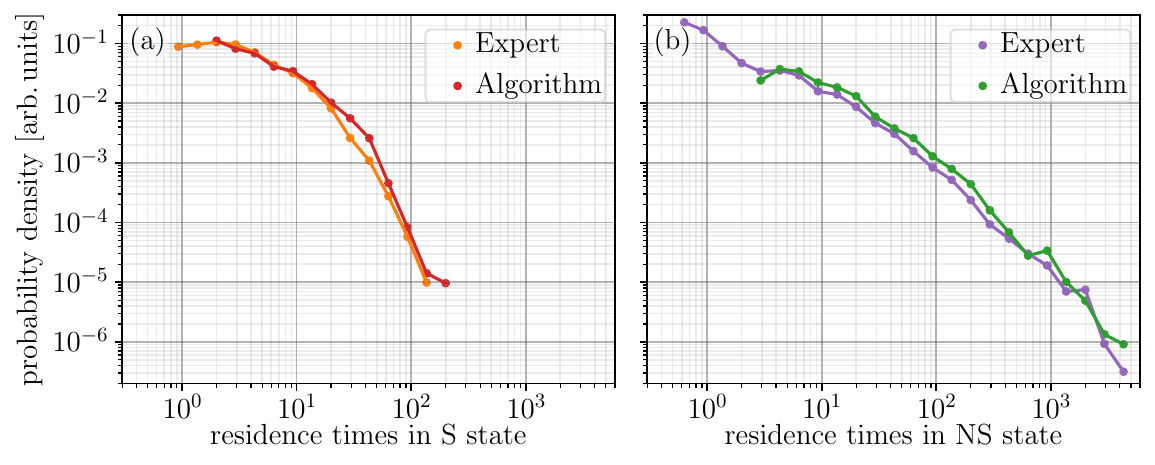}
    \caption{Probability density of residence times in the (a) S state and (b) NS state for rat T according to the (orange and purple) expert's annotations and (red and green) our CT detection algorithm.}
    \label{fig:ResTimes_RatT_Cian_vs_Algorithm}
\end{figure}

\begin{figure}[t!]
    \centering
    \includegraphics[width=0.8\textwidth]{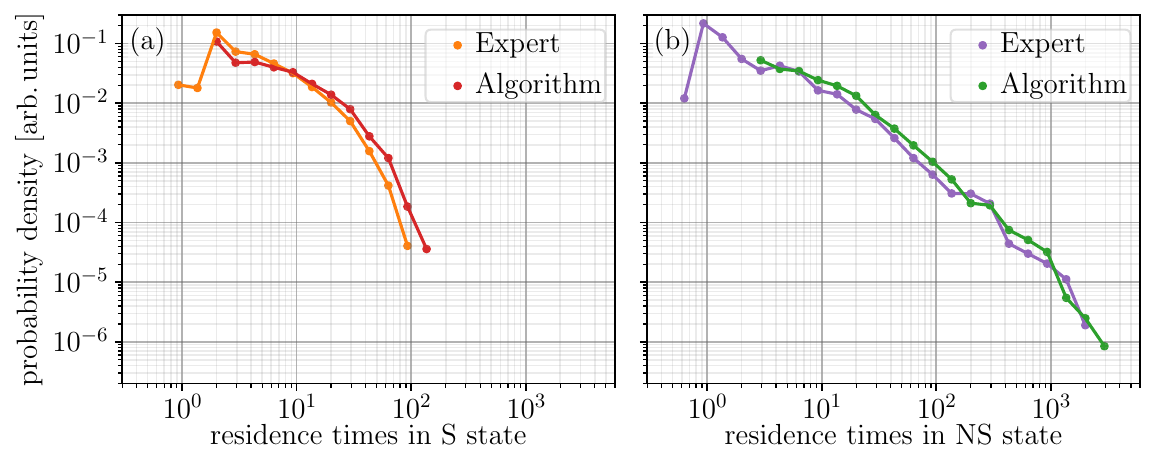}
    \caption{Probability density of residence times in the (a) S state and (b) NS state for rat K according to the (orange and purple) expert's annotations and (red and green) our CT detection algorithm.}
    \label{fig:ResTimes_RatK_Cian_vs_Algorithm}
\end{figure}

From several experiments we found that the above agreement is most sensitive to different values of $\alpha$. 
Therefore, for each recording session, we select a different $\alpha$ from $[0.03, 0.1]$ to obtain the best agreement with the expert annotations. 
$\beta$ is adjusted accordingly so that $\beta = \alpha - 0.01$. 
The values of $\tau_{w}$, $\Delta$, $\tau_{\text{S}}$, and $\tau_{\text{NS}}$ specified in Table~4 in \hyperref[seref:MainText]{[1]} provide a strong agreement between the algorithm and expert across all the analysed recording sessions from rats S, T, and K. 
This strong agreement is illustrated in Fig.~3 in \hyperref[seref:MainText]{[1]}, Figs.~\ref{fig:ResTimes_RatT_Cian_vs_Algorithm} and \ref{fig:ResTimes_RatK_Cian_vs_Algorithm} in terms of the probability distributions of residence times in the NS and S states for rats S, T, and K. 
These residence times are computed based on (i) the expert's annotations of all recording sessions from each rat and (ii) the values of $t_{1}$ and $t_{2}$ obtained from using the optimal choice of algorithm parameters to detect CTs between NS and S states in each of the same recording sessions. 

Figure~3 in \hyperref[seref:MainText]{[1]}, \ref{fig:ResTimes_RatT_Cian_vs_Algorithm}, and \ref{fig:ResTimes_RatK_Cian_vs_Algorithm} show there is relatively little difference between the probability distributions of residence times in the NS and S states according to the expert and the algorithm. 
Small differences occur at the end points between the corresponding distributions, this is mainly due to the choice of $\tau_{\text{NS}}$ and $\tau_{\text{S}}$ being less than the smallest residence times in the NS and S states according to the expert's annotations. 
As a result, a given seizure/non-seizure interval according to the algorithm may consist of multiple seizure and non-seizure intervals according to the expert. 
However, we find that when $\tau_{\text{NS}}$ and $\tau_{\text{S}}$ are chosen in line with the expert's annotations, we find this significantly reduces the agreement between the algorithm and the expert. 
{Further details on our algorithm's performance are provided in Flynn \textit{et al.}~\hyperref[seref:AlgorithmROC]{[2]}.}



\section*{S6: Additional clarification on how TSPs and \texorpdfstring{$m$(TSPs, $t_m$)}{TEXT} are computed for CTs in the model's output \label{si:Fig7_clarification}}

We obtain 100 BCTs, 100 BNCTs, and 100 NCTs from the model's NS to the S state as follows. For NCTs, we use the model parameters specified in Table~3 in \hyperref[seref:MainText]{[1]}, 
algorithm parameters specified in Table~4 in \hyperref[seref:MainText]{[1]}, 
and integrate the model until our algorithm detects 100 CTs between the NS and S states that satisfy the selection criteria stated in \hyperref[seref:MainText]{[1]}.
For BCTs and BNCTs, we use the model parameters specified in Table~3 in \hyperref[seref:MainText]{[1]}, different noise realisations for each simulation of the model, and detect CTs using algorithm parameters specified in Table~4 in \hyperref[seref:MainText]{[1]}. In Fig.~\ref{fig:probdens_CTtime_SuperSub_sigma_0_1_2_} in \hyperref[si:shear_t1times_bif]{S4} we show that all these BCTs and BNCTs satisfy the selection criteria as each $t_{1}>30$ and $t_{2}-t_{1}>10$.

Based on the above we obtain 100 examples of BCTs, NCTs, and BNCTs from the model. When computing the TSPs and $m$(TSPs, $t_m$) near each CT, it is convenient to work with a new time $T=t-T-{1}$ where $T\in\left[T^{-},T^{+}\right]=\left[-30,10\right]$ and each CT is detected at $T=0$. 
We then follow the steps in part {M5} of the Methods section in \hyperref[seref:MainText]{[1]} to compute the TSPs and $m$(TSPs, $t_m$) for a given $t_{w}=1$ and chosen {duration} $t_{m}$. 
This choice of $T^{-}$ and $T^{+}$ results in the GV, log$_{10}$GV, and AC being defined from $T=T^{-}+t_{w}=-29$, and log$_{10}$GV(AC) is defined from $T=T^{-}+2t_{w}=-28$. 
Further, the first values of $m(\text{GV},~t_{m})$, $m(\text{log$_{10}$GV},~t_{m})$, and $m(\text{AC},~t_{m})$ are defined from $T=T^{-}+t_{w}+t_{m}$ and $m(\text{log$_{10}$GV(AC)},~t_{m})$ is defined from $T=T^{-}+2t_{w}+t_{m}$. 
Thus, if $t_{m}=8$ then $m(\text{GV},~8)$ is defined from $T=-21$ and $m(\text{log$_{10}$GV(AC)},~8)$ is defined from $T=-20$, and so on.

\section*{S7: TSPs for a single example of each CT type \label{si:TSPs_single_example}}

\begin{figure}[t]
\centering
\includegraphics[width=0.85\textwidth]{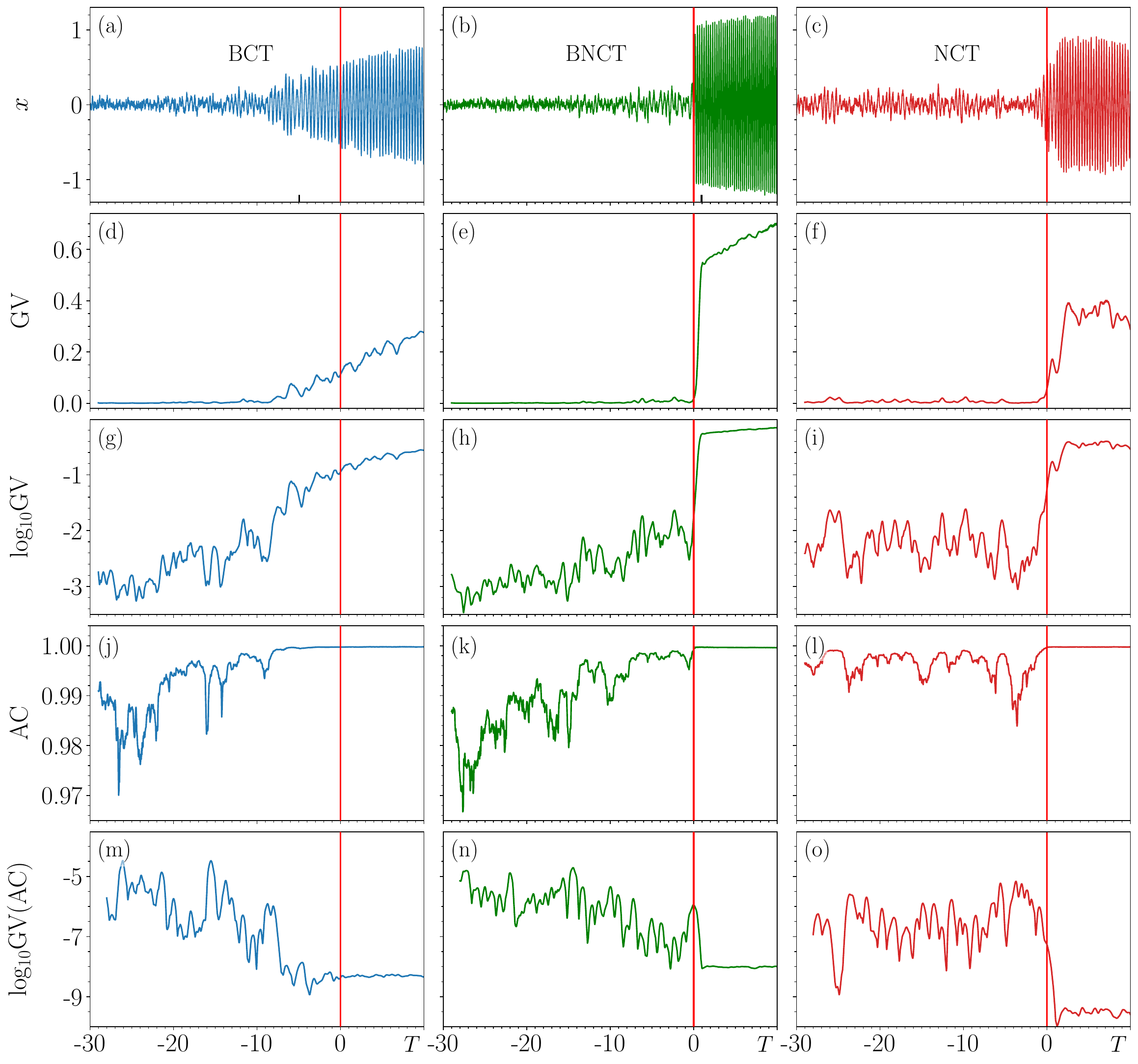}
\caption{TSPs of the model's output near an example of a (left-hand column) BCT, (middle column) BNCT, and (right-hand column) NCT. For each CT-type, the top row shows the model's output versus $T=t-t_{1}$ with $t_{1}$ detected by the algorithm. Inward black ticks in (a) and (b) indicate when the bifurcation occurs. Remaining rows show, for each CT-type, the corresponding behaviour of the TSPs specified in the vertical axes labels. Vertical red line in each panel emphasises $T=0$.}
\label{fig:EWSvalues}
\end{figure}

In this section we show how the time series properties (TSPs) of the model's output behave nearby an example of a BCT, BNCT, and NCT from the NS to S state. 
Each of these TSPs are specified in Table~5 in \hyperref[seref:MainText]{[1]}. 
We generate and detect the BCT and BNCT in Figs.~\ref{fig:EWSvalues}~(a) and (b) using the same parameters used in Fig.~\ref{fig:mu_sub_sup_vs_t_CTplots}. 
For the NCT case in Fig.~\ref{fig:EWSvalues}~(c) we generate and detect a CT for the same parameters used in Fig.~\ref{fig:SeizureDetectionExample1} and consider a similar example of an NCT to Fig.~\ref{fig:SeizureDetectionExample1}~(a). 
Each column of Fig.~\ref{fig:EWSvalues} shows, for different CT-types, the model's output and corresponding TSPs versus $T = t - t_{1}$ for the values of $t_{1}$ detected by our algorithm. 
The TSPs are computed from $T=-30$ in small moving windows of {duration} $t_{w} = 1$. 
Thus, for the reasons outlined in part {M5} of the Methods section in \hyperref[seref:MainText]{[1]}), the GV, log$_{10}$GV, and AC are defined from $T=-29$, and log$_{10}$GV(AC) is defined from $T=-28$. 
\\
We now discuss how each of these TSPs behave in relation to the different CT types. Our results are in agreement with the concept of `critical slowing down'; see \hyperref[seref:MainText]{[1]}).
\\
\textbf{GV: } 
{
Figures.~\ref{fig:EWSvalues}~(d)
}
to (f) show that the GV of the model's output in the S state is much larger than the NS state for each CT-type. One would naturally expect this to happen when a transition from a small to a large-amplitude of oscillation occurs. 
In the BCT case (Fig.~\ref{fig:EWSvalues}~(d)), there is a noticeably large and relatively constant increase in GV from $T \approx -7$ onwards, which is just before the bifurcation point indicated by the inward black tick. For the BNCT and NCT cases (Figs.~\ref{fig:EWSvalues}~(e) and (f)), no significant change in GV is seen before the CTs are detected. However, after the CTs are detected, the increase in GV is much greater for the BNCT case than the NCT case. This difference arises because the radius of the limit cycle for the BNCT case is larger than for the NCT case (see Eq.~(3) in \hyperref[seref:MainText]{[1]}), thereby resulting in a more dramatic change in the model's output.
\\
\textbf{log$_{10}$GV: } Figures~\ref{fig:EWSvalues}~(g) and (h) show log$_{10}$GV increases at a relatively small but constant rate before the CT in the BCT and BNCT cases, whereas Fig.~\ref{fig:EWSvalues}~(i) shows no such increase for the NCT case. Similar to the GV results, there is a sudden and large increase after the CT for both the BNCT and NCT cases. Notably, for the BCT case, log$_{10}$GV starts to increase at an earlier time than GV in Fig.~\ref{fig:EWSvalues}~(a), thus demonstrating the benefit of monitoring log$_{10}$GV in addition to GV.
\\
\textbf{AC: } Figures~\ref{fig:EWSvalues}~(j) and (k) show AC increases at a relatively constant rate before the BCTs and BNCTs. For the BCT case, AC remains close to $1$ from $T \approx -7$ (mirroring the trend seen in Fig.~\ref{fig:EWSvalues}~(d)). On the other hand, it is only after the CT that AC remains near $1$ for the BNCT case. In stark contrast, Fig.~\ref{fig:EWSvalues}~(l) shows that the AC remains relatively close to $1$ for the NCT case, with no trends present long before the CT. Interestingly, there is a local minimum in AC present before each CT. This occurs because, regardless of the AC's value before the CT, the model's output becomes highly correlated with itself after the CT due to the inherent periodicity of the limit cycle, thus causing the AC to quickly increase to $1$. We also see less fluctuations in the AC as the CT is approached for the BCT and BNCT cases.
\\
\textbf{log$_{10}$GV(AC): } Figures.~\ref{fig:EWSvalues}~(m)-(o) provide a more informative picture on the fluctuations mentioned above by examining how log$_{10}$GV(AC) changes over time. For the BCT and BNCT cases, there is an overall constant decrease in log$_{10}$GV(AC) before the CT. Log$_{10}$GV(AC) remains at a similar and constant value after both CTs. For the NCT case, log$_{10}$GV(AC) remains relatively constant before the CT and quickly decreases to a smaller value after the CT than the BCT and BNCT.


\section*{S8: TSPs of model's output computed for \texorpdfstring{$t_{w}=1$}{TEXT} and \texorpdfstring{$2$}{TEXT}\label{si:compare_tw}}

\begin{figure}[t!]
\centering
\includegraphics[width=0.85\textwidth]{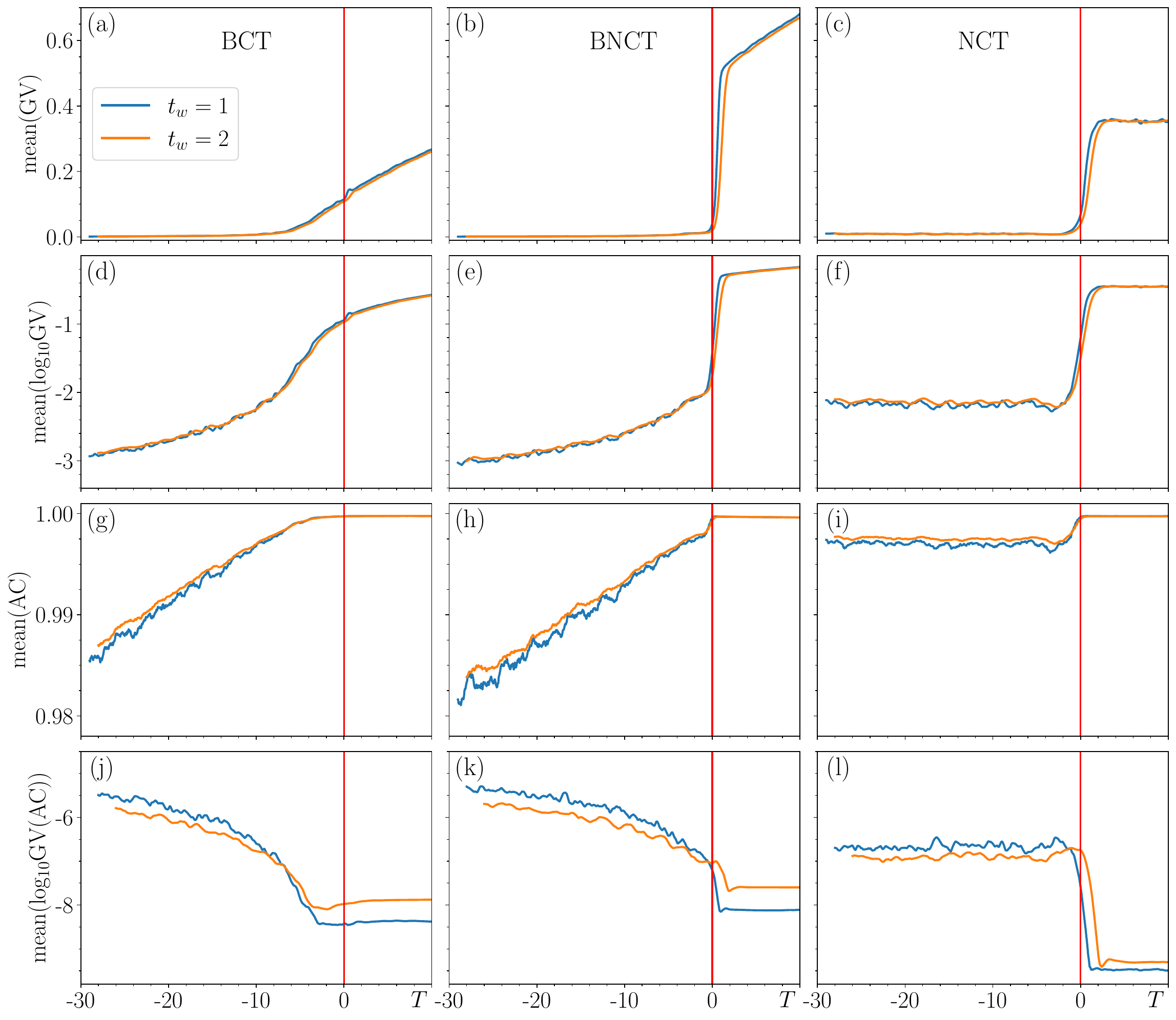}
\caption{Comparing different choices of $t_{w}$ when computing the time series properties (TSPs) of model's output near BCTs (left-hand column), BNCTs (middle column), and NCTs (right-hand column). Each panel shows the mean of each TSP versus $T=t-t_{1}$ for the $t_{1}$ values detected by the algorithm. Each curve represents the mean computed from 100 examples of each CT-type for a given $T$.}
\label{fig:compare_tw}
\end{figure}

In this section we provide clarification on comments made in \hyperref[seref:MainText]{[1]} in relation to tuning the model's parameters so that its output mimics the three important characteristics of real seizure activity in the voltage recordings.
In particular, we mentioned that we compromise on characteristic (iii), ``Intrinsic time scales of the NS and S states.''. We now show that there are little consequences from this in terms of the main aim of this paper, to identify the CT-type responsible for generating a given seizure in voltage recordings. 

In Fig.~10~(c) and (d) in \hyperref[seref:MainText]{[1]}, we showed that the timescale of oscillation of the model's S state is approximately half that of the dominant timescale of real seizure activity. 
Thus, when choosing a $t_{w}$ to compute the TSPs (the time duration of the moving window), the {duration} $t_{w}$ for the model should be twice that of the value used for the voltage recordings. 
However, we show in Fig.~\ref{fig:compare_tw} that there are no major differences in the TSPs of the model's output when computed using either $t_{w}=1$ (in blue) or $t_{w}=2$ (in orange). We also find that whether the SVM is trained using TSPs that are calculated with $t_{w}=1$ or $2$, there is little impact on classifying seizure generation mechanisms in the voltage recordings. For convenience, we choose $t_{w}=1$ when computing TSPs from both the model and the voltage recordings. 


\section*{S9: Computing SVM feature importance via the mean permutation importance (MPI) \label{si:MPI}}

In this section we follow the steps outlined in part {M6} of the Methods section in \hyperref[seref:MainText]{[1]} on how to compute the `mean permutation importance' (MPI). The MPI provides us with a means to quantify the importance of each feature used to construct SVM type-1, 2, and 3 (see Table~6 in \hyperref[seref:MainText]{[1]} for details on which features are used to construct each SVM). 
The results discussed in this section serve to accompany the results shown in Fig.~6 in \hyperref[seref:MainText]{[1]} where SVM accuracy versus $T$ is shown for SVM type-1, 2, and 3 with $t_{m}= 4, 8,$ and $12$. 
Since we trained and tested a new SVM for different values of $T$ in Fig.~6 in \hyperref[seref:MainText]{[1]}, we compute the MPI the corresponding values of $T$. 

\begin{figure}[t]
\centering
\begin{tikzpicture}
    \node (image) at (current page.center) {\includegraphics[width=0.95\textwidth]{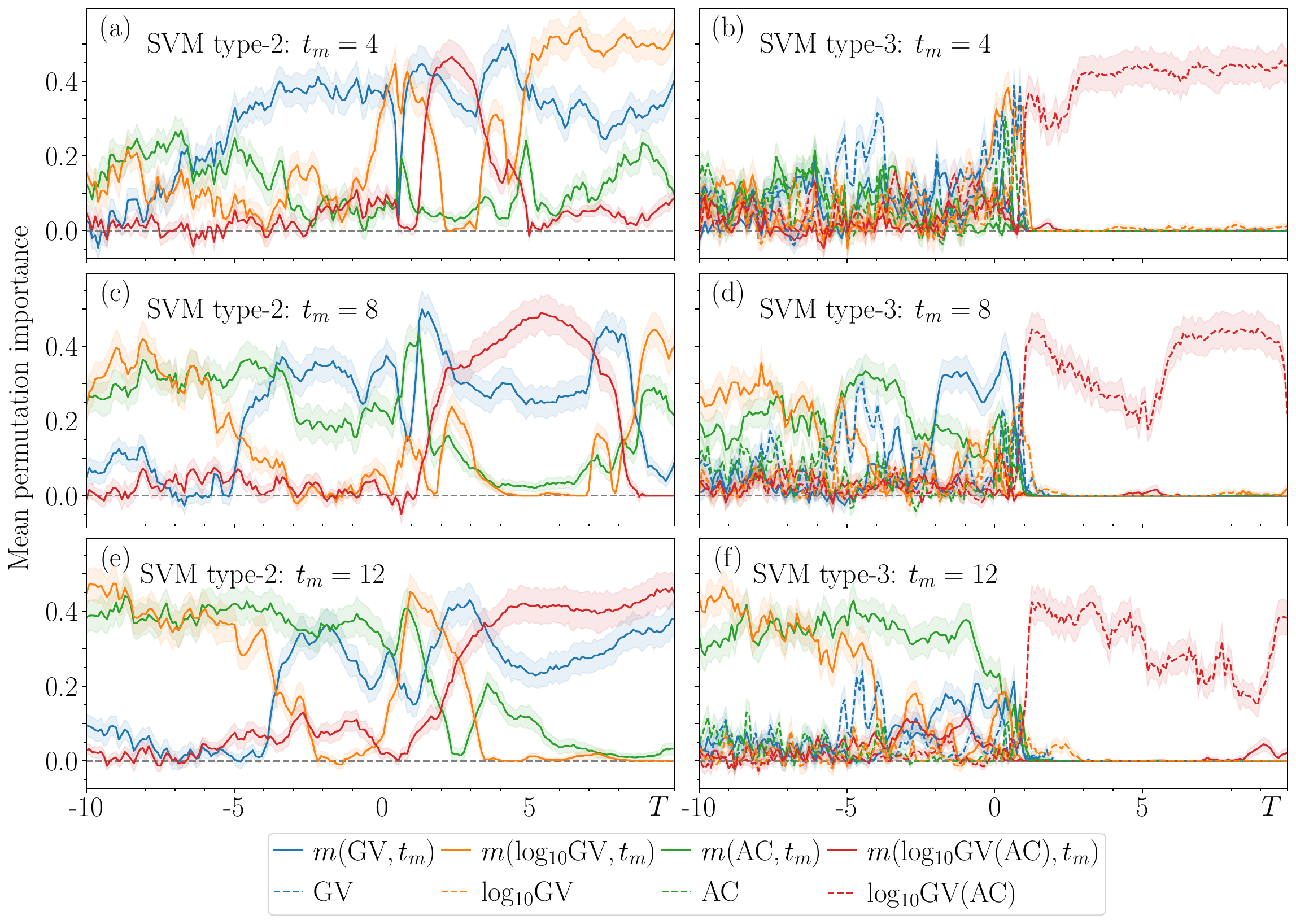}};
    \begin{scope}[shift={(image.south west)},x={(image.south east)},y={(image.north west)}]     
    {\node[fill=white,font=\large,inner sep=2pt] (Rd) at (0.53,0.57) {$\,$};}
    \end{scope}
\end{tikzpicture}
\caption{Feature importance for SVMs type-2 and 3: Mean permutation importance (solid and dashed curves) and its standard deviation (cloud surrounding each curve) at time $T=t-t_{1}$ for a given feature specified in the plot legend for (left-hand column) SVM type-2 and (right-hand column) SVM type-3 with (a) and (b) $t_{m}=4$, (c) and (d) $t_{m}=8$, and (e) and (f) $t_{m}=12$.}
\label{fig:SVM_MPI_slopes_slopesANDvalues}
\end{figure}

\begin{figure}[t!]
\centering
\includegraphics[width=0.9\textwidth]{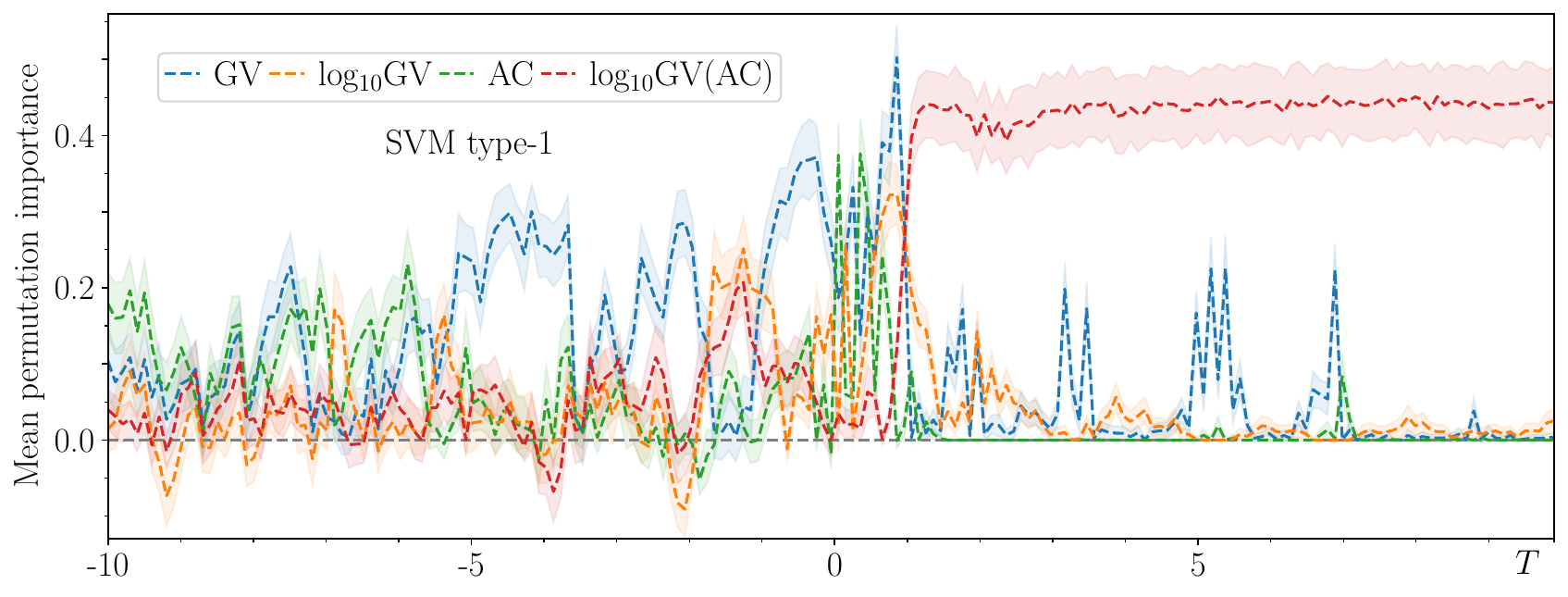}
\caption{Feature importance for SVM type-1: Mean permutation importance (dashed curves) and its standard deviation (cloud surrounding each curve) at time $T=t-t_{1}$ for a given feature specified in the plot legend.}
\label{fig:SVM_MPI_values}
\end{figure}

Figures~\ref{fig:SVM_MPI_slopes_slopesANDvalues} and \ref{fig:SVM_MPI_values} show the resulting MPI versus $T$. 
Each of the solid and dashed curves are coloured to correspond with the features specified in the legend. 
These curves are surrounded by a cloud of the corresponding colour which describes the standard deviation of the MPI for a given feature. 
Note, the results shown in Fig.~\ref{fig:SVM_MPI_slopes_slopesANDvalues}~(c) informed the choice of features that were used to construct Fig.~5 in \hyperref[seref:MainText]{[1]}. 


Figures~\ref{fig:SVM_MPI_slopes_slopesANDvalues} and \ref{fig:SVM_MPI_values} show that for each SVM, (i) some features are much more important than others, (ii) the importance of each feature varies depending on the value of $T$, and (iii) the importance also varies depending on the {duration} $t_{m}$ in the case of SVM type-2 and 3. 
For instance, when comparing the results shown in Figs.~\ref{fig:SVM_MPI_slopes_slopesANDvalues}~(a) and (e) for $T>5$, Fig.~\ref{fig:SVM_MPI_slopes_slopesANDvalues}~(a) shows that when $t_{m}=4$ the most important feature is $m(\text{log$_{10}$GV},~t_{m})$ and least important feature is $m(\text{log$_{10}$GV(AC)},~t_{m})$, 
whereas Fig.~\ref{fig:SVM_MPI_slopes_slopesANDvalues}~(c) shows the opposite for $t_{m}=12$. 

It is also worth noting the significant changes in the MPI of some features near $T=0$ for SVM type-2 and 3. 
For example, Fig.~\ref{fig:SVM_MPI_slopes_slopesANDvalues}~(a) shows that $m(\text{GV},~t_{m})$ suddenly changes from being the most important to least important and back to the most important feature over a small range of $T$ values near $0$. 
On the other hand, Fig.~\ref{fig:SVM_MPI_slopes_slopesANDvalues}~(a) shows that $m(\text{log$_{10}$GV},~t_{m})$ and $m(\text{log$_{10}$GV(AC)},~t_{m})$ change from being the least important to most important and back to least important feature over a small range of $T$ values between $0$ and $5$. 
Similar behaviour is seen in Figs.~\ref{fig:SVM_MPI_slopes_slopesANDvalues}~(c) and (e) but to a much less significant extent.

While Figs.~\ref{fig:SVM_MPI_slopes_slopesANDvalues}~(a), (c), and (e) show that the most important feature for SVM type-2 often varies for $T>0$, 
a different picture emerges for SVM type-3 in Figs.~\ref{fig:SVM_MPI_slopes_slopesANDvalues}~(b), (d), and (f) as $\text{log$_{10}$GV(AC)}$ is the most importance feature for each choice of $t_{m}$ when $T > 0$. 
Similar behaviour is shown in Fig.~\ref{fig:SVM_MPI_values} for SVM type-1.
Figures~\ref{fig:SVM_MPI_slopes_slopesANDvalues}~(b), (d), and (f) also show that the $m$(TSPs, $t_m$) play a more significant role in the accuracy of SVM type-3 for larger $t_{m}$ before the CT. 
This further illustrates the effects of combining the TSPs and their slopes in SVM type-3. 
More specifically, for $T<0$, the features which are most important for SVM type-2 are the most important for SVM type-3, 
for $T>0$ the features which are most important in SVM type-1 are the most important for SVM type-3. 


\section*{S10: TSPs and \texorpdfstring{$m$(TSPs, $t_m$)}{TEXT} for CTs not suitable for classification \label{si:TSPs_artefacts}}

In this section we show why it is not suitable to classify CTs from the NS to S state in the voltage recordings when almost-occurring CTs, artefacts, or exaggerated spikings happen in close proximity to the CT. More specifically, we show that some of the corresponding TSPs and $m$(TSPs, $t_m$) listed in Table~5 in \hyperref[seref:MainText]{[1]} become obscured and deviate significantly from those of a typical CT from the NS to S state when the above happen nearby the CT. This is important to show because if the SVM classifies CTs using these obscured TSPs and $m$(TSPs, $t_m$) then the classifications are much less reliable and may skew the main results presented in Table~2 in \hyperref[seref:MainText]{[1]}. 

\begin{figure}[t!]
\centering
\includegraphics[width=0.88\textwidth]{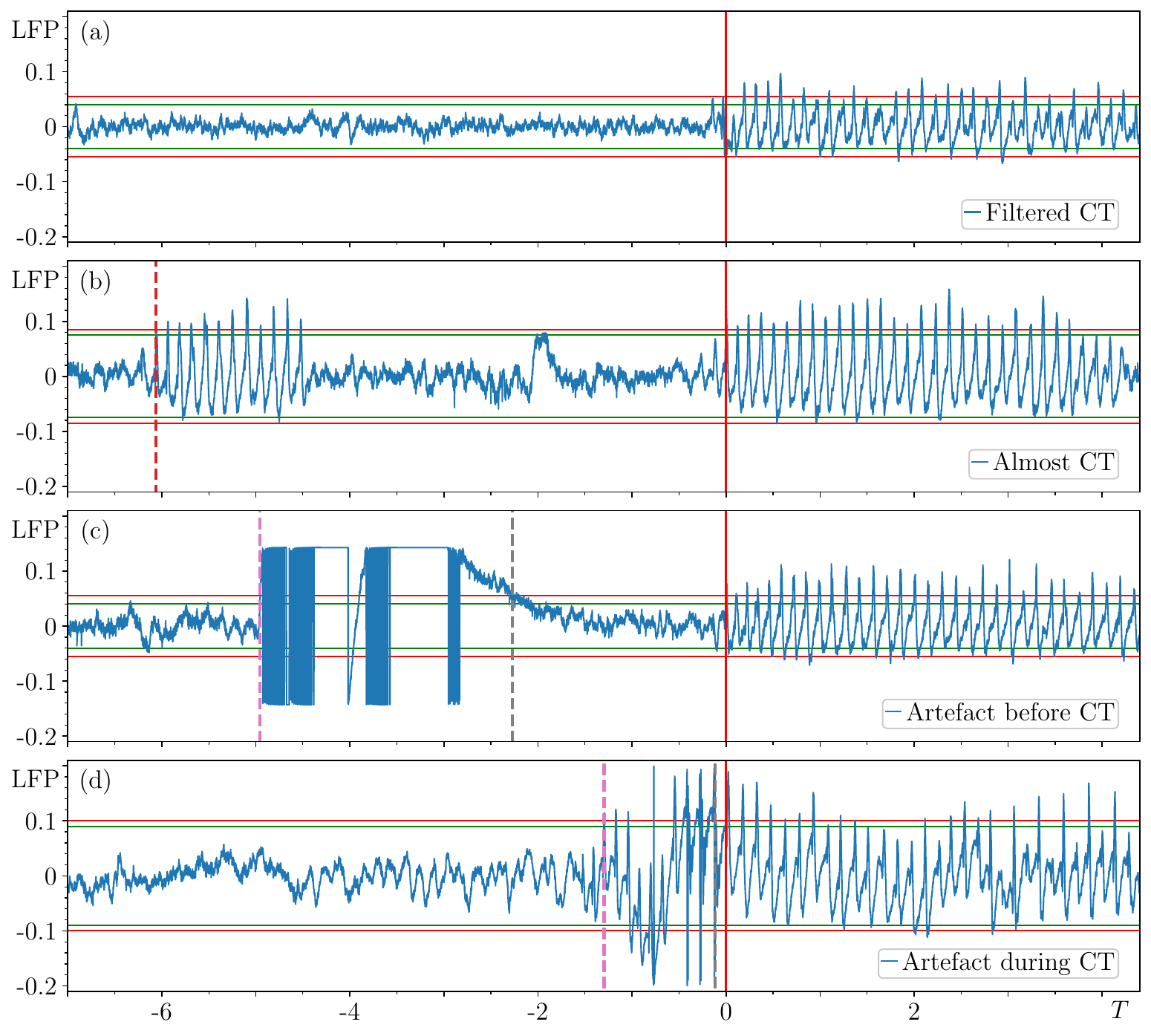}
\caption{Examples of detected CTs in voltage recordings that (a) are suitable for classification, (b)-(d) are not suitable for reasons given in each plot legend, and are plotted versus $T=t-t_{1}$ for $t_{1}$ values detected by our algorithm. 
The vertical lines indicate (solid red) when $T=0$, (dashed red) an almost-occurring CT from the NS to S state, (dashed pink and grey) the beginning and end of artefact activity or exaggerated spikings, according to our algorithm. 
Horizontal solid lines indicate thresholds of (red) $\alpha$ and (green) $\beta$ used by the algorithm.}
\label{fig:_traj_Artefact1_Artefact2_CTalmost_NoArtefact_}
\end{figure}

In Fig.~\ref{fig:_traj_Artefact1_Artefact2_CTalmost_NoArtefact_} we show the examples of CTs that we compute the TSPs and $m$(TSPs, $t_m$) of in this section. 
In each example we plot the LFP versus $T=t-t_{1}$, where $t_{1}$ corresponds to the time where the algorithm detects a CT from the NS to S state. 
The example in Fig.~\ref{fig:_traj_Artefact1_Artefact2_CTalmost_NoArtefact_}~(a) (same as in Fig.~\ref{fig:Compare_ModelOutput_to_RatData_detection}) is a typical example of a CT that satisfies the classification conditions specified in \hyperref[seref:MainText]{[1]}, {C1}-{C4} and is chosen as a reference point to compare with the following examples of CTs that satisfy {C1}-{C3} but dot not satisfy {C4}: Fig.~\ref{fig:_traj_Artefact1_Artefact2_CTalmost_NoArtefact_}~(b) (taken from recording session `T8C') shows an example of an almost-occurring CT before a CT from the NS to S state, Fig.~\ref{fig:_traj_Artefact1_Artefact2_CTalmost_NoArtefact_}~(c) (same as in Fig.~\ref{fig:Compare_detection_RatData_artefact}) shows an example of when artefact activity occurs before a CT from the NS to S state, and Fig.~\ref{fig:_traj_Artefact1_Artefact2_CTalmost_NoArtefact_}~(d) (same as in Fig.~\ref{fig:Compare_detection_RatData_artefact_duringCT_AlgUpdate}) shows an example of when exaggerated spikings occur in conjunction with a CT from the NS to S state. 
Note also that the example in Fig.~\ref{fig:_traj_Artefact1_Artefact2_CTalmost_NoArtefact_}~(b) provides similar insight to CTs that do not satisfy {C1} and CTs generated by the model that does not satisfy the corresponding selection criteria. 
The horizontal red and green lines in each panel of Fig.~\ref{fig:_traj_Artefact1_Artefact2_CTalmost_NoArtefact_} indicate the values of $\alpha$ and $\beta$ used to detect CTs between the NS and S states, these values are chosen according to the procedure outlined in \hyperref[{si:restimes_alg_expert}]{S5}. 
The vertical red solid lines indicates a CT from the NS to S state. Vertical dashed lines indicate (in red) an almost-occurring CT from the NS to S state, (in pink) the beginning and (in grey) end of artefact activity or exaggerated spikings. 

\begin{figure}[t]
\centering
\includegraphics[width=\textwidth]{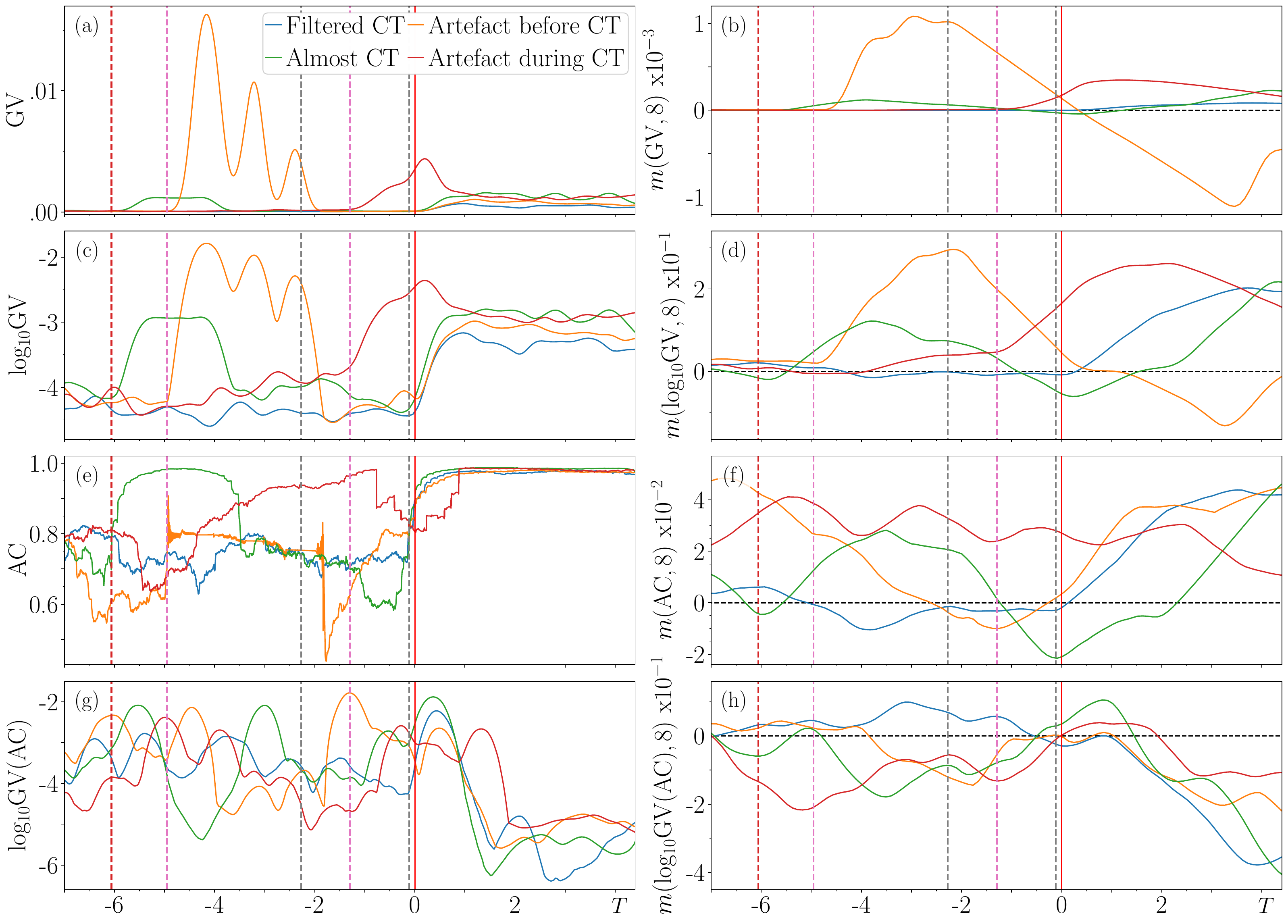}
\caption{Obscured TSPs and $m$(TSPs, $t_m$) of the voltage recording examples shown in Fig.~\ref{fig:_traj_Artefact1_Artefact2_CTalmost_NoArtefact_} and plotted versus the same $T$ values where the colour of each curve corresponds to the different examples specified in the plot legend. 
The vertical lines indicate (solid red) when $T=0$, (dashed red) an almost-occurring CT from the NS to S state, (dashed pink and grey) the beginning and end of artefact activity or exaggerated spikings, according to our algorithm. }
\label{fig:_Compare_TSPs_mTSPs_for_Artefact1_Artefact2_CTalmost_NoArtefact_}
\end{figure}

In Fig.~\ref{fig:_Compare_TSPs_mTSPs_for_Artefact1_Artefact2_CTalmost_NoArtefact_} we show how each of the TSPs and $m$(TSPs, $t_m$) with $t_{m}=8$ behave versus $T$ for the voltage recordings shown in Fig.~\ref{fig:_traj_Artefact1_Artefact2_CTalmost_NoArtefact_}. 
More specifically, the different coloured curves in each panel of Fig.~\ref{fig:_Compare_TSPs_mTSPs_for_Artefact1_Artefact2_CTalmost_NoArtefact_} correspond to one of the different TSPs or $m$(TSPs, $t_m$) of the examples shown in Fig.~\ref{fig:_traj_Artefact1_Artefact2_CTalmost_NoArtefact_}; the {blue} curves correspond to Fig.~\ref{fig:_traj_Artefact1_Artefact2_CTalmost_NoArtefact_}~(a), the {green} curves correspond to Fig.~\ref{fig:_traj_Artefact1_Artefact2_CTalmost_NoArtefact_}~(b), the {orange} curves correspond to Fig.~\ref{fig:_traj_Artefact1_Artefact2_CTalmost_NoArtefact_}~(c), and the {red} curves correspond to Fig.~\ref{fig:_traj_Artefact1_Artefact2_CTalmost_NoArtefact_}~(d). 
The different coloured vertical lines have the same meaning as in Fig.~\ref{fig:_traj_Artefact1_Artefact2_CTalmost_NoArtefact_}. 
The pink and grey vertical dashed lines at $T\approx-5$ and $-2.25$ correspond with the example shown in Fig.~\ref{fig:_traj_Artefact1_Artefact2_CTalmost_NoArtefact_}~(c), the other pink and grey vertical dashed lines at $T\approx-1.25$ and $-0.1$ correspond with the example shown in Fig.~\ref{fig:_traj_Artefact1_Artefact2_CTalmost_NoArtefact_}~(d).  

Figure~\ref{fig:_Compare_TSPs_mTSPs_for_Artefact1_Artefact2_CTalmost_NoArtefact_} shows that when $T=0$ there is an almost immediate change in the TSPs and a slightly delayed and longer lasting change in the $m$(TSPs, $t_m$) in response to the presence of almost-occurring CTs and the different artefacts. 
This delayed and longer lasting change in the $m$(TSPs, $t_m$) results in major differences from the $m$(TSPs, $t_m$) of a typical CT with no feature obscuring activity nearby, i.e., when comparing the {blue} curves to the others. 
These differences are visible at several values of $T$ and, crucially for our SVM setup, are noticeable for $T=2$ which is the value of $T$ we base our final results on in Table~2 in \hyperref[seref:MainText]{[1]}. 
In the following headings we highlight these differences. 
\\
\textbf{Almost-occurring CT: } (Comparing {blue} and {green} curves at $T=2$) Considerable differences particularly in $m(\text{log$_{10}$GV},~t_{m})$ and $m(\text{AC},~t_{m})$.
\\
\textbf{Artefact before CT: } (Comparing {blue} and {orange} curves at $T=2$) Major differences particularly in $m(\text{GV},~t_{m})$ and $m(\text{log$_{10}$GV},~t_{m})$ but quite similar in $m(\text{AC},~t_{m})$ and $m(\text{log$_{10}$GV(AC)},~t_{m})$.
\\
\textbf{Exaggerated spikings near CT: } (Comparing {blue} and {red} curves at $T=2$) 
Major differences in all $m$(TSPs, $t_m$), there appears as if no changes happens in $m(\text{AC},~t_{m})$ and $m(\text{log$_{10}$GV(AC)},~t_{m})$. 

As a final word, we would like to add the following for clarification. Since the differences highlighted above are only present in the $m$(TSPs, $t_m$) when $T=2$ one may argue that in the interest of classifying as many CTs in the voltage recordings as possible, one should do so using SVM type-1 instead of SVM type-3 so that less CTs like those shown in Fig.~\ref{fig:_traj_Artefact1_Artefact2_CTalmost_NoArtefact_} need to be excluded via conditions {C1}-{C4}. However, we found that SVM type-1 provides classifications which are in total disagreement with the fundamental characteristics of the different CTs that are captured by the $m$(TSPs, $t_m$); there were instances when some CTs were classified as BCTs without the typical the increase in GV and AC before the CT, and other CTs were classified as NCTs with increases in GV and AC before the CT that are typical of BCTs. 
Therefore it is necessary to consider both the TSPs and the $m$(TSPs, $t_m$) when classifying CTs in voltage recordings of real seizure activity even if this restricts the amount of CTs that we can classify. 


\section*{S11: SVM classifications and mean feature fit error vs. \texorpdfstring{T}{TEXT} \label{si:SVM_classifications_vs_T}}

In this section we show how the proportion of filtered CTs classified as a particular CT-type for a given $T$, denoted by $N_{\text{type}}(T)/N_{\text{filt}}(T^{-}, T^{+}) \in [0,1]$, varies for $T \in \left[ -10, 10 \right]$ when setting $T^{-}=-20$ and $T^{+}=10$. These classifications are provided by a type-3 SVM with $t_{m}=8$, i.e., using TSPs and $m$(TSPs, $t_m$) as features. 
The TSPs and $m$(TSPs, $t_m$) with $t_{m}=8$ we use for this SVM are computed from the voltage recordings of rat K. 
CTs between the NS and S states are detected by our algorithm (using parameters in Table~4 in \hyperref[seref:MainText]{[1]}) and all CTs that are classified satisfy {C1}-{C4}. 
In total there are $N_{\text{filt}}(-20, 10) = 95$ CTs eligible for classification, $\approx 10 \%$ of the detected CTs. 
Despite only classifying such a small amount of the detected CTs, the proportion of CTs classified as a particular CT-type when $T=2$ in this section are consistent with those presented in Table~2 in \hyperref[seref:MainText]{[1]} for rat K ($\approx 25 \%$ BCT, $\approx 25 \%$ BNCT, $\approx 50 \%$ NCT). 

\begin{figure}[t]
    \centering
    \includegraphics[width=0.9\textwidth]{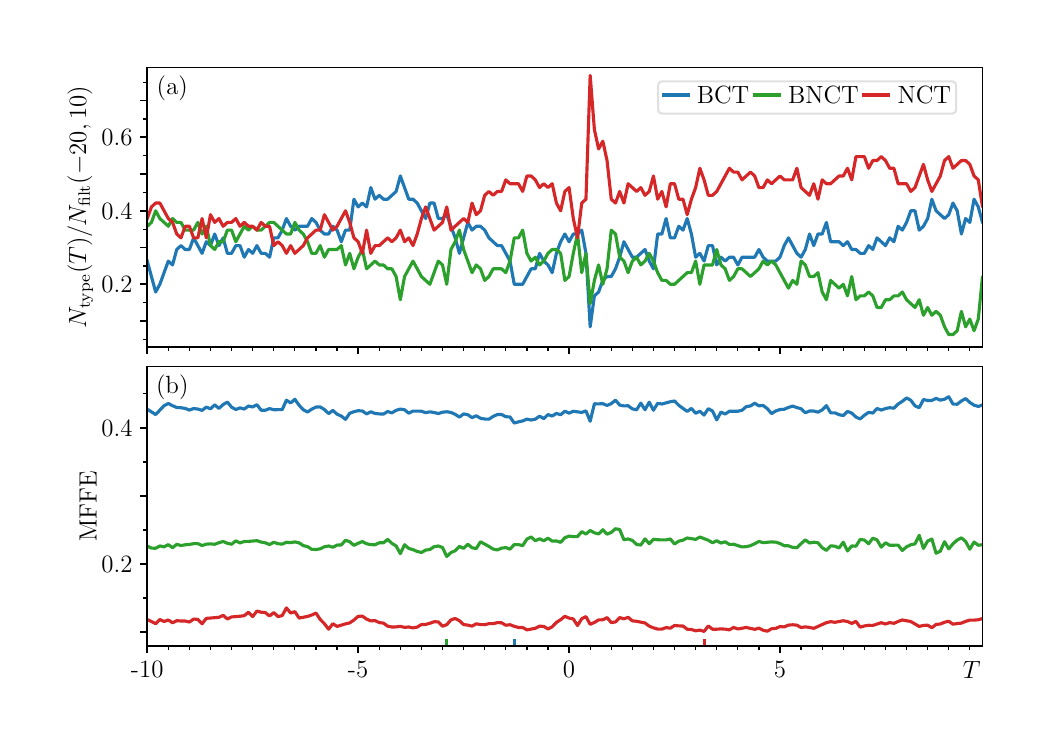}
    \caption{
    Classification of seizure generation mechanisms and MFFE (mean feature fit error) versus $T$. (a) shows $N_{\text{type}}(T)/N_{\text{filt}}(-20, 10)$ versus $T$ for rat K, i.e., proportion of filtered CTs that are classified as ({blue}) BCT, ({green}) BNCT, or ({red}) NCT, for a given $T$ with $T^{-}=-20$ and $T^{+}=10$. (b) shows the MFFE computed for each CT-type versus $T$, inward coloured ticks indicate the smallest MFFE for the corresponding CT-type.}
\label{fig:RatK_NoOfCTs_vs_T_}
\end{figure}

\textbf{SVM classifications: } 
Figure~\ref{fig:RatK_NoOfCTs_vs_T_}~(a) shows how $N_{\text{type}}(T)/N_{\text{filt}}(-20, 10)$ varies for different values of $T \in \left[ -10, 10 \right]$, this tells us what proportion of the filtered CTs are classified as a particular CT-type; BCT in {blue}, BNCT in {green}, and NCT in {red}. 

\begin{figure}[t]
    \centering
    \includegraphics[width=0.9\textwidth]{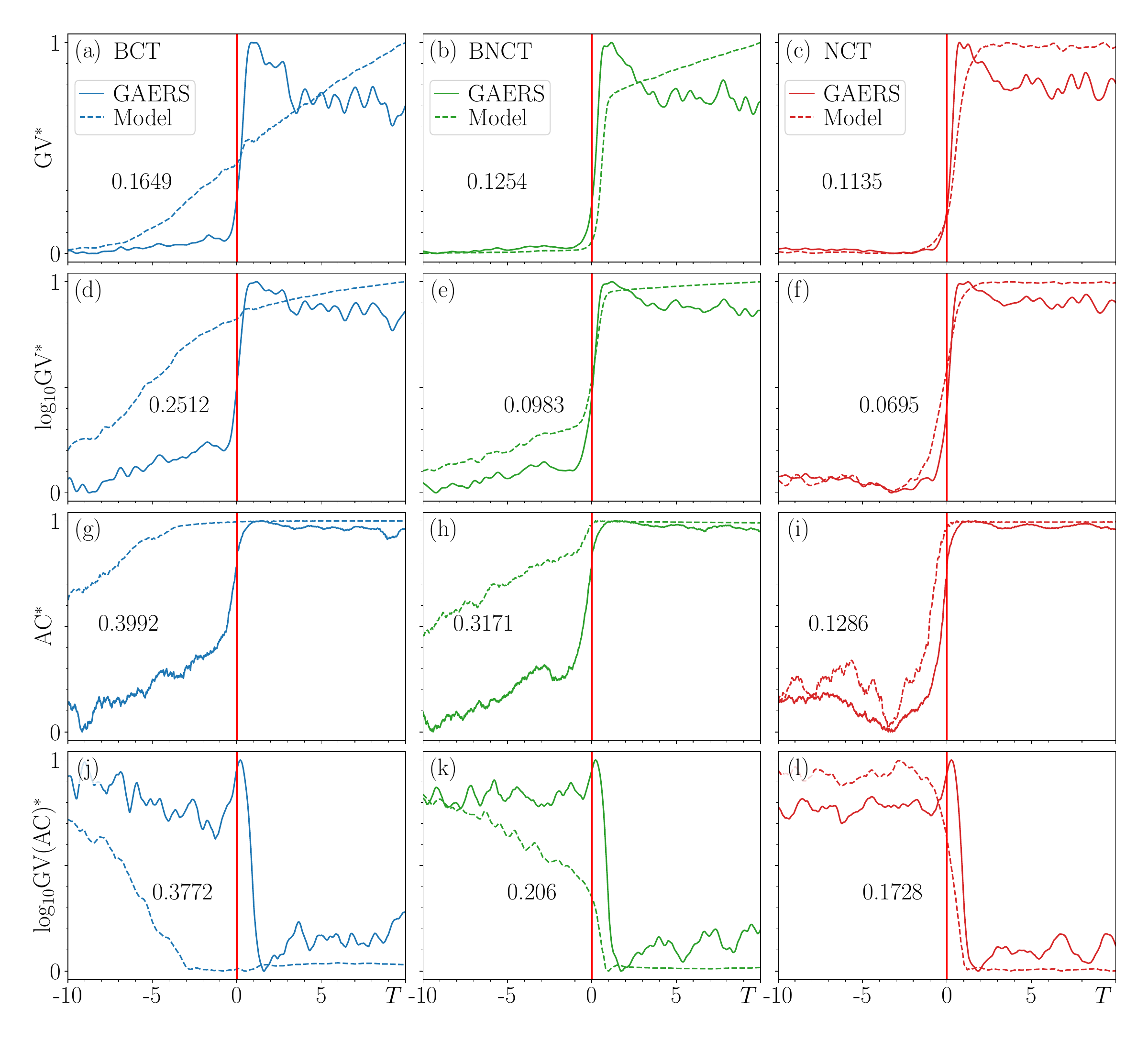}
    \caption{Comparing normalised values of TSPs of (dashed curves) the model to the corresponding quantities obtained from the (solid curves) voltage recordings that were classified by the SVM for $T=-1.3$ (the value of $T$ where the best fit of the BCTs occur, shown in Fig.~\ref{fig:RatK_NoOfCTs_vs_T_}~(b)). Each curve is coloured according to the class specified in the top left corner of panels (a)-(c). RMSE between the curves is quoted beside each pair of curves.}
\label{fig:ModelFit_SuperC_smallest_error_Values_}
\end{figure}

\begin{figure}[t]
    \centering
    \includegraphics[width=0.9\textwidth]{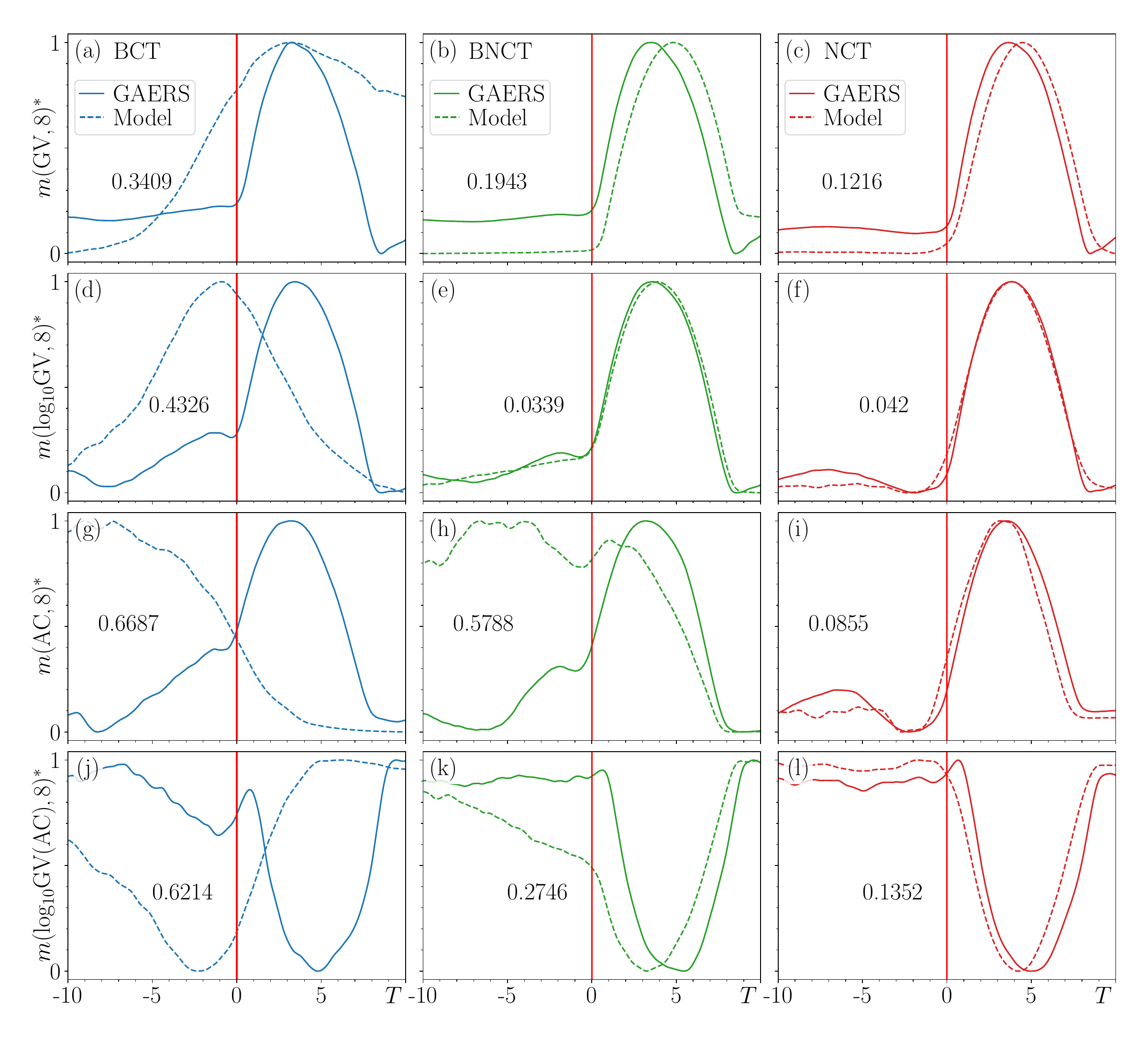}
    \caption{Comparing normalised values of $m$(TSPs, $t_m$) with $t_{m}=8$ of the (dashed curves) model to the corresponding quantities obtained from (solid curves) the voltage recordings that were classified by the SVM for $T=-1.3$ (the value of $T$ where the best fit of the BCTs occur, shown in Fig.~\ref{fig:RatK_NoOfCTs_vs_T_}~(b)). Each curve is coloured according to the class specified in the top left corner of panels (a)-(c). RMSE between the curves is quoted beside each pair of curves.}
\label{fig:ModelFit_SuperC_smallest_error_Slopes_}
\end{figure}

Figure~\ref{fig:RatK_NoOfCTs_vs_T_}~(a) shows that long before the CTs, i.e., for $T<-5$, $N_{\text{type}}(T)/N_{\text{filt}}(-20, 10)$ is relatively similar for each CT-type. 
For $T \in \left[-5, 0\right)$, Fig.~\ref{fig:RatK_NoOfCTs_vs_T_}~(a) shows that $N_{\text{type}}(T)/N_{\text{filt}}(-20, 10)$ varies significantly and there are times where different CT types appear to be the dominant seizure generation mechanism. 
However, just after the CT takes place, specifically, for $T \in \left[0, 1\right]$, there is a major increase in the proportion of CTs classified as NCT, balanced by a significant decrease in the proportion of CTs classified as either BCT or BNCT. 
For $T\in \left(1, 5 \right]$ we note that $N_{\text{type}}(T)/N_{\text{filt}}(-20, 10)$ remains relatively constant, thus providing a range of $T$ values where more reliable classifications can be obtained. 
Similar to the results presented in Fig.~7 in \hyperref[seref:MainText]{[1]}, we see that for $T=2$ in Fig.~\ref{fig:RatK_NoOfCTs_vs_T_}~(a), NCTs are also the dominant seizure generation mechanism for this portion of rat K's CTs. 
For $T>5$, the proportion of CTs classified as NCT remains relatively close to $0.5$, BCT starts to increase and is balanced by a near equal decrease in BNCT. Figure~4 in \hyperref[seref:MainText]{[1]} shows this happens because the TSPs and $m$(TSPs, $t_{m}$) of each CT-type are much less distinguishable from each other for $T>5$. Thus, the classification are less reliable for $T>5$.  
\\
\textbf{Mean feature fit error: } 
While our main results are computed using an SVM that achieves perfect accuracy when classifying the model's CTs, we are not able to quantify the accuracy of the SVM when used to classify CTs in the voltage recordings in the same way since their generation mechanisms are unknown. 
Instead, we take the following steps to compare the classifications of CTs in the voltage recordings to the corresponding CTs in the model's output in terms of the TSPs and $m$(TSPs, $t_m$):
\begin{enumerate}
    \item After classifying CTs in the voltage recordings, organise their features (i.e., the TSPs and $m(\text{TSPs},~t_{m})$) into groups based on the classifications, i.e., groups of TSPs and $m(\text{TSPs},~t_{m})$ that correspond to BCTs, BNCTs, and NCTs.
    \item Compute the mean value of each feature in the different groups for a given $T$, e.g., compute the mean GV at a given $T$ for the different groups of CTs classified as BCTs, BNCTs, and NCTs and repeat for each feature. 
    \item Compare each mean to the equivalent quantities computed from the model's output (first and third columns of Fig.~4 in \hyperref[seref:MainText]{[1]}) by:
    \begin{enumerate}
        \item[(i)] Performing a min-max normalisation of the means of each feature.
        \item[(ii)] Computing the root mean square error (RMSE) between corresponding normalised curves obtained from (i), see Figs.~\ref{fig:ModelFit_SuperC_smallest_error_Values_} and \ref{fig:ModelFit_SuperC_smallest_error_Slopes_} for an example.
        \item[(iii)] Computing the `mean feature fit error' (MFFE): the mean of the RMSEs over all features for a given $T$.
    \end{enumerate}
\end{enumerate}
In Fig.~\ref{fig:RatK_NoOfCTs_vs_T_}~(b) we plot the MFFE versus $T$ having applied steps 1-3 listed above to the data used to generate Fig.~\ref{fig:RatK_NoOfCTs_vs_T_}~(a). 
Figure~\ref{fig:RatK_NoOfCTs_vs_T_}~(b) shows there is a relatively strong agreement between NCTs in the model and CTs in the voltage recordings that are classified as NCT. Figure~\ref{fig:RatK_NoOfCTs_vs_T_}~(b) also shows there is an apparent weaker agreement between CTs that are classified as BNCT and an even weaker agreement for BCT.

To more closely investigate how well the features from the classified CTs in the voltage recordings mimic their counterparts in the model's output, we choose $T=-1.3$ as a point to visually compare the normalised means of the features from the model and the voltage recordings in Figs.~\ref{fig:ModelFit_SuperC_smallest_error_Values_} and \ref{fig:ModelFit_SuperC_smallest_error_Slopes_}. This time corresponds to when the MFFE is nearest to $0$ for the BCT case.


Figures~\ref{fig:ModelFit_SuperC_smallest_error_Values_} and \ref{fig:ModelFit_SuperC_smallest_error_Slopes_} show there are several qualitative similarities between the classified CTs in the voltage recordings and the corresponding CTs generated by the model. For instance, the first column in both Figs.~\ref{fig:ModelFit_SuperC_smallest_error_Values_} and \ref{fig:ModelFit_SuperC_smallest_error_Slopes_} show that CTs classified as BCTs exhibit strong signatures of critical slowing down (CSD), i.e., increases in GV and AC before the CT, thus indicating that the rat's brain has drifted towards and then past a bifurcation point that generates the seizure.

Figures~\ref{fig:ModelFit_SuperC_smallest_error_Values_} and \ref{fig:ModelFit_SuperC_smallest_error_Slopes_} also show that, in comparison to the BCT case, there are significantly smaller RMSEs between the curves in the BNCT and NCT cases.
For instance, panels (e) and (f) in Figs.~\ref{fig:ModelFit_SuperC_smallest_error_Values_} and \ref{fig:ModelFit_SuperC_smallest_error_Slopes_} illustrate the strong agreement between the model and the voltage recordings in terms of the normalised log$_{10}$GV and $m(\text{log$_{10}$GV}, 8)$, further indicated by the small value of the RMSE quoted beside these curves. 
However, the same cannot be said for all features. For instance, panel (h) in Figs.~\ref{fig:ModelFit_SuperC_smallest_error_Values_} and \ref{fig:ModelFit_SuperC_smallest_error_Slopes_} show much poorer agreement between BNCTs in terms of the normalised AC and $m(\text{AC}, 8)$ curves. 
The large RMSE quoted beside these curves provides the greatest contribution to the MFFE shown in Fig.~\ref{fig:RatK_NoOfCTs_vs_T_}~(b). 
The same issue arises for the BCT case as panel (g) in Figs.~\ref{fig:ModelFit_SuperC_smallest_error_Values_} and \ref{fig:ModelFit_SuperC_smallest_error_Slopes_} show the RMSE for the normalised AC and $m(\text{AC}, 8)$ curves is also significantly large. 
Importantly in both BCT and BNCT cases is that there is a clear increase in AC before the CT, however, it is the rate of increase that differs from the model's. 
Furthermore, panel (i) in Figs.~\ref{fig:ModelFit_SuperC_smallest_error_Values_} and \ref{fig:ModelFit_SuperC_smallest_error_Slopes_} show no gradual increase in AC before NCTs in the voltage recordings, which is consistent with the corresponding CTs generated by the model.

\begin{figure}[t]
\centering
\includegraphics[width=0.85\textwidth]{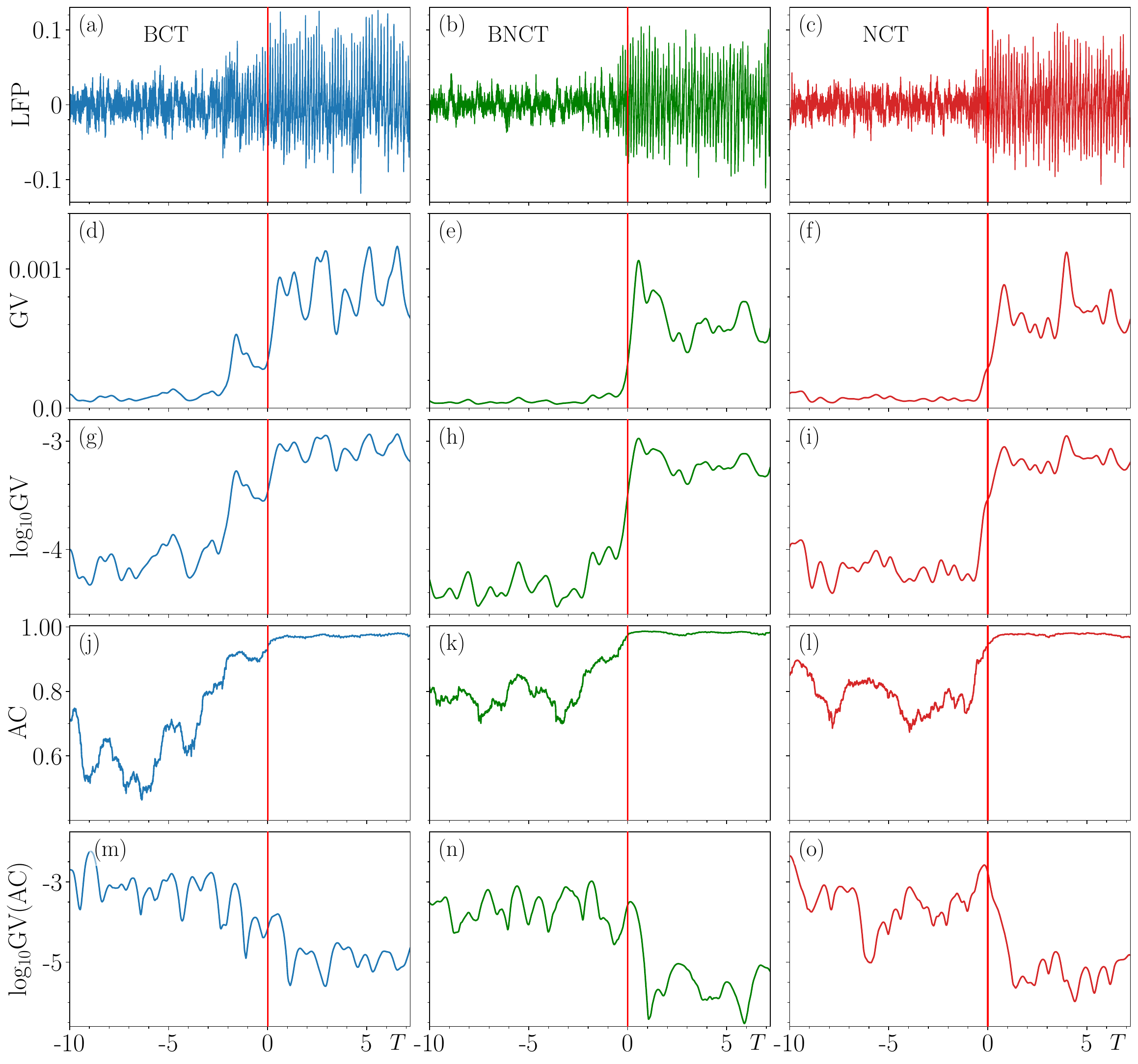}
\caption{TSPs of voltage recordings for rat K near an example of a (left-hand column) BCT, (middle column) BNCT, and (right-hand column) NCT as classified by the same model-trained SVM used to generate Figs.~\ref{fig:ModelFit_SuperC_smallest_error_Values_} and \ref{fig:ModelFit_SuperC_smallest_error_Slopes_}. For each CT-type, the top row shows the LFP versus $T=t-t_{1}$ with $t_{1}$ detected by the algorithm. The remaining rows show how each of the TSPs, specified in the vertical axes labels, behave for the same values of $T$. Vertical red line in each panel emphasises when $T=0$.}
\label{fig:EWSvalues_voltages}
\end{figure}

Figure~\ref{fig:EWSvalues_voltages} shows the TSPs for a single example of each CT-type in the voltage recordings of rat K which were classified with the same SVM as in Figs.~\ref{fig:ModelFit_SuperC_smallest_error_Values_} and \ref{fig:ModelFit_SuperC_smallest_error_Slopes_}. In these examples, one can see a gradual increase in variance and autocorrelation prior to the BCT and BNCT, and no such increase in variance and autocorrelation prior to the NCT. When compared to Fig.~\ref{fig:EWSvalues}, one can observe in greater detail how TSPs corresponding to different CT-types in the voltage recordings align with those of the model's output.

\section*{S12: Applying filtering conditions to detected CTs \label{si:ApplyingFilteringConditions}}

\begin{figure}[t]
    \centering
    \includegraphics[width=0.9\textwidth]{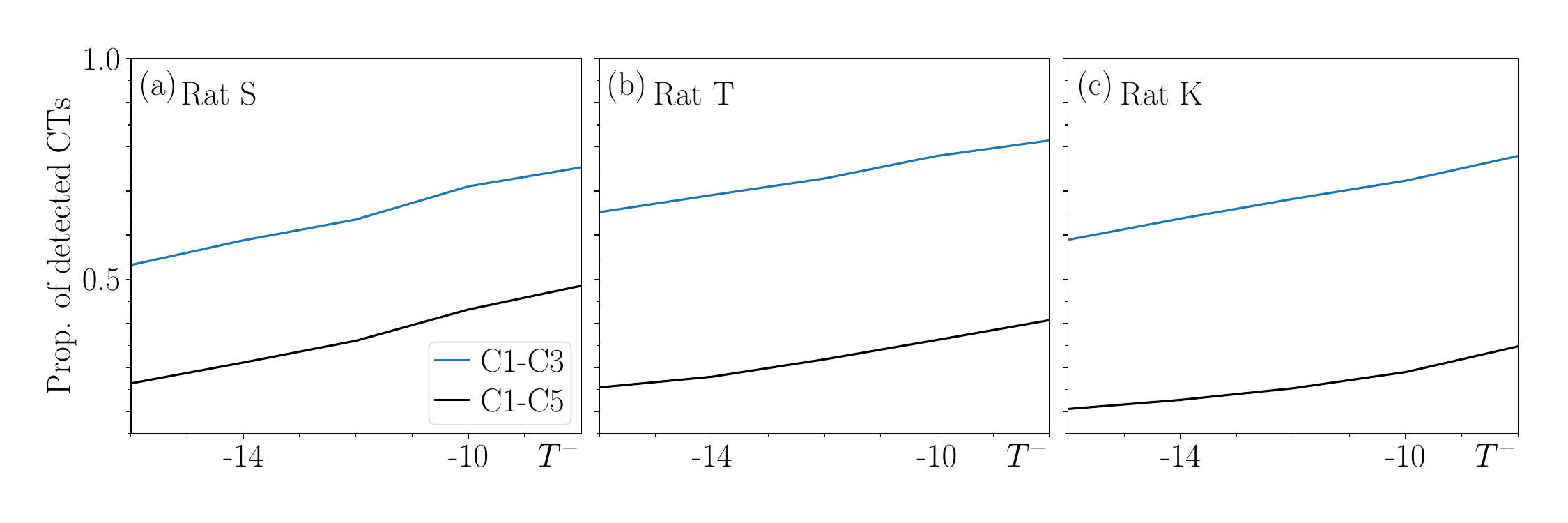}
    \caption{Complimentary figure to Figs.~7~(d)-(f) in \hyperref[seref:MainText]{[1]} showing the influence of {C4} when filtering the detected CTs of (a) Rat S, (b) Rat T, and (c) Rat K. Each panel shows the proportion of detected CTs that satisfy conditions (blue) {C1}-{C3} and (black) {C1}-{C4}, i.e., comparing the proportion of filtered CTs to detected CTs with artefacts, exaggerated spikings, and almost-occurring CTs happening nearby the CT.}
\label{fig:Change_in_Nfilt_for_diff_C}
\end{figure}

In this section we provide a more detailed account on how the conditions {C1}-{C4} contribute to the information presented in Figs.~7~(d)-(f) in \hyperref[seref:MainText]{[1]}.  
More specifically, we compare the proportion of detected CTs that satisfy {C1}-{C3} to those satisfying {C1}-{C4}. 
\\
In Fig.~\ref{fig:Change_in_Nfilt_for_diff_C} we plot the proportion of detected CTs versus $T^{-}$ for each rat for detected CTs that satisfy (i) {C1}-{C3} (in blue) and (ii) {C1}-{C4} (in black). 
Figure.~\ref{fig:Change_in_Nfilt_for_diff_C} shows that rat K is the most affected by {C4}, followed closely by rat T. Rat S is least affected. 
\\
From further calculations we found that by changing $T^{-}=$ from $-8$ to $-\tau_{\text{NS}}=-3$, the proportion of CTs that satisfy {C1}-{C3} increased from $75-85\%$ to $99.7-100\%$ across each rat. 
By imposing {C4} we found that these proportions of detected CTs fell to $~80\%$ for rat S and $~66\%$ for rats T and K.

\section*{References}

\textbf{[1]. \label{seref:MainText}} Flynn, A. \textit{et al.} “Classifying absence seizure generation mechanisms: A critical transitions framework”. (See above).\\
\textbf{[2]. \label{seref:AlgorithmROC}} Flynn, A. \textit{et al.} “Detecting seizure onset and offset times using human intelligence: A critical-transitions-based approach
”. \textit{arXiv preprint arXiv:2607.27105} (2026).


\section*{}


\begin{thebibliography}{10}
\urlstyle{rm}
\expandafter\ifx\csname url\endcsname\relax
  \def\url#1{\texttt{#1}}\fi
\expandafter\ifx\csname urlprefix\endcsname\relax\def\urlprefix{URL }\fi
\expandafter\ifx\csname doiprefix\endcsname\relax\def\doiprefix{DOI: }\fi
\providecommand{\bibinfo}[2]{#2}
\providecommand{\eprint}[2][]{\url{#2}}

\bibitem{beniczky2025_LeagueAgainstEpil_def}
\bibinfo{author}{Beniczky, S.} \emph{et~al.}
\newblock \bibinfo{journal}{\bibinfo{title}{Updated classification of epileptic seizures: {P}osition paper of the {I}nternational {L}eague {A}gainst {E}pilepsy}}.
\newblock {\emph{\JournalTitle{Epilepsia}}} \textbf{\bibinfo{volume}{66}}, \bibinfo{pages}{1804--1823} (\bibinfo{year}{2025}).

\bibitem{kaye2025_seizimpact}
\bibinfo{author}{Kaye, D.} \emph{et~al.}
\newblock \bibinfo{title}{Impact of {P}rolonged {S}eizures on {P}atients’ and {C}aregivers’ {Q}uality of {L}ife {(P1-9.014)}}.
\newblock In \emph{\bibinfo{booktitle}{Neurology}}, vol. \bibinfo{volume}{104}, \bibinfo{pages}{2708} (\bibinfo{organization}{Lippincott Williams \& Wilkins Hagerstown, MD}, \bibinfo{year}{2025}).

\bibitem{gonzalez2019_WhySeizHappen}
\bibinfo{author}{Gonz{\'a}lez, O.~C.}, \bibinfo{author}{Krishnan, G.~P.}, \bibinfo{author}{Timofeev, I.} \& \bibinfo{author}{Bazhenov, M.}
\newblock \bibinfo{journal}{\bibinfo{title}{Ionic and synaptic mechanisms of seizure generation and epileptogenesis}}.
\newblock {\emph{\JournalTitle{Neurobiol. {D}is}}} \textbf{\bibinfo{volume}{130}}, \bibinfo{pages}{104485} (\bibinfo{year}{2019}).

\bibitem{ashwinwieczorek2017CTdef}
\bibinfo{author}{Ashwin, P.}, \bibinfo{author}{Perryman, C.} \& \bibinfo{author}{Wieczorek, S.}
\newblock \bibinfo{journal}{\bibinfo{title}{Parameter shifts for nonautonomous systems in low dimension: bifurcation-and rate-induced tipping}}.
\newblock {\emph{\JournalTitle{Nonlinearity}}} \textbf{\bibinfo{volume}{30}}, \bibinfo{pages}{2185} (\bibinfo{year}{2017}).

\bibitem{lenton2008tippingpolicymaking}
\bibinfo{author}{Lenton, T.~M.} \emph{et~al.}
\newblock \bibinfo{journal}{\bibinfo{title}{Tipping elements in the {E}arth's climate system}}.
\newblock {\emph{\JournalTitle{Proceedings of the National Academy of Sciences}}} \textbf{\bibinfo{volume}{105}}, \bibinfo{pages}{1786--1793} (\bibinfo{year}{2008}).

\bibitem{lenton2020tippingpositivechange}
\bibinfo{author}{Lenton, T.~M.}
\newblock \bibinfo{journal}{\bibinfo{title}{Tipping positive change}}.
\newblock {\emph{\JournalTitle{Philosophical {T}rans. of the {R}. {S}ociety {B}}}} \textbf{\bibinfo{volume}{375}}, \bibinfo{pages}{20190123} (\bibinfo{year}{2020}).

\bibitem{armstrong2022exceeding}
\bibinfo{author}{Armstrong~McKay, D.~I.} \emph{et~al.}
\newblock \bibinfo{journal}{\bibinfo{title}{Exceeding {1.5 C} global warming could trigger multiple climate tipping points}}.
\newblock {\emph{\JournalTitle{Science}}} \textbf{\bibinfo{volume}{377}}, \bibinfo{pages}{eabn7950} (\bibinfo{year}{2022}).

\bibitem{lenton2025global}
\bibinfo{author}{Lenton, T.~M.} \emph{et~al.}
\newblock \bibinfo{title}{Global tipping points report 2025} (\bibinfo{year}{2025}).
\newblock \bibinfo{note}{University of {E}xeter}.

\bibitem{scheffer2009early}
\bibinfo{author}{Scheffer, M.} \emph{et~al.}
\newblock \bibinfo{journal}{\bibinfo{title}{Early-warning signals for critical transitions}}.
\newblock {\emph{\JournalTitle{Nature}}} \textbf{\bibinfo{volume}{461}}, \bibinfo{pages}{53} (\bibinfo{year}{2009}).

\bibitem{mcsharry2003prediction}
\bibinfo{author}{McSharry, P.~E.}, \bibinfo{author}{Smith, L.~A.} \& \bibinfo{author}{Tarassenko, L.}
\newblock \bibinfo{journal}{\bibinfo{title}{Prediction of epileptic seizures: are nonlinear methods relevant?}}
\newblock {\emph{\JournalTitle{Nature {M}ed.}}} \textbf{\bibinfo{volume}{9}}, \bibinfo{pages}{241--242} (\bibinfo{year}{2003}).

\bibitem{meiselkuehn2012scaling}
\bibinfo{author}{Meisel, C.} \& \bibinfo{author}{Kuehn, C.}
\newblock \bibinfo{journal}{\bibinfo{title}{Scaling effects and spatio-temporal multilevel dynamics in epileptic seizures}}.
\newblock {\emph{\JournalTitle{PLoS One}}} \textbf{\bibinfo{volume}{7}}, \bibinfo{pages}{e30371} (\bibinfo{year}{2012}).

\bibitem{milanowski2016seizures}
\bibinfo{author}{Milanowski, P.} \& \bibinfo{author}{Suffczynski, P.}
\newblock \bibinfo{journal}{\bibinfo{title}{Seizures start without common signatures of critical transition}}.
\newblock {\emph{\JournalTitle{Int. {J}. of {N}eural {S}yst.}}} \textbf{\bibinfo{volume}{26}}, \bibinfo{pages}{1650053} (\bibinfo{year}{2016}).

\bibitem{Jirsa_2014_NatureOfSeizDyn}
\bibinfo{author}{Jirsa, V.~K.}, \bibinfo{author}{Stacey, W.~C.}, \bibinfo{author}{Quilichini, P.~P.}, \bibinfo{author}{Ivanov, A.~I.} \& \bibinfo{author}{Bernard, C.}
\newblock \bibinfo{journal}{\bibinfo{title}{On the nature of seizure dynamics}}.
\newblock {\emph{\JournalTitle{Brain}}} \textbf{\bibinfo{volume}{137}}, \bibinfo{pages}{2210--2230} (\bibinfo{year}{2014}).

\bibitem{wilkat2019noCSDseizure}
\bibinfo{author}{Wilkat, T.}, \bibinfo{author}{Rings, T.} \& \bibinfo{author}{Lehnertz, K.}
\newblock \bibinfo{journal}{\bibinfo{title}{No evidence for critical slowing down prior to human epileptic seizures}}.
\newblock {\emph{\JournalTitle{Chaos}}} \textbf{\bibinfo{volume}{29}}, \bibinfo{pages}{091104} (\bibinfo{year}{2019}).

\bibitem{maturana2020yesCSDseizure}
\bibinfo{author}{Maturana, M.~I.} \emph{et~al.}
\newblock \bibinfo{journal}{\bibinfo{title}{Critical slowing down as a biomarker for seizure susceptibility}}.
\newblock {\emph{\JournalTitle{Nat. {C}ommun.}}} \textbf{\bibinfo{volume}{11}}, \bibinfo{pages}{2172} (\bibinfo{year}{2020}).

\bibitem{DaSilva03EpilepsyDynDis}
\bibinfo{author}{da~Silva, F.~L.} \emph{et~al.}
\newblock \bibinfo{journal}{\bibinfo{title}{Epilepsies as {D}ynamical {D}iseases of {B}rain {S}ystems: {B}asic {M}odels of the {T}ransition {B}etween {N}ormal and {E}pileptic {A}ctivity}}.
\newblock {\emph{\JournalTitle{Epilepsia}}} \textbf{\bibinfo{volume}{44}}, \bibinfo{pages}{72--83} (\bibinfo{year}{2003}).

\bibitem{Dakos_2024_ESD}
\bibinfo{author}{Dakos, V.} \emph{et~al.}
\newblock \bibinfo{journal}{\bibinfo{title}{Tipping point detection and early warnings in climate, ecological, and human systems}}.
\newblock {\emph{\JournalTitle{Earth System Dynamics}}} \textbf{\bibinfo{volume}{15}}, \bibinfo{pages}{1117--1135} (\bibinfo{year}{2024}).

\bibitem{mccafferty2018cortical}
\bibinfo{author}{McCafferty, C.} \emph{et~al.}
\newblock \bibinfo{journal}{\bibinfo{title}{Cortical drive and thalamic feed-forward inhibition control thalamic output synchrony during absence seizures}}.
\newblock {\emph{\JournalTitle{Nature neuroscience}}} \textbf{\bibinfo{volume}{21}}, \bibinfo{pages}{744--756} (\bibinfo{year}{2018}).

\bibitem{crunelli2020roots_of_epil}
\bibinfo{author}{Crunelli, V.} \emph{et~al.}
\newblock \bibinfo{journal}{\bibinfo{title}{Clinical and experimental insight into pathophysiology, comorbidity and therapy of absence seizures}}.
\newblock {\emph{\JournalTitle{Brain}}} \textbf{\bibinfo{volume}{143}}, \bibinfo{pages}{2341--2368} (\bibinfo{year}{2020}).

\bibitem{zhou2024generalized_vs_focal}
\bibinfo{author}{Zhou, Z.} \emph{et~al.}
\newblock \bibinfo{journal}{\bibinfo{title}{A generalized seizure type: {M}yoclonic-to-tonic seizure}}.
\newblock {\emph{\JournalTitle{Clinical {N}europhysiology}}} \textbf{\bibinfo{volume}{164}}, \bibinfo{pages}{24--29} (\bibinfo{year}{2024}).

\bibitem{depaulis2006models}
\bibinfo{author}{Depaulis, A.}, \bibinfo{author}{van Luijtelaar, G.}, \bibinfo{author}{Pitkanen, A.}, \bibinfo{author}{Schwartzkroin, P.} \& \bibinfo{author}{Moshe, S.}
\newblock \bibinfo{title}{Models of seizures and epilepsy} (\bibinfo{year}{2006}).

\bibitem{mccafferty2025interrupt_seiz}
\bibinfo{author}{McCafferty, C.~P.}, \bibinfo{author}{Zheng, X.}, \bibinfo{author}{Tung, R.}, \bibinfo{author}{Gruenbaum, B.~F.} \& \bibinfo{author}{Blumenfeld, H.}
\newblock \bibinfo{journal}{\bibinfo{title}{Interruption of rat absence seizures by auditory stimulation}}.
\newblock {\emph{\JournalTitle{bioRxiv}}} \bibinfo{pages}{2025--10} (\bibinfo{year}{2025}).

\bibitem{DaSilva03dynamicaldisease}
\bibinfo{author}{da~Silva, F. H.~L.} \emph{et~al.}
\newblock \bibinfo{journal}{\bibinfo{title}{Dynamical diseases of brain systems: different routes to epileptic seizures}}.
\newblock {\emph{\JournalTitle{IEEE {T}rans. {B}iomedical {E}ngineering}}} \textbf{\bibinfo{volume}{50}}, \bibinfo{pages}{540--548} (\bibinfo{year}{2003}).

\bibitem{suffczynski04EpilepsyDynamics}
\bibinfo{author}{Suffczynski, P.}, \bibinfo{author}{Kalitzin, S.} \& \bibinfo{author}{da~Silva, F.~L.}
\newblock \bibinfo{journal}{\bibinfo{title}{Dynamics of non-convulsive epileptic phenomena modeled by a bistable neuronal network}}.
\newblock {\emph{\JournalTitle{Neuroscience}}} \textbf{\bibinfo{volume}{126}}, \bibinfo{pages}{467--484} (\bibinfo{year}{2004}).

\bibitem{suffczynski05EpilepticTransitions}
\bibinfo{author}{Suffczynski, P.}, \bibinfo{author}{da~Silva, F.~L.}, \bibinfo{author}{Parra, J.}, \bibinfo{author}{Velis, D.} \& \bibinfo{author}{Kalitzin, S.}
\newblock \bibinfo{journal}{\bibinfo{title}{Epileptic transitions: model predictions and experimental validation}}.
\newblock {\emph{\JournalTitle{Journal of {C}linical {N}europhysiology}}} \textbf{\bibinfo{volume}{22}}, \bibinfo{pages}{288--299} (\bibinfo{year}{2005}).

\bibitem{jungesTerry20epilepsySurgery_SubCSupC_Network}
\bibinfo{author}{Junges, L.}, \bibinfo{author}{Woldman, W.}, \bibinfo{author}{Benjamin, O.~J.} \& \bibinfo{author}{Terry, J.~R.}
\newblock \bibinfo{journal}{\bibinfo{title}{Epilepsy surgery: Evaluating robustness using dynamic network models}}.
\newblock {\emph{\JournalTitle{Chaos}}} \textbf{\bibinfo{volume}{30}}, \bibinfo{pages}{113106} (\bibinfo{year}{2020}).

\bibitem{harringtonTerry24_SubCSupC_Network}
\bibinfo{author}{Harrington, E.~G.}, \bibinfo{author}{Kissack, P.}, \bibinfo{author}{Terry, J.~R.}, \bibinfo{author}{Woldman, W.} \& \bibinfo{author}{Junges, L.}
\newblock \bibinfo{journal}{\bibinfo{title}{Treatment effects in epilepsy: a mathematical framework for understanding response over time}}.
\newblock {\emph{\JournalTitle{Frontiers in {N}etw. {P}hysiol.}}} \textbf{\bibinfo{volume}{4}}, \bibinfo{pages}{1308501} (\bibinfo{year}{2024}).

\bibitem{qinBassett24_SubSupC_Network}
\bibinfo{author}{Qin, Y.}, \bibinfo{author}{El-Gazzar, A.}, \bibinfo{author}{Bassett, D.~S.}, \bibinfo{author}{Pasqualetti, F.} \& \bibinfo{author}{Van~Gerven, M.}
\newblock \bibinfo{title}{Analytical characterization of epileptic dynamics in a bistable system}.
\newblock In \emph{\bibinfo{booktitle}{2024 IEEE 63rd Conference on Decision and Control (CDC)}}, \bibinfo{pages}{583--588} (\bibinfo{organization}{IEEE}, \bibinfo{year}{2024}).

\bibitem{byrne2022modelEEG}
\bibinfo{author}{Byrne, {\'A}.}, \bibinfo{author}{Ross, J.}, \bibinfo{author}{Nicks, R.} \& \bibinfo{author}{Coombes, S.}
\newblock \bibinfo{journal}{\bibinfo{title}{Mean-field models for {EEG}/{MEG}: from oscillations to waves}}.
\newblock {\emph{\JournalTitle{Brain {T}opogr.}}} \textbf{\bibinfo{volume}{35}}, \bibinfo{pages}{36--53} (\bibinfo{year}{2022}).

\bibitem{Flynn25_Supplement}
\bibinfo{author}{Flynn, A.} \emph{et~al.}
\newblock \bibinfo{title}{Supplementary material from: ``{C}lassifying absence seizure generation mechanisms: {A} critical transitions framework''}.
\newblock \bibinfo{note}{{(See below from page 27 onwards)}}.

\bibitem{Pokrovskii12systemswHysteresis}
\bibinfo{author}{Krasnosel'skii, M.~A.} \& \bibinfo{author}{Pokrovskii, A.~V.}
\newblock \emph{\bibinfo{title}{Systems with {H}ysteresis}} (\bibinfo{publisher}{Springer Science \& Business Media}, \bibinfo{year}{2012}).

\bibitem{flynn2026detecting}
\bibinfo{author}{Flynn, A.} \emph{et~al.}
\newblock \bibinfo{journal}{\bibinfo{title}{Detecting seizure onset and offset times using human intelligence: {A} critical-transitions-based approach}}.
\newblock {\emph{\JournalTitle{arXiv preprint arXiv:2607.27105}}}  (\bibinfo{year}{2026}).

\bibitem{ditlevsen2023ews}
\bibinfo{author}{Ditlevsen, P.} \& \bibinfo{author}{Ditlevsen, S.}
\newblock \bibinfo{journal}{\bibinfo{title}{Warning of a forthcoming collapse of the {A}tlantic meridional overturning circulation}}.
\newblock {\emph{\JournalTitle{Nature {C}ommun.}}} \textbf{\bibinfo{volume}{14}}, \bibinfo{pages}{4254} (\bibinfo{year}{2023}).

\bibitem{vandijkstra2024ewsAMOC}
\bibinfo{author}{Van~Westen, R.~M.}, \bibinfo{author}{Kliphuis, M.} \& \bibinfo{author}{Dijkstra, H.~A.}
\newblock \bibinfo{journal}{\bibinfo{title}{Physics-based early warning signal shows that {AMOC} is on tipping course}}.
\newblock {\emph{\JournalTitle{Science {A}dv.}}} \textbf{\bibinfo{volume}{10}}, \bibinfo{pages}{1189} (\bibinfo{year}{2024}).

\bibitem{ashwin2025_ews_skill}
\bibinfo{author}{Ashwin, P.}, \bibinfo{author}{Bastiaansen, R.}, \bibinfo{author}{von~der {H}eydt, A.~S.} \& \bibinfo{author}{Ritchie, P.~D.}
\newblock \bibinfo{journal}{\bibinfo{title}{Early warning skill, extrapolation and tipping for accelerating cascades}}.
\newblock {\emph{\JournalTitle{Proc. {R}oy. {S}oc. A}}} \textbf{\bibinfo{volume}{481}}, \bibinfo{pages}{20250405} (\bibinfo{year}{2025}).

\bibitem{bury2021deepEWS}
\bibinfo{author}{Bury, T.~M.} \emph{et~al.}
\newblock \bibinfo{journal}{\bibinfo{title}{Deep learning for early warning signals of tipping points}}.
\newblock {\emph{\JournalTitle{Proceedings of the {N}ational {A}cademy of {S}ciences}}} \textbf{\bibinfo{volume}{118}}, \bibinfo{pages}{e2106140118} (\bibinfo{year}{2021}).

\bibitem{mccormick2001cellular}
\bibinfo{author}{McCormick, D.~A.} \& \bibinfo{author}{Contreras, D.}
\newblock \bibinfo{journal}{\bibinfo{title}{On the cellular and network bases of epileptic seizures}}.
\newblock {\emph{\JournalTitle{Annu. {R}ev. {P}hysiol.}}} \textbf{\bibinfo{volume}{63}}, \bibinfo{pages}{815--846} (\bibinfo{year}{2001}).

\bibitem{vezzani2011role}
\bibinfo{author}{Vezzani, A.}, \bibinfo{author}{French, J.}, \bibinfo{author}{Bartfai, T.} \& \bibinfo{author}{Baram, T.~Z.}
\newblock \bibinfo{journal}{\bibinfo{title}{The role of inflammation in epilepsy}}.
\newblock {\emph{\JournalTitle{Nat. Rev. {N}eurol.}}} \textbf{\bibinfo{volume}{7}}, \bibinfo{pages}{31--40} (\bibinfo{year}{2011}).

\bibitem{kuhlmannlehnertz2018_seizurepred_newera}
\bibinfo{author}{Kuhlmann, L.}, \bibinfo{author}{Lehnertz, K.}, \bibinfo{author}{Richardson, M.~P.}, \bibinfo{author}{Schelter, B.} \& \bibinfo{author}{Zaveri, H.~P.}
\newblock \bibinfo{journal}{\bibinfo{title}{Seizure prediction — ready for a new era}}.
\newblock {\emph{\JournalTitle{Nature {R}eviews {N}eurology}}} \textbf{\bibinfo{volume}{14}}, \bibinfo{pages}{618--630} (\bibinfo{year}{2018}).

\bibitem{gonzalez2019ionic}
\bibinfo{author}{Gonz{\'a}lez, O.~C.}, \bibinfo{author}{Krishnan, G.~P.}, \bibinfo{author}{Timofeev, I.} \& \bibinfo{author}{Bazhenov, M.}
\newblock \bibinfo{journal}{\bibinfo{title}{Ionic and synaptic mechanisms of seizure generation and epileptogenesis}}.
\newblock {\emph{\JournalTitle{Neurobiol. {D}is.}}} \textbf{\bibinfo{volume}{130}}, \bibinfo{pages}{104485} (\bibinfo{year}{2019}).

\bibitem{jiruska2023update}
\bibinfo{author}{Jiruska, P.}, \bibinfo{author}{Freestone, D.}, \bibinfo{author}{Gnatkovsky, V.} \& \bibinfo{author}{Wang, Y.}
\newblock \bibinfo{journal}{\bibinfo{title}{An update on the seizures beget seizures theory}}.
\newblock {\emph{\JournalTitle{Epilepsia}}} \textbf{\bibinfo{volume}{64}}, \bibinfo{pages}{S13--S24} (\bibinfo{year}{2023}).

\bibitem{lehnertz2023seizcontrol}
\bibinfo{author}{Lehnertz, K.}, \bibinfo{author}{Broehl, T.} \& \bibinfo{author}{von Wrede, R.}
\newblock \bibinfo{journal}{\bibinfo{title}{Epileptic-network-based prediction and control of seizures in humans}}.
\newblock {\emph{\JournalTitle{Neurobiol. {D}is.}}} \textbf{\bibinfo{volume}{181}}, \bibinfo{pages}{106098} (\bibinfo{year}{2023}).

\bibitem{polack2007seiz_neuron}
\bibinfo{author}{Polack, P.-O.} \emph{et~al.}
\newblock \bibinfo{journal}{\bibinfo{title}{Deep layer somatosensory cortical neurons initiate spike-and-wave discharges in a genetic model of absence seizures}}.
\newblock {\emph{\JournalTitle{Journal of Neuroscience}}} \textbf{\bibinfo{volume}{27}}, \bibinfo{pages}{6590--6599} (\bibinfo{year}{2007}).

\bibitem{meeren2002GAERS}
\bibinfo{author}{Meeren, H.~K.}, \bibinfo{author}{Pijn, J. P.~M.}, \bibinfo{author}{Van~Luijtelaar, E.~L.}, \bibinfo{author}{Coenen, A.~M.} \& \bibinfo{author}{da~Silva, F. H.~L.}
\newblock \bibinfo{journal}{\bibinfo{title}{Cortical focus drives widespread corticothalamic networks during spontaneous absence seizures in rats}}.
\newblock {\emph{\JournalTitle{Journal of {N}euroscience}}} \textbf{\bibinfo{volume}{22}}, \bibinfo{pages}{1480--1495} (\bibinfo{year}{2002}).

\bibitem{jalife1979phaseresetting}
\bibinfo{author}{Jalife, J.} \& \bibinfo{author}{Antzelevitch, C.}
\newblock \bibinfo{journal}{\bibinfo{title}{Phase resetting and annihilation of pacemaker activity in cardiac tissue}}.
\newblock {\emph{\JournalTitle{Science}}} \textbf{\bibinfo{volume}{206}}, \bibinfo{pages}{695--697} (\bibinfo{year}{1979}).

\bibitem{glass1988clocks}
\bibinfo{author}{Glass, L.} \& \bibinfo{author}{Mackey, M.}
\newblock \emph{\bibinfo{title}{From {C}locks to {C}haos: {T}he {R}hythms of {L}ife}}.
\newblock Princeton paperbacks (\bibinfo{publisher}{Princeton University Press}, \bibinfo{year}{1988}).

\bibitem{winfree1980geometry}
\bibinfo{author}{Winfree, A.~T.}
\newblock \emph{\bibinfo{title}{The {G}eometry of {B}iological {T}ime}}, vol.~\bibinfo{volume}{2} (\bibinfo{publisher}{Springer}, \bibinfo{year}{1980}).

\bibitem{osorio2009seizure}
\bibinfo{author}{Osorio, I.} \& \bibinfo{author}{Frei, M.~G.}
\newblock \bibinfo{journal}{\bibinfo{title}{Seizure abatement with single {DC} pulses: {I}s phase resetting at play?}}
\newblock {\emph{\JournalTitle{International {J}ournal of {N}eural {S}ystems}}} \textbf{\bibinfo{volume}{19}}, \bibinfo{pages}{149--156} (\bibinfo{year}{2009}).

\bibitem{mccafferty2023decreased}
\bibinfo{author}{McCafferty, C.} \emph{et~al.}
\newblock \bibinfo{journal}{\bibinfo{title}{Decreased but diverse activity of cortical and thalamic neurons in consciousness-impairing rodent absence seizures}}.
\newblock {\emph{\JournalTitle{Nature {C}ommunications}}} \textbf{\bibinfo{volume}{14}}, \bibinfo{pages}{117} (\bibinfo{year}{2023}).

\bibitem{bai2010GAERS}
\bibinfo{author}{Bai, X.} \emph{et~al.}
\newblock \bibinfo{journal}{\bibinfo{title}{Dynamic time course of typical childhood absence seizures: {EEG}, behavior, and functional magnetic resonance imaging}}.
\newblock {\emph{\JournalTitle{Journal of {N}euroscience}}} \textbf{\bibinfo{volume}{30}}, \bibinfo{pages}{5884--5893} (\bibinfo{year}{2010}).

\bibitem{scalmani2026GABA}
\bibinfo{author}{Scalmani, P.} \emph{et~al.}
\newblock \bibinfo{journal}{\bibinfo{title}{Non-synaptic {GABA} release as a trigger for synchronous epileptiform network patterns}}.
\newblock {\emph{\JournalTitle{Epilepsia}}} \textbf{\bibinfo{volume}{67}}, \bibinfo{pages}{3141--3157} (\bibinfo{year}{2026}).

\bibitem{gluckman2001seizcontrol}
\bibinfo{author}{Gluckman, B.~J.}, \bibinfo{author}{Nguyen, H.}, \bibinfo{author}{Weinstein, S.~L.} \& \bibinfo{author}{Schiff, S.~J.}
\newblock \bibinfo{journal}{\bibinfo{title}{Adaptive electric field control of epileptic seizures}}.
\newblock {\emph{\JournalTitle{Journal of Neuroscience}}} \textbf{\bibinfo{volume}{21}}, \bibinfo{pages}{590--600} (\bibinfo{year}{2001}).

\bibitem{grziwotz2023_dyn_e_value}
\bibinfo{author}{Grziwotz, F.} \emph{et~al.}
\newblock \bibinfo{journal}{\bibinfo{title}{Anticipating the occurrence and type of critical transitions}}.
\newblock {\emph{\JournalTitle{Science {A}dvances}}} \textbf{\bibinfo{volume}{9}}, \bibinfo{pages}{eabq4558} (\bibinfo{year}{2023}).

\bibitem{Kimia24_HeartSeizDetect}
\bibinfo{author}{Rezaei, K.} \emph{et~al.}
\newblock \bibinfo{title}{Assessing the effectiveness of heart rate variability as a diagnostic tool for brain injuries in infants}.
\newblock In \emph{\bibinfo{booktitle}{2024 46th Annual International Conference of the IEEE Engineering in Medicine and Biology Society (EMBC)}}, \bibinfo{pages}{1--4} (\bibinfo{year}{2024}).

\bibitem{zhang2024review}
\bibinfo{author}{Zhang, X.}, \bibinfo{author}{Zhang, X.}, \bibinfo{author}{Huang, Q.} \& \bibinfo{author}{Chen, F.}
\newblock \bibinfo{journal}{\bibinfo{title}{A review of epilepsy detection and prediction methods based on {EEG} signal processing and deep learning}}.
\newblock {\emph{\JournalTitle{Frontiers in {N}euroscience}}} \textbf{\bibinfo{volume}{18}}, \bibinfo{pages}{1468967} (\bibinfo{year}{2024}).

\bibitem{Aravind25_preprint}
\bibinfo{author}{Kamaraj, A.~K.} \& \bibinfo{author}{Szuromi, M.~P.}
\newblock \bibinfo{journal}{\bibinfo{title}{Modelling dysfunction-specific interventions for seizure termination in epilepsy}}.
\newblock {\emph{\JournalTitle{NPJ systems biology and applications}}} \textbf{\bibinfo{volume}{12}}, \bibinfo{pages}{9} (\bibinfo{year}{2025}).

\bibitem{kramer2012_SeizTerm}
\bibinfo{author}{Kramer, M.~A.} \emph{et~al.}
\newblock \bibinfo{journal}{\bibinfo{title}{Human seizures self-terminate across spatial scales via a critical transition}}.
\newblock {\emph{\JournalTitle{Proceedings of the {N}ational {A}cademy of {S}ciences}}} \textbf{\bibinfo{volume}{109}}, \bibinfo{pages}{21116--21121} (\bibinfo{year}{2012}).

\bibitem{schindler2007assessing}
\bibinfo{author}{Schindler, K.}, \bibinfo{author}{Leung, H.}, \bibinfo{author}{Elger, C.~E.} \& \bibinfo{author}{Lehnertz, K.}
\newblock \bibinfo{journal}{\bibinfo{title}{Assessing seizure dynamics by analysing the correlation structure of multichannel intracranial {EEG}}}.
\newblock {\emph{\JournalTitle{Brain}}} \textbf{\bibinfo{volume}{130}}, \bibinfo{pages}{65--77} (\bibinfo{year}{2007}).

\bibitem{ritchie2023rate}
\bibinfo{author}{Ritchie, P.~D.}, \bibinfo{author}{Alkhayuon, H.}, \bibinfo{author}{Cox, P.~M.} \& \bibinfo{author}{Wieczorek, S.}
\newblock \bibinfo{journal}{\bibinfo{title}{Rate-induced tipping in natural and human systems}}.
\newblock {\emph{\JournalTitle{Earth {S}ystem {D}ynamics}}} \textbf{\bibinfo{volume}{14}}, \bibinfo{pages}{669--683} (\bibinfo{year}{2023}).

\bibitem{sujith2017thermoacoustic}
\bibinfo{author}{Etikyala, S.} \& \bibinfo{author}{Sujith, R.}
\newblock \bibinfo{journal}{\bibinfo{title}{Change of criticality in a prototypical thermoacoustic system}}.
\newblock {\emph{\JournalTitle{Chaos}}} \textbf{\bibinfo{volume}{27}} (\bibinfo{year}{2017}).

\bibitem{kurths2025aircraft}
\bibinfo{author}{Ma, J.} \emph{et~al.}
\newblock \bibinfo{journal}{\bibinfo{title}{Predicting tipping phenomenon in a conceptual airfoil structure under extreme flight environment}}.
\newblock {\emph{\JournalTitle{Journal of Sound and Vibration}}} \textbf{\bibinfo{volume}{618}}, \bibinfo{pages}{119306} (\bibinfo{year}{2025}).

\bibitem{powell2014seizure}
\bibinfo{author}{Powell, K.~L.} \emph{et~al.}
\newblock \bibinfo{journal}{\bibinfo{title}{Seizure expression, behavior, and brain morphology differences in colonies of {G}enetic {A}bsence {E}pilepsy {R}ats from {S}trasbourg}}.
\newblock {\emph{\JournalTitle{Epilepsia}}} \textbf{\bibinfo{volume}{55}}, \bibinfo{pages}{1959--1968} (\bibinfo{year}{2014}).

\bibitem{gregorvcivc2022difficult}
\bibinfo{author}{Gregor{\v{c}}i{\v{c}}, S.} \emph{et~al.}
\newblock \bibinfo{journal}{\bibinfo{title}{Difficult to treat absence seizures in children: {A} single-center retrospective study}}.
\newblock {\emph{\JournalTitle{Frontiers in {N}eurology}}} \textbf{\bibinfo{volume}{13}}, \bibinfo{pages}{958369} (\bibinfo{year}{2022}).

\bibitem{depaulis2016genetic}
\bibinfo{author}{Depaulis, A.}, \bibinfo{author}{David, O.} \& \bibinfo{author}{Charpier, S.}
\newblock \bibinfo{journal}{\bibinfo{title}{The genetic absence epilepsy rat from strasbourg as a model to decipher the neuronal and network mechanisms of generalized idiopathic epilepsies}}.
\newblock {\emph{\JournalTitle{Journal {N}eurosci. {M}ethods}}} \textbf{\bibinfo{volume}{260}}, \bibinfo{pages}{159--174} (\bibinfo{year}{2016}).

\bibitem{kuznetsov2013elements}
\bibinfo{author}{Kuznetsov, Y.~A.}
\newblock \emph{\bibinfo{title}{Elements of {A}pplied {B}ifurcation {T}heory}} (\bibinfo{publisher}{Springer Science \& Business Media}, \bibinfo{year}{2013}).

\bibitem{LinYoung2008shear}
\bibinfo{author}{Lin, K.~K.} \& \bibinfo{author}{Young, L.-S.}
\newblock \bibinfo{journal}{\bibinfo{title}{Shear-induced chaos}}.
\newblock {\emph{\JournalTitle{Nonlinearity}}} \textbf{\bibinfo{volume}{21}}, \bibinfo{pages}{899} (\bibinfo{year}{2008}).

\bibitem{Wieczorek2009shear}
\bibinfo{author}{Wieczorek, S.}
\newblock \bibinfo{journal}{\bibinfo{title}{Stochastic bifurcation in noise-driven lasers and {H}opf oscillators}}.
\newblock {\emph{\JournalTitle{Physical Review E}}} \textbf{\bibinfo{volume}{79}}, \bibinfo{pages}{036209} (\bibinfo{year}{2009}).

\bibitem{BlackbeardWieczorek2011shear}
\bibinfo{author}{Blackbeard, N.}, \bibinfo{author}{Erzgr{\"a}ber, H.} \& \bibinfo{author}{Wieczorek, S.}
\newblock \bibinfo{journal}{\bibinfo{title}{Shear-induced bifurcations and chaos in models of three coupled lasers}}.
\newblock {\emph{\JournalTitle{SIAM {J}ournal on {A}pplied {D}ynamical {S}ystems}}} \textbf{\bibinfo{volume}{10}}, \bibinfo{pages}{469--509} (\bibinfo{year}{2011}).

\bibitem{dakos2008detrend}
\bibinfo{author}{Dakos, V.} \emph{et~al.}
\newblock \bibinfo{journal}{\bibinfo{title}{Slowing down as an early warning signal for abrupt climate change}}.
\newblock {\emph{\JournalTitle{Proceedings of the {N}ational {A}cademy of {S}ciences}}} \textbf{\bibinfo{volume}{105}}, \bibinfo{pages}{14308--14312} (\bibinfo{year}{2008}).

\bibitem{lenton2012detrend}
\bibinfo{author}{Lenton, T.~M.}, \bibinfo{author}{Livina, V.}, \bibinfo{author}{Dakos, V.}, \bibinfo{author}{van Nes, E.~H.} \& \bibinfo{author}{Scheffer, M.}
\newblock \bibinfo{journal}{\bibinfo{title}{Early warning of climate tipping points from critical slowing down: comparing methods to improve robustness}}.
\newblock {\emph{\JournalTitle{Philos. {T}rans. {R}. {S}oc. A: Mathematical, Physical and Engineering Sciences}}} \textbf{\bibinfo{volume}{370}}, \bibinfo{pages}{1185--1204} (\bibinfo{year}{2012}).

\end{thebibliography}
\end{document}